\documentclass{amsart}
\usepackage[utf8]{inputenc}
\usepackage{amsmath, amssymb, amsfonts, amsthm}
\usepackage{hyperref}
\usepackage{enumerate}
\usepackage{bbold}
\usepackage{color}
\usepackage{graphicx}
\usepackage{fullpage}
\usepackage{tikz}
\usepackage[boxsize = .3em]{ytableau}
\usepackage{blkarray, bigstrut} 
\usepackage[shortlabels]{enumitem}

\usepackage{tikz-cd}
\usepackage{hyperref}
\usepackage{lipsum} 
\usepackage{comment}
\usepackage{bbm}

\allowdisplaybreaks
\DeclareMathOperator{\Rep}{\textrm{Rep}}

\DeclareMathOperator{\Z}{\mathbb{Z}}

\DeclareMathOperator{\End}{\textrm{End}}

\DeclareMathOperator{\id}{\textrm{id}}

\DeclareMathOperator{\Hom}{\textrm{Hom}}

\DeclareMathOperator{\im}{\textrm{Im}}

\newcommand{\arxiv}[1]{\href{http://arxiv.org/abs/#1}{\tt arXiv:\nolinkurl{#1}}}

\newcommand{\att}[2]{\raisebox{-.5\height}{ \includegraphics[scale = #2]{diagrams/#1.pdf}}}
\usepackage{array}   
\newcolumntype{L}{>{$}l<{$}}

\def\End{{\rm End}}

\def\al{\alpha}

\def\Gn0{\bar G_n[0]}

\def\Z{{\mathbb Z}}

\def\E{{\mathcal E}}

\def\Gno{G_n[0]}

\def\tr{\tilde{\mathrm{tr}}}

\theoremstyle{plain}
\newtheorem{thm}{Theorem}[section]
\newtheorem{lem}[thm]{Lemma}
\newtheorem{prop}[thm]{Proposition}
\newtheorem{cor}[thm]{Corollary}

\newtheorem{con}[thm]{Conjecture}

\theoremstyle{definition}
\newtheorem{defn}[thm]{Definition}
\newtheorem{question}[thm]{Question}
\newtheorem{warning}[thm]{Warning}
\newtheorem{ex}[thm]{Example}
\newtheorem{remark}[thm]{Remark}

\newcommand{\edit}[1]{#1}

\title{Interpolation categories for Conformal Embeddings}
\author{Cain Edie-Michell and Noah Snyder}
\address{Cain Edie-Michell\\
University of New Hampshire\\
Kingsbury Hall, W383 \\
33 Academic Way \\
Durham, NH 03824}
\email{cain.edie-michell@unh.edu}
\address{Noah Snyder\\
Indiana University\\
Rawles Hall \\
831 East 3rd St.\\
Bloomington, IN 47405-7106}
\email{nsnyder1@\edit{iu}.edu}
\date{}

\subjclass{18M15, 18M20, 17B37}

\begin{document}

\begin{abstract}
In this paper we give a diagrammatic description of the categories of modules coming from the conformal embeddings $\mathcal{V}(\mathfrak{sl}_N,N) \subset  \mathcal{V}(\mathfrak{so}_{N^2-1},1)$. A small variant on this construction (morally corresponding to a conformal embedding of $\mathfrak{gl}_N$ level $N$ into $\mathfrak{o}_{N^2-1}$ level $1$) has uniform generators and relations which are rational functions in $q = e^{2 \pi i/4N}$, which allows us to construct a new continuous family of tensor categories at non-integer level which interpolate between these categories. This is the second example of such an interpolation category for families of conformal embeddings after Zhengwei Liu's interpolation categories $\mathcal{V}(\mathfrak{sl}_N, N\pm 2) \subset \mathcal{V}(\mathfrak{sl}_{N(N\pm 1)/2},1)$ which he constructed using his classification \edit{of} Yang-Baxter planar algebras. Our approach is different from Liu's, we build a two-color skein theory, with one strand coming from $X$ the image of defining representation of $\mathfrak{sl}_N$ and the other strand coming from an invertible object $g$ in the category of local modules, and a trivalent vertex coming from a map $X \otimes X^* \rightarrow g$. We anticipate small variations on our approach will yield interpolation categories for every infinite discrete family of conformal embeddings.
\end{abstract}
\maketitle

\section{Introduction}

A major theme in the study of tensor categories is to take a discrete family of categories, and construct a continuous family of categories which ``interpolates'' between them \cite{MR654325, MR1106898, Deligne, MR2349713,MR3200430,MR3200430}. The famous examples of this construction are Deligne's interpolation categories $GL_t$, $O_t$, and $S_t$. Here interpolation means that if you specialize $t$ to a positive integer in say Deligne's $O_t$ and semisimplify, you recover $\Rep(O(n))$. In the case of $O_t$, the endomorphism algebras of the interpolated ``vector representation'' are the Brauer algebras \cite{MR1503378}. These interpolation categories capture certain stability in the representation theory of the discrete parameter (for example, the fusion rules stabilize as $n$ grows, and the dimensions of representations are given by polynomials \edit{in} $n$).

One important source of tensor categories are the ``quantum subgroups'' which are categories of $A$-modules for $A$ an \textit{\'etale} algebra object in the semisimplified category of representations of a quantum group at a root of unity. The main source of quantum subgroups comes from conformal field theory, where a conformal embedding of VOAs $\mathcal{V}(\mathfrak{g},k) \subset  \mathcal{V}(\mathfrak{h},1)$ yields an \'etale algebra object $A$ in the braided tensor category of modules for $\mathcal{V}(\mathfrak{g},k)$ \cite{Xu, Kril,HKL} (which is closely related to \edit{the} semisimplified quantum group category of $\mathfrak{g}$ at root of 1). Among the conformal embeddings (listed on \cite[p. 34]{LagrangeUkraine}) are several discrete families yielding non-pointed \textit{\'etale} algebra objects:

\begin{itemize}
    \item $\mathcal{V}(\mathfrak{sl}_N, N \pm 2) \subset \mathcal{V}(\mathfrak{sl}_{N(N\pm 1)/2},1)$
    \item $\mathcal{V}(\mathfrak{sl}_N,N) \subset  \mathcal{V}(\mathfrak{so}_{N^2-1},1)$
    \item $\mathcal{V}(\mathfrak{so}_N,N\pm 2) \subset  \mathcal{V}(\mathfrak{so}_{(N+1\pm 1)(N-1)/2},1)$
    \item $\mathcal{V}(\mathfrak{sp}_{2N},N\pm 1) \subset  \mathcal{V}(\mathfrak{so}_{(2N-1\pm 1)(2N+1)/2},1)$.
\end{itemize}
The sign differences correspond to level-rank duality pairs \cite{MR675108, MR1710999, MR1854694, Mirror,lrdual}. Note that the last three of these families fit into a more general setup $\mathcal{V}(\mathfrak{g},h^\vee) \subset  \mathcal{V}(\mathfrak{so}_{\dim g},1)$ where $h^\vee$ is the dual Coxeter number.

For the first family on the above list, Zhengwei Liu gave a diagrammatic description for the categories of $A$-modules for the corresponding \textit{\'etale} algebra $A$, and a diagrammatic description of a continuous family of tensor categories interpolating between them. This diagrammatic description was given as an equivariantization of a Yang-Baxter relation planar algebra \cite{LiuYB}. This approach was related to the Bisch-Jones project of understanding singly generated planar algebras \cite{MR1733737, MR1972635, MR3592517, MR1950890, MR4002229, MR4008521}. In this paper we develop a general technique for finding a diagrammatic presentation for the categories interpolating the quantum subgroup categories associated to the above discrete families of conformal embeddings. We restrict our attention in this paper to the second family $\mathcal{V}(\mathfrak{sl}_N,N) \subset  \mathcal{V}(\mathfrak{so}_{N^2-1},1)$, though we expect that slight variations of our approach will work for the remaining discrete families (More precisely, these categories interpolate between determinant-free versions of the quantum subgroup categories, paralleling that there's an $O(t)$ interpolation category but no $SO(t)$ interpolation category.)

Our key observation in this paper is the following. A quantum subgroup category always contains as a (non-full) subcategory the original quantum group category via the free module functor. So the remaining problem is describing the additional morphisms in the quantum subgroup which are not in the image of the free module functor. In our scenario, the free module functor gives us the well-known HOMFLYPT skein theory \cite{MR0776477,MR945888, turaev} for $\operatorname{Rep}(U_q(\mathfrak{sl}_N))$, where an oriented black strand denotes the image of the defining representation. We denote the object corresponding to this black strand as $X$. On the other hand, the local modules (which are completely understood for any conformal embedding) give us an invertible object $g$ generating a $\operatorname{Vec}(\mathbb{Z}_2)$ subcategory which we denote with a red strand. This allows us to give a diagrammatic presentation for our quantum subgroup category as a combined skein theory with both colors. The additional generator of this combined skein theory is a trivalent vertex with two black strands and one red strand. This comes from the fact that $g$ appears in $X \otimes X^*$, which we can deduce from the modular invariant of the conformal embedding.

The presentation we obtain is uniform with respect to $N$ and $q$ in the sense that the relations involve coefficients which are evaluations at $q= e^{2 \pi i \frac{1}{4N}}$ of rational functions in the variable $q$. This allows us to apply the techniques of Deligne interpolation to show existence of the category in the formal variable $q$.

We now state our main definitions and theorems.

As usual, in order to have a $1$-parameter family of tensor categories where we can specialize to particular complex numbers, we will work with an integral form over a localization of a polynomial ring, namely $R := \mathbb{C}[q, q^{-1},(q-q^{-1})^{-1}]$. The key point is that, unlike $\mathbb{C}(q)$, this ring has specializations maps $R \rightarrow \mathbb{C}$ setting $q$ to any complex number other than $0$, $1$, or $-1$.

\begin{defn}\label{def:main}
    Let $\mathcal{E}_R$ be the pivotal $R$-linear monoidal category with objects strings in $\{+,-\}$ and morphisms generated by the two morphisms
    \[  \att{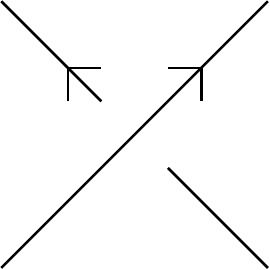}{.25}\in \End_{\mathcal{E}_R}(++)\quad,\quad\att{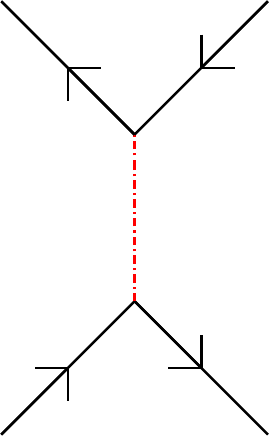}{.2} \in \End_{\mathcal{E}_R}(+-) \]
    satisfying the relations:
    \begin{align*}
    &(\text{Loop}) \att{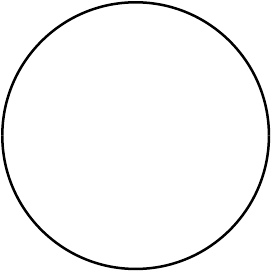}{.15} = \frac{2\mathbf{i}}{q-q^{-1}} \quad
    (\text{R1}) \att{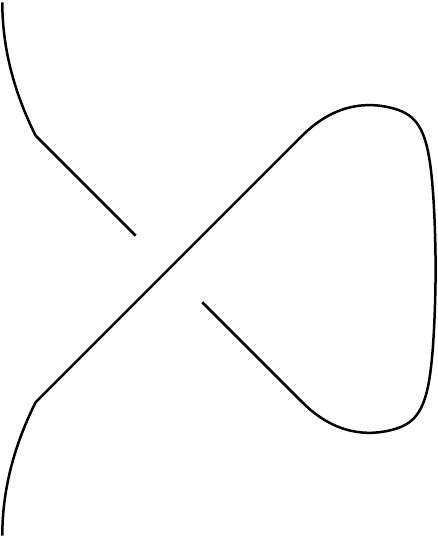}{.15} = \mathbf{i}\att{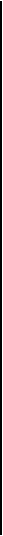}{.15}  \quad (\text{R2}) \att{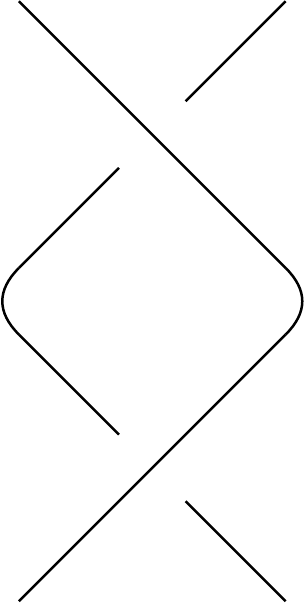}{.15} = \att{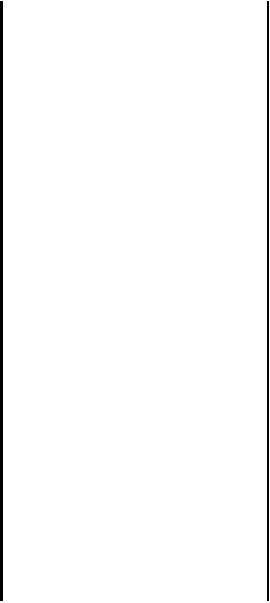}{.15}   \quad (\text{R3}) \att{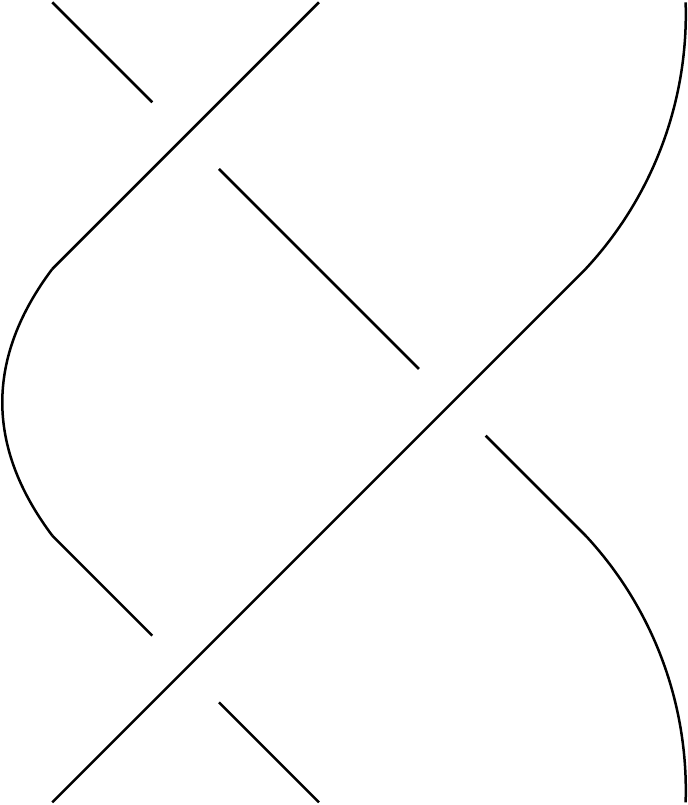}{.10} = \att{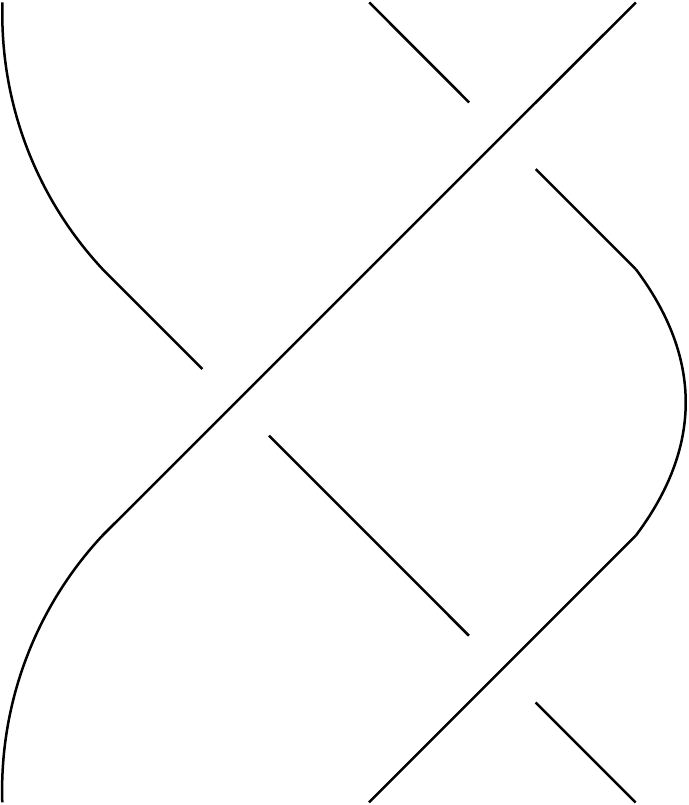}{.10}\\
    &(\text{Hecke})   \att{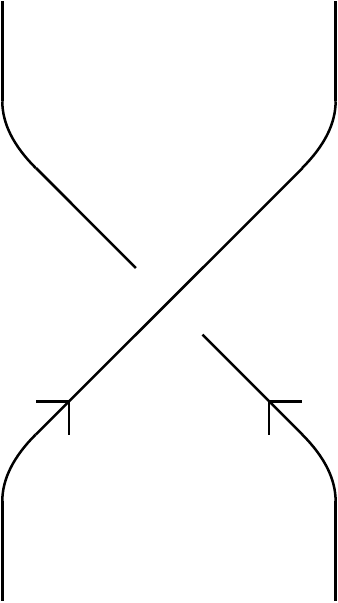}{.15} - \att{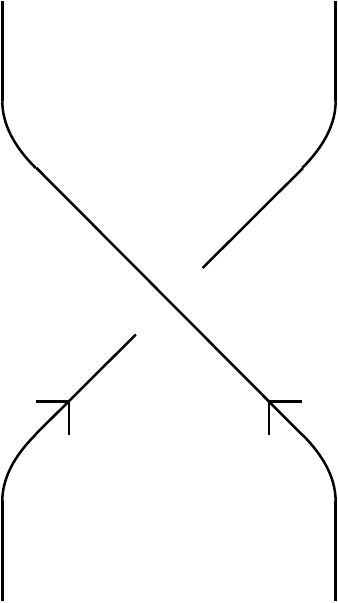}{.15} = (q-q^{-1})\att{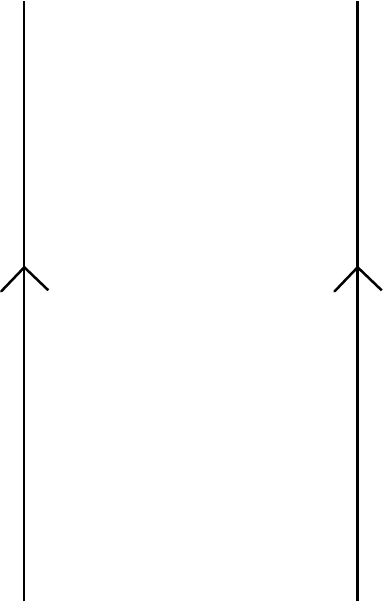}{.15}
   \quad  (\text{Trace}) \quad \att{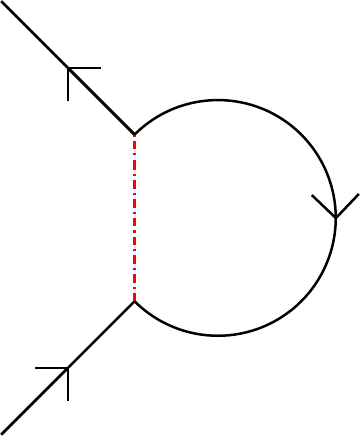}{.15} = \frac{q-q^{-1}}{2\mathbf{i}}\att{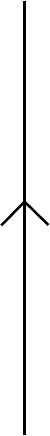}{.15}\quad (\text{Dual})  \left(\att{splittingEq}{.15}\right)^* =\att{splittingEq}{.15}\\
   &(\text{Half-Braid}) \quad \att{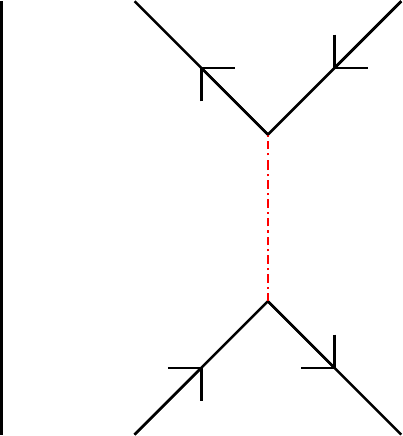}{.15} = \att{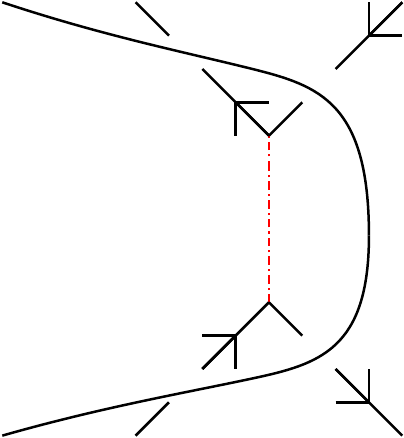}{.15} \quad (\text{Tadpole}) \att{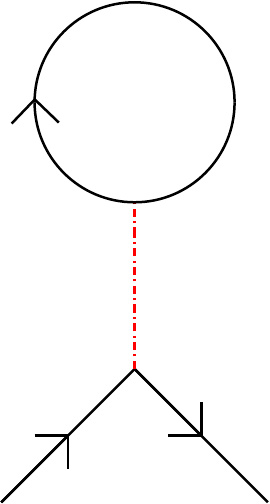}{.15} = 0\quad 
  (\mathbb{Z}_2)  \att{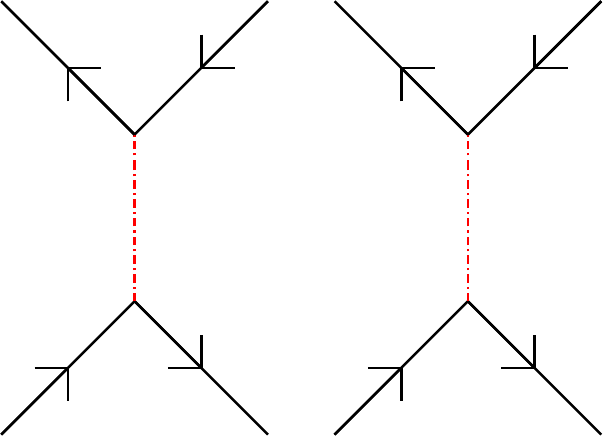}{.15} =  \att{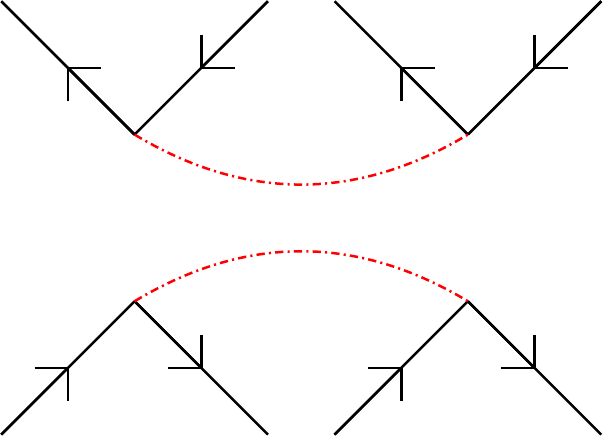}{.15}
    \end{align*}
\end{defn}
Here $*$ in the Dual relation represents the dual morphism, i.e. $180$-degree rotation applied to the diagram. We use the convention that a relation drawn using an un-oriented black strand holds for all possible orientations of that strand. Note that the horizontal reflection of (Half-Braid) holds, due to the relation (Dual). An immediate consequence of the relations (Hecke) and (R2) is the quadratic relation
\[\att{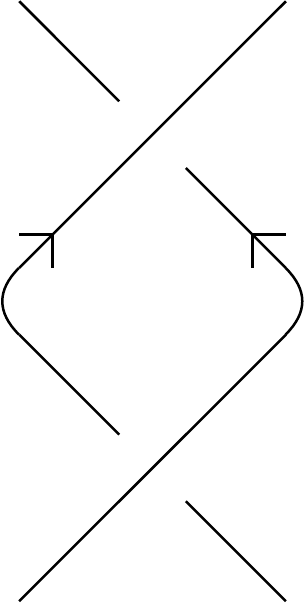}{.15} = \att{H12}{.15} + (q-q^{-1})\att{H11}{.15}\]
which we will use several times throughout this paper.

Note that although our description is written using crossings, the category $\mathcal{E}_R$ is not braided. In particular, the Half-Braid relation only holds for overcrossings and not for undercrossings, and so the crossing does not give a braiding natural in both variables. Geometrically, one should think of the trivalent vertices and red strands as glued to the ground, while black strands live in a 3-dimensional upper half-space. Thus black strands can pass over or under other black strands, but can only pass over (and not under) trivalent vertices (see \cite{MR2577673} and \cite[Remark 3.12]{MR3578212}, though the idea is implicit at least as early as \cite[\textsection 3.3]{MR1729094}). We do not rely in this paper on this geometric description, the crossing is simply treated as a generating morphism. This over-braiding phenomenon occurs in a tensor category $\mathcal{C}$ whenever you have a central functor \edit{(see Def. \ref{def:CentralFunctor}} from a braided tensor category $\mathcal{B}$ to $\mathcal{C}$ \cite[Def. 4.16]{braidedI}. The asymmetry between overcrossings and undercrossings comes from the definition of half-braidings, which are only natural with respect to morphisms on the strand passing under (see Section \ref{sec:etale}).

One unusual feature of this family is that it has no symmetric specialisations. This is because $q=\pm 1$ is the only symmetric specialization of the Hecke category, and our family is undefined at $q=\pm 1$ as the denominator of the loop relation vanishes.
\begin{defn}
    Let $t \in \mathbb{C} - \{0,1,-1\}$, and let $\mathcal{E}_t$ denote the $\mathbb{C}$-linear category given by specializing $q$ to $t$, i.e. $\mathbb{C} \otimes_R \mathcal{E}_R$ where $R$ acts on $\mathbb{C}$ by setting $q$ to $t$.
\end{defn}

We will often abuse notation by using $q$ both as a formal variable in $R$ and also as a complex number when we write $\mathcal{E}_q$.

When the parameter $q$ is of the form $e^{2\pi i \frac{1}{4N}}$ we define an extension of $\mathcal{E}_q$ by an additional generator corresponding to the determinant map.

\begin{defn}\label{def:main2}
    Let $N\in \mathbb{N}_{\geq 2}$. We define $\mathcal{SE}_N$ as the extension of $\mathcal{E}_{e^{2\pi i \frac{1}{4N}}}$ by the additional generators
    \[  \att{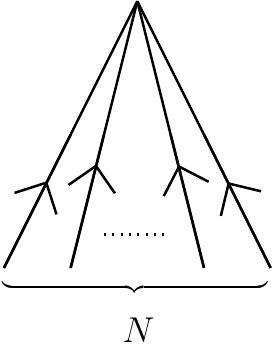}{.25},\edit{\qquad  \att{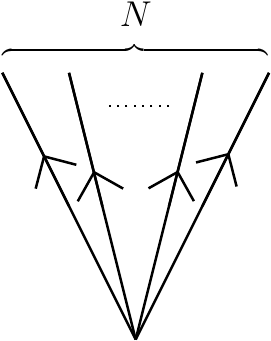}{.25}}  \]
    satisfying the relations
    \[  (\text{$q$-Braid}) \att{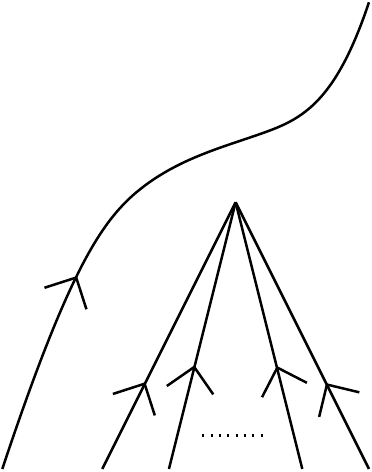}{.25} = q \att{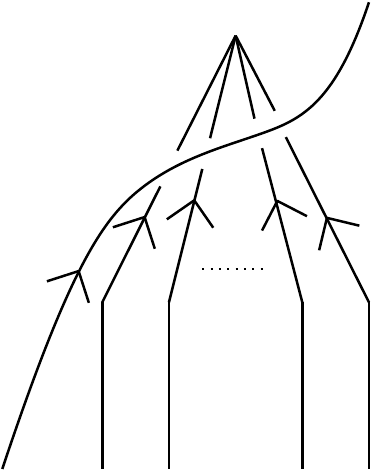}{.25}\qquad \text{ and }\qquad   (\text{Pair}) \att{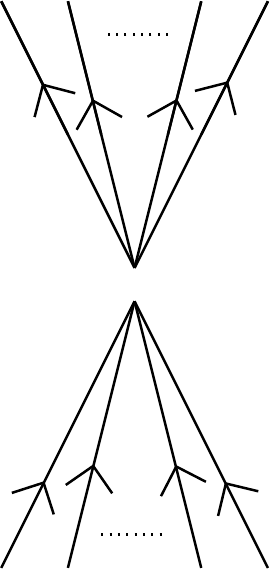}{.25} =  \att{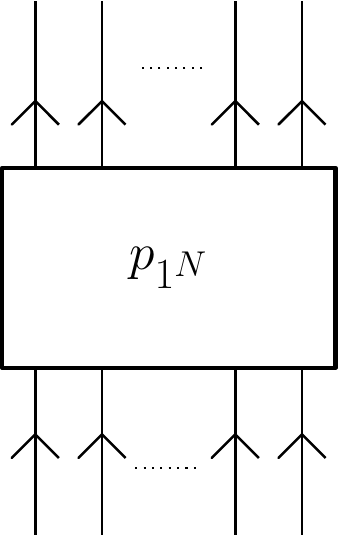}{.25}, \]
    where the box in the (Pair) relation is the projection onto the $N$-th quantum anti-symmetrizer (which \edit{is} expressed as a linear combination of braid elements).
\end{defn}

\begin{remark} \label{rem:ivsminusi}
    There is a (possibly non-unitary) Galois conjugate version of $\mathcal{SE}_N$ for primitive $4N$-th roots of unity such that $q^N = i$, but primitive $4N$-th roots of unity where $q^N = -\mathbf{i}$ behave differently because of the appearance of $\mathbf{i}$ in the defining relations of $\mathcal{E}_R$. 
\end{remark}

Our first main result shows that the Cauchy completion of the semisimplification of $\mathcal{SE}_N$ is equivalent to the tensor category corresponding to the conformal embedding $ \mathcal{V}(\mathfrak{sl}_N,N) \subset  \mathcal{V}(\mathfrak{so}_{N^2-1},1)$. This result justifies our claim that the family $\mathcal{E}_q$ interpolates between these conformal embedding categories.
\begin{thm}\label{thm:embedding}
    For each $N\in \mathbb{N}_{\geq 2}$ there is a monoidal equivalence between \[ 
 \operatorname{Ab}(\overline{\mathcal{SE}_{N}})\simeq \overline{\operatorname{Rep}(U_q(\mathfrak{sl}_N))}_A\]
 where $A$ is the \textit{\'etale} algebra object corresponding to the conformal embedding
    \[     \mathcal{V}(\mathfrak{sl}_N,N) \subset  \mathcal{V}(\mathfrak{so}_{N^2-1},1).      \]
\end{thm}
Here $\operatorname{Ab}$ denotes Cauchy completion (see Subsection~\ref{sec:Cauchy}) and overline denotes semisimplification (see Subsection~\ref{sec:semisimp}).

Our second main theorem gives a basis for $\mathcal{E}_R$. In particular, $\mathcal{E}_q$ is non-trivial for all $q \not\in \{0, \pm 1\}$.
\begin{thm}\label{thm:basisdim}
   Let $s_1,s_2$ be two strings in the alphabet $\{+,-\}$. We have that $\Hom_{\mathcal{E}_R}(s_1  \to s_2)$ is a free $R$-module of dimension $0$ if $\sum s_1 \neq \sum s_2$ and a free $R$-module of dimension 
    \[2^{\frac{|s_1|+|s_2|}{2}-1}\left(\frac{|s_1|+|s_2|}{2}\right)!\]  if $\sum s_1 = \sum s_2$.
\end{thm}
We achieve this result by giving explicit bases for the hom spaces of $\mathcal{E}_R$ in Theorem~\ref{thm:existence}. In fact, we use Theorem \ref{thm:embedding} as a tool to prove that the basis used in the proof of Theorem~\ref{thm:existence} is linearly independent. This is a typical argument in interpolation categories, where the existence of non-zero quotients at infinitely many specializations is used to show that the whole family is non-zero.

We expect that our techniques will work to give a diagrammatic description of the category associated to any infinite family of conformal embeddings. In the case of the embeddings of types $B$, $C$, and $D$, one has to work with the Kauffman skein theory \cite{MR958895,MR992598} instead of the HOMFLYPT skein theory. For Liu's family, you can use the HOMFLYPT skein theory together with an oriented red strand forming a $\mathrm{Vec}(\mathbb{Z})$ skein theory, and a trivalent vertex coming from a map $X \otimes X \rightarrow g$.

\begin{question}
    For generic values of the parameter Deligne interpolation categories are semisimple, but at special values they typically have both a semisimplification and an abelian envelope. Is there an abelian envelope of $\mathcal{E}_{e^{2\pi i \frac{1}{4N}}}$? One approach to this problem would be to construct versions of conformal embeddings over finite fields so one could take a limit in characterstic following \cite{harman}. Another approach would be to use the results of \cite{Kev, MR4533718}.
\end{question}

Our techniques for proving Theorem~\ref{thm:embedding} are non-constructive in the sense that we do not obtain the structure of $\operatorname{Ab}(\overline{\mathcal{SE}_{N}})$. In a subsequent paper, CE and Hans Wenzl will classify the simple objects, calculate the dimensions of simple objects, and describe the fusion rules for tensoring with $+$, for both $\operatorname{Ab}(\overline{\mathcal{E}_R})$ and $\operatorname{Ab}(\overline{\mathcal{SE}_{N}})$ \cite{HansNew}. Together with our Theorem~\ref{thm:embedding} this will provide a combinatorial description of the category $\overline{\operatorname{Rep}(U_q(\mathfrak{sl}_N))}_A$. 

We end the paper with a key result which is used in \cite{HansNew}, which is that the endomorphism algebras $\End_{\mathcal{E}_R}(+^n)$ are the even \edit{subalgebras} of the \textit{Hecke-Clifford} or \textit{affine Sergeev} algebras \cite{Ser,Ol, JN}. For generic $q$, these algebras were already studied \cite{JN} as they also appear as endomorphism algebras for the quantum group $U_q(\mathfrak{q}_N)$ where $\mathfrak{q}_N$ is the isomeric (or queer) Lie superalgebra \cite{Ol}. The approach in \cite{HansNew} is to analyze representations of the Hecke-Clifford algebras at roots of unity, and to study the novel trace on these algebras coming from the categorical trace on $\mathcal{SE}_N$.

This connection to isomeric Lie superalgebras is intriguing but remains mysterious. The paper \cite{Savage} introduces a tensor supercategory $Q(z)$, which is the large $N$ limit of the categories of rational representations of $U_q(\mathfrak{q}_N)$. The subcategory of $Q(z)$ spanned by diagrams with an even number of Clifford tokens is generated by the crossing and cup-cap diagram with two Clifford tokens, and if you change variables $z = q-q^{-1}$ and swap over/under crossing, the even part of $Q(z)$ has exactly the same generators and relations as our $\mathcal{E}_q$, except that the Dual relation has a minus sign, the R1 relation has $1$ instead of $i$, and the Loop value is $0$ instead of $2i/(q-q^{-1})$. Despite looking innocuous, this altered loop relation has a massive effect on the structure of the category. For example, the semisimplicification of the even part of $Q(q-q^{-1})$ is the trivial category, while the semisimplicification of $\mathcal{E}_q$ is much richer. 

The paper is outlined as follows.

In Section~\ref{sec:prelim} we review the required background material for this paper, and prove some basic results. We first review several categorical constructions: \textit{\'etale} algebra objects, Cauchy completion, and semisimplification. We then review the combinatorics of Young diagrams, and recall the definition of the HOMFLYPT skein categories. Finally we give the definition of the $q$-deformation of the Sergeev algebra, called the Hecke-Clifford algebras.

In Section~\ref{sec:con} we review the construction of an \textit{\'etale} algebra object $A\in \overline{\operatorname{Rep}(U_q(\mathfrak{sl}_N))}$ from the conformal embedding $\mathcal{V}(\mathfrak{sl}_{N},N) \subset \mathcal{V}(\mathfrak{so}_{N^2-1},1)$. We collate results from the literature to give an explicit description of the object $A$. Using this explicit description of $A$, we are able to show that $g$ is a summand of $X \otimes X^*$.
This information is the key to obtaining our presentation for the $A$-module category. Further, we also use the explicit description of $A$ to give explicit formulae for the dimensions of certain hom space in these categories. This result allows us to show the linear independence of our basis for these hom spaces later in the paper.

In Section~\ref{sec:sen} we show that the categories $\mathcal{SE}_N$ defined in Definition~\ref{def:main2} are non-trivial. We do this by constructing a monoidal equivalence between the Cauchy completion of $\overline{\mathcal{SE}_N}$, and the category of $A$-modules in $\overline{\operatorname{Rep}(U_q(\mathfrak{sl}_N))}$. Our definition of $\mathcal{SE}_N$ was obtained via studying this category of $A$-modules, and so a monoidal functor out of $\mathcal{SE}_N$ essentially comes for free. The hard part of the argument is showing that the functor is full. We reduce the question of fullness down to a question regarding the sub-algebra objects of $A$. Using classification results in the literature on these sub-algebras, we obtain the proof of Theorem~\ref{thm:embedding}. The nature of our proof is non-constructive, and so at this stage of the paper, we have little knowledge about the category $\operatorname{Ab}(\overline{\mathcal{SE}_N})$.

In Section~\ref{sec:eq} we construct a diagrammatic spanning set for the hom spaces of $\mathcal{SE}_N$. When the number of strands is small relative to $N$, we use our previous knowledge of the dimensions of these hom spaces to deduce that these spanning sets are in fact bases. Further, we obtain formulae for the tensor product and composition of these basis elements in terms of rational functions in $q$, evaluated at specific roots of unity. This set-up allows us to use the techniques of Deligne interpolation, to show that our diagrammatic sets are bases of the hom spaces of $\mathcal{E}_q$ for generic $q$. This gives the proof of Theorem~\ref{thm:basisdim}. Furthermore, we obtain that the endomorphism algebras $\operatorname{End}_{\mathcal{E}_q}(+^n)$ are isomorphic to the even part the Hecke-Clifford algebras. In particular this answers a question of Xu on the endomorphism algebras of $\overline{\operatorname{Rep}(U_q(\mathfrak{sl}_N))}_A$ \cite{Xuquestion}.

\subsection*{Acknowledgments}
CE was supported by NSF DMS grant 2245935 and 2400089. NS was supported by NSF DMS grant 2000093 and Simons Foundation grant MPS-TSM-00007608. Both authors would like to thank Hans Wenzl for his contributions to this project, and the anonymous referee for a careful reading which greatly improved the paper.
CE would like to thank Daniel Copeland for helpful conversations. NS thanks Zhengwei Liu for explaining the construction of his interpolation category and for suggesting this problem, Andr\'e Henriques for some pointers regarding free fermion algebras, and Harshit Yadav for some helpful pointers on Section 2.2. NS would like to thank CE and the University of New Hampshire for support and hospitality during a monthlong visit where this work was begun. This material is based upon work supported by the National Science Foundation under Grant No. DMS-1928930, while the authors were in residence at the Mathematical Sciences Research Institute in Berkeley, California, during the Summer of 2024. Both authors would like to thank BIRS for hosting them while part of this project was completed.

\section{Preliminaries}\label{sec:prelim}
We direct the reader to \cite{Book} for the basics on tensor categories, to \cite{sawin} for background on the categories $\overline{\operatorname{Rep}(U_q(\mathfrak{sl}_N))}$, to \cite{RigidCStar} for rigid $C^*$ tensor categories, and to \cite{WZW} for background on the categories $ \operatorname{Rep}(\mathcal{V}(\mathfrak{g},k))$.

\subsection{Cauchy Completion} \label{sec:Cauchy} 
In this subsection we define the \textit{Cauchy completion} of a tensor category. Informally, this construction completes the category to include direct sums and sub-objects. The Cauchy completion is defined as the additive envelope of the idempotent completion of a category. We refer the reader to \cite[Section 3]{Tuba} and \cite[Section 4.1]{d2n} for additional details.

\begin{defn} \label{def:KaroubiMotivation}
Suppose that $\mathcal{C}$ is a category, $X$ and $Y$ are objects in $\mathcal{C}$, and that $p_X: X \rightarrow X$ and $p_Y: Y \rightarrow Y$ is an idempotent endomorphism with images $i_X: X' \rightarrow X$ and $i_Y: Y' \rightarrow Y$, respectively. Let $H(X,Y)$ be the vector space \[\{f\in \Hom_\mathcal{C}(X\to Y) : p_X \circ f = f = f\circ p_Y \}.\] Then, by the universal property of images, the map $H(X,Y) \rightarrow \Hom_\mathcal{C}(X',Y)$ defined by $f \mapsto f \circ i_X$ descends to a map $K: H(X,Y) \rightarrow \Hom_\mathcal{C}(X',Y')$.
\end{defn}

\begin{lem} \label{lem:KaroubiMotivation}
    The map $K: H(X,Y) \rightarrow \Hom_\mathcal{C}(X',Y')$ in Definition \ref{def:KaroubiMotivation} is an isomorphism.
\end{lem}

The above lemma suggests the following definition of the idempotent completion of a tensor category.
\begin{defn}
    Let $\mathcal{C}$ be a pivotal tensor category. The objects of $\operatorname{Idem}(\mathcal{C})$ are pairs $(X,p_X)$ where $X\in \mathcal{C}$, and $p_X\in \operatorname{End}_\mathcal{C}(X)$ is an idempotent. The morphisms are defined by 
    \[ \operatorname{Hom}_{\operatorname{Idem}(\mathcal{C})}( 
  (X,p_X) \to (Y,p_Y)  ):= \{  f\in \Hom_\mathcal{C}(X\to Y) : p_X \circ f = f = f\circ p_Y        \}.  \]
  The identity morphism on $X$ is $p_X$.
  The tensor product, composition, and pivotal structure are inherited from the base category $\mathcal{C}$.
\end{defn}

The additive envelope is defined as follows.
\begin{defn}
    Let $\mathcal{C}$ be a pivotal $\mathbb{C}$-linear category. We define $\operatorname{Add}(\mathcal{C})$ as the category with objects formal finite direct sums
    \[  \bigoplus_i X_i   \]
    where each $X_i\in \mathcal{C}$, and morphisms matrices 
    \[  \begin{blockarray}{ccccc}
&Y_1  & Y_2 &\cdots & Y_m    \\
\begin{block}{c[cccc]}
X_1 & f_{1,1} & f_{1,2} & \cdots & f_{1,m} \\
X_2 & f_{2,1} & f_{2,2} & \cdots & f_{2,m} \\
\vdots & \vdots & \vdots& \ddots & \vdots \\
X_n &f_{n,1} & f_{n,2} & \cdots & f_{n,m} \\
\end{block}
\end{blockarray}\in \Hom_{\operatorname{Add}(\mathcal{C})}\left(\bigoplus_{i=1}^n X_i \to \bigoplus_{j=1}^m Y_j\right) \]
where $f_{i,j}\in \Hom_\mathcal{C}(X_i\to Y_j)$. The composition of morphisms is given by matrix composition, and the tensor product of morphisms is given by the Kronecker product. The pivotal structure is inherited from the pivotal structure on $\mathcal{C}$.
\end{defn}
With the above definitions we can now define the Cauchy completion.
\begin{defn}
The Cauchy completion $\operatorname{Ab}(\mathcal{C})$ of a category $\mathcal{C}$ is defined as the additive envelope of the idempotent completion of $\mathcal{C}$. We let $\iota_\mathcal{C}: \mathcal{C} \rightarrow \operatorname{Ab}(\mathcal{C})$ denote the canonical completion functor.
\end{defn}

The construction of taking the Cauchy completion is functorial, in the sense of the following lemma
\begin{lem}
    Let $\mathcal{F}: \mathcal{C}\to \mathcal{D}$ be a pivotal tensor functor, and $\mathcal{D}$ Cauchy complete. Then $\mathcal{F}$ extends uniquely to a pivotal tensor functor $\operatorname{Ab}(\mathcal{F}): \operatorname{Ab}(\mathcal{C})\to \mathcal{D}$. 
\end{lem}
\begin{proof}
This follows directly from the fact that the restriction functor
   \begin{equation}\label{eq:res}  \operatorname{Fun}_\otimes(\operatorname{Ab}(\mathcal{C})\to \mathcal{D}) \to \operatorname{Fun}_\otimes(\mathcal{C}\to \mathcal{D})        ,\quad  \mathcal{F} \mapsto \iota_\mathcal{C} \circ \mathcal{F}  \end{equation}
   is an equivalence \cite[Equations 7 and 9]{comes}.
\end{proof}
In particular, since $\operatorname{Ab}(\mathcal{D})$ is Cauchy complete, we have a functor \[\operatorname{Ab}(\iota_\mathcal{D} \circ \mathcal{F}): \operatorname{Ab}(\mathcal{C}) \rightarrow \operatorname{Ab}(\mathcal{C}) \rightarrow \operatorname{Ab}(\mathcal{D}).\]
We will often encounter functors which satisfy the following surjectivity-like condition on objects.
\begin{defn}
    A functor $\mathcal{C} \rightarrow \mathcal{D}$ is called Cauchy-dominant if any object in $\mathcal{D}$ is isomorphic to a direct sum of direct summands of objects in the image of $\mathcal{C}$.
\end{defn}

The following Lemma follows directly from the definitions.

\begin{lem}
    The completion functor $\iota_\mathcal{D}: \mathcal{D} \rightarrow \operatorname{Ab}(\mathcal{D})$ is faithful and Cauchy dominant.
\end{lem}

Since the compositon of faithful Cauchy-dominant functors is faithful and Cauchy dominant, we also have the following consequence.

\begin{lem} \label{lem:AbFaithful}
    If $\mathcal{F}: \mathcal{C} \rightarrow \mathcal{D}$ is faithful and Cauchy-dominant, then the completion $\operatorname{Ab}(\iota_\mathcal{D}\circ \mathcal{F}): \operatorname{Ab}(\mathcal{C}) \rightarrow \operatorname{Ab}(\mathcal{D})$ is faithful and Cauchy-dominant.
\end{lem}

A key property of Cauchy complete categories is the following characterization of equivalences.

\begin{lem}
    Suppose that $\mathcal{C}$ and $\mathcal{D}$ are Cauchy-complete, then a functor $\mathcal{F}: \mathcal{C} \rightarrow \mathcal{D}$ is an equivalence iff it is full, faithful, and Cauchy-dominant.
\end{lem}
\begin{proof}
    The only-if direction is obvious, since essentially surjective implies  Cauchy-dominant. We turn to the if direction.
    Since the functor is full and faithful, any idempotent in $\mathcal{D}$ can be pulled back to one in $\mathcal{C}$, and so the functor is essentially surjective.
\end{proof}

\subsection{\textit{\'Etale} algebras and their tensor categories of modules}\label{sec:etale}

In this subsection we briefly review the theory of \textit{\'etale} algebra objects in a braided fusion category $\mathcal{C}$, and the correspondence between \textit{\'etale} algebra in $\mathcal{C}$ and fusion categories over $\mathcal{C}$. Additional details can be found in \cite{LagrangeUkraine,Exact,FancyEtaleCorrespondence, ShimizuYadav, Kril}.

\begin{defn}
Let $\mathcal{C}$ be a braided fusion category. An algebra object $A\in \mathcal{C}$ is said to be \textit{\'etale} if it is commutative and separable (i.e. the multiplication map $A\otimes A\to A$ splits as a map of $A$-$A$ bimodules). We say that $A$ is connected if $\dim\Hom_\mathcal{C}(A\to \mathbf{1})=1$. Let $\mathrm{Alg}_\mathcal{C}$ denote the category of algebra objects in $\mathcal{C}$ with algebra morphisms, and let $\mathrm{Alg}^\text{\'et}_\mathcal{C}$ denote the full subcategory of \textit{\'etale} algebras.
\end{defn}
Given an \`etale algebra object, we have the following key map.
\begin{lem}\label{lem:EtaleToFrob} \cite[Prop. 3.1(ii)]{OstMod}
    Suppose that $A$ is a connected \'etale algebra object in a braided fusion category, and let $\tr:A \rightarrow 1$ be a splitting of the unit map. Then $\tr \circ \mu: A\otimes A\to \mathbf{1}$ gives a non-degenerate pairing between $A$ and itself.
\end{lem}

\begin{remark} \label{rem:Frobenius}
Lemma \ref{lem:EtaleToFrob} says that any \'etale algebra satsifies one of many equivalent definitions of a commutative Frobenius algebra object. We will use the usual diagram calculus for commutative Frobenius algebras objects, where a trivalent vertex denotes multiplication, and cup and cap are the pairing and copairing in Lemma \ref{lem:EtaleToFrob}. See \cite[\S 1]{Kril} and \cite[\S 2.3]{MR2187404}.
\end{remark}

Given an \textit{\'etale} algebra object $A\in \mathcal{C}$ we can consider the category $\mathcal{C}_A$ of left $A$-modules in $\mathcal{C}$. We can equip $\mathcal{C}_A$ with the structure of a tensor category as follows.

\begin{defn}
Let $A\in \mathcal{C}$ be an \textit{\'etale} algebra object, and let $M\in \mathcal{C}_A$ be a left $A$-module. Using the braiding on $\mathcal{C}$, we can 
equip $M$ with the structure of a right $A$-module via
\[  M\otimes A \xrightarrow{c_{M,A}} A\otimes M \to M .\]
We can then define a tensor product on $\mathcal{C}_A$ as the relative tensor product over $A$ between a left $A$-module and a right $A$-module, Explicitly, for left $A$-modules $M_1,M_2$, we can define $M_1\otimes_A M_2$ as the image of the projection
\[ \att{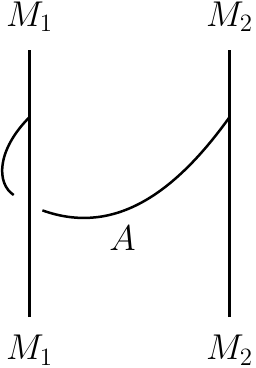}{.4}.  \]

Here the cup is the copairing from Remark \ref{rem:Frobenius}. That this map gives projection onto the relative tensor product is \cite[Lemma 1.21]{Kril}.
\end{defn}

\begin{remark}
    There's two choices of conventions required in defining $\mathcal{C}_A$: you could consider left $A$-modules or right $A$-modules, and could use the braiding or its inverse in the definition of tensor product. Left $A$-modules with $c_{M,A}$ is equivalent as using right $A$-modules with $c_{A,M}^{-1}$. But choosing left $A$-modules with $c_{M,A}$ or $c_{M,A}^{-1}$ leads to slight differences in the choice of crossings in other definitions. For example, because we use $c_{M,A}$ in our definition of $\mathcal{C}_A$, the crossing in the projection onto $M_1 \otimes_A M_2$ has to go \emph{under} in order to coequalize the $A$ actions.
\end{remark}

Recall from Definition \ref{def:KaroubiMotivation} and Lemma \ref{lem:KaroubiMotivation}, that for projections $\pi_1: X \rightarrow X$ and $\pi_2: Y \rightarrow Y$, we have a canonical identification of $\Hom(\im(\pi_1)\rightarrow \im(\pi_2))$ with the subset of $\Hom(X\rightarrow Y)$ satisfying $\pi_2 f = f = f \pi_1$. When we write maps between relative tensor products we will write the map from the ordinary tensor product, and then there's a routine check to see that it does descend to the relative tensor product. Keep in mind, that under this correspondence the identity morphism on $X$ corresponds to $\pi_1$.

The following structure results on $\mathcal{C}_A$ are well-known.
\begin{lem}
    If $A$ is a connected \textit{\'etale} algebra object in $\mathcal{C}$, then $\mathcal{C}_A$ with the above tensor product is a fusion category.
\end{lem}
\begin{proof}
    The condition that $A$ is separable implies that $\mathcal{C}_A$ is semisimple \cite[Proposition 2.7]{LagrangeUkraine}, and the condition that $A$ is connected implies that $\mathcal{C}_A$ has simple unit \cite[Remark 3.4]{LagrangeUkraine}.
\end{proof}

\begin{lem}
    \cite[Rem. 1.19]{Kril} If $\mathcal{C}$ is a ribbon fusion category, and the twist $\theta_A = \id_A$, then $\mathcal{C}_A$ is spherical.
\end{lem}

We then have the free module functor $\mathcal{F}_A : \mathcal{C}\to \mathcal{C}_A$.
\begin{defn}
    Let $A\in \mathcal{C}$ be an \textit{\'etale} algebra object. The free module functor $\mathcal{F}_A : \mathcal{C}\to \mathcal{C}_A$ is defined by $X\mapsto A\otimes X$ where the $A$-module action is by multiplying with the first tensor factor. 
\end{defn}

A direct computation shows the following lemma.

\begin{lem}
Let $\mathcal{A}$ be an \`etale algebra object, then $\mathcal{F}_A: \mathcal{C}\to \mathcal{C}_A$ is a strong tensor functor with tensorator
\[ \tau_{X,Y} := \att{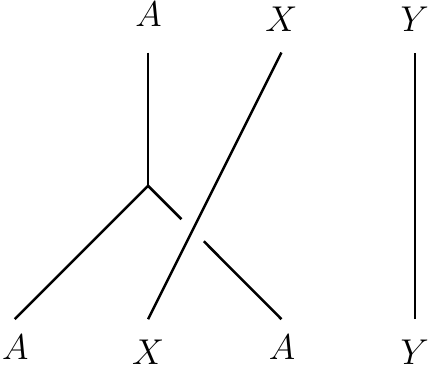}{.3}  \] 
and inverse tensorator 
\[\tau^{-1}_{X,Y} = \att{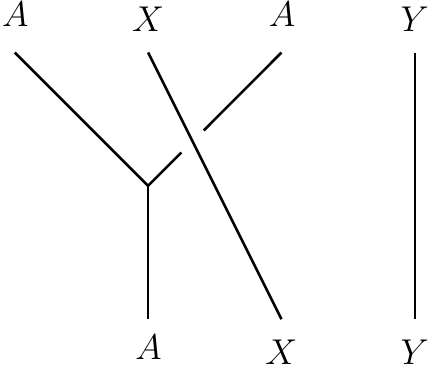}{.3}.\]
\end{lem}

The following is well-known to experts.
\begin{lem}\label{lem:dom}
    The monoidal functor $\mathcal{F}_A: \mathcal{C} \to \mathcal{C}_A$ is Cauchy-dominant.
\end{lem}
\begin{proof}
    Every module is a quotient of the free module on its underlying object, but since $\mathcal{C}_A$ is semisimple, any quotient is a direct summand.
\end{proof}


In general, the tensor category $\mathcal{C}_A$ will not have a braided structure. However, we can show that $\mathcal{C}_A$ will admit a ``half-braiding''. This will follow from known results regarding central functors, and central structures.

\begin{defn} \label{def:CentralFunctor}
    Let $\mathcal{C}$ be a braided fusion category, and $\mathcal{D}$ a fusion category. A \emph{central functor} is a functor $\mathcal{F}: \mathcal{C} \rightarrow \mathcal{D}$ endowed with a compatible lift to the center $\mathcal{C} \rightarrow Z(\mathcal{D})$. A fusion category \emph{over $\mathcal{C}$} is a fusion category $\mathcal{D}$ equipped with a Cauchy-dominant central functor $\mathcal{C} \rightarrow \mathcal{D}$ called the action functor.
    
    If $\mathcal{D}_1$ and $\mathcal{D}_2$ are fusion categories over $\mathcal{C}$ (with action functors $\mathcal{F}_1$ and $\mathcal{F}_2$), then a $\mathcal{C}$-linear functor is a functor $\mathcal{G}: \mathcal{D}_1 \rightarrow \mathcal{D}_1$ endowed with a natural isomorphism of central functors $\mathcal{G} \circ \mathcal{F}_1 \cong \mathcal{F}_2.$

    Let $\mathrm{Over}_\mathcal{C}$ denote the category of fusion categories over $\mathcal{C}$ with morphisms the natural isomorphism classes of $\mathcal{C}$-linear functors.
\end{defn}

The following lemma directly expands the definition of central functor.

\begin{lem}
    Endowing $\mathcal{F}: \mathcal{C} \rightarrow \mathcal{D}$ with a central structure is equivalent to giving for every object $c$ in $\mathcal{C}$ and every object $d$ in $\mathcal{D}$ a half-braiding $\sigma_{c,d}: \mathcal{F}(c) \otimes d \rightarrow d \otimes \mathcal{F}(c)$ such that \[\att{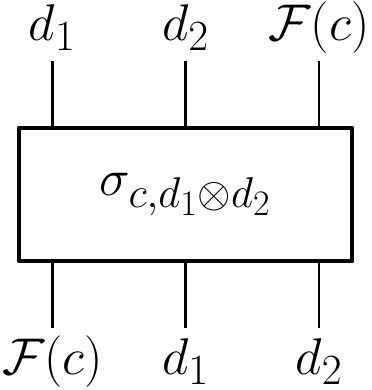}{.3} \quad= \quad \att{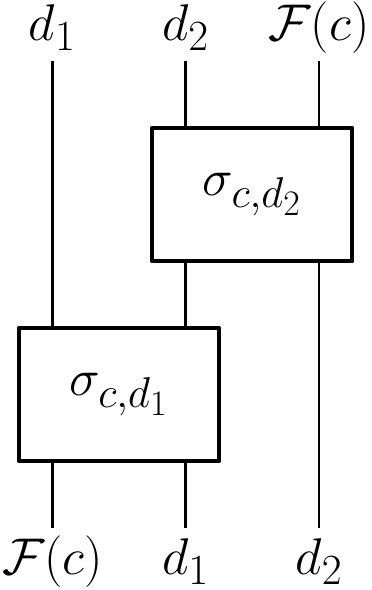}{.3}\] and such that for every $g\in \Hom_{\mathcal{D}}(d_1 \rightarrow d_2)$, we have naturality in the second variable: \[\att{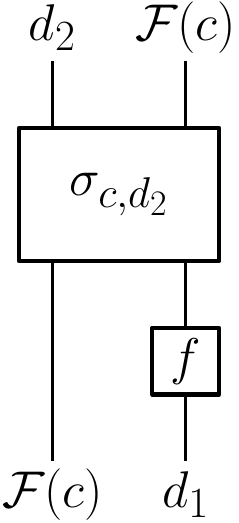}{.3} \quad = \quad \att{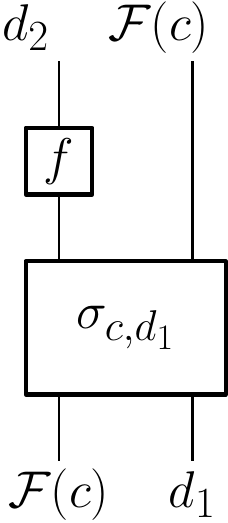}{.3}\]
\end{lem}

\begin{remark}\label{rem:OvervsUnder1}
If $\mathcal{D}$ is a tensor category over $\mathcal{C}$ and $x$ is an object in $\mathcal{C}$, we will use a crossing in our graphical calculus for the half-braiding $\sigma$. Some care must be taken, if $g\in \Hom_\mathcal{D}(y \rightarrow y')$, then naturality lets us pull $g$ through a crossing, i.e. $\sigma_{x,y} \circ (\id \otimes g) = (g \otimes \id) \circ \sigma_{x,y'}$. However, we cannot pull morphisms \emph{over} crossings. That is, if $f \in \Hom_\mathcal{D}(\mathcal{F}(x) \rightarrow \mathcal{F}(x'))$ then there is no relationship between $\sigma_{x,y} \circ (f \otimes \id)$ and $(\id \otimes f) \circ \sigma_{x',y}$. (We only have for $f' \in \Hom_\mathcal{C}(x \rightarrow x')$ then $\sigma_{c,d} \circ (\mathcal{F}(f') \otimes \id)$ and $(\id \otimes \mathcal{F}(f')) \circ \sigma_{c',d}$.) 
\end{remark}
A key source of central functors come from \`etale algebras in braided tensor categories.
\begin{lem}
    Let $A\in \mathcal{C}$ be an \`etale algebra object. Then $\mathcal{F}_A: \mathcal{C} \rightarrow \mathcal{C}_A$ can be endowed with a central structure via 
    \[\sigma_{X, M} := \att{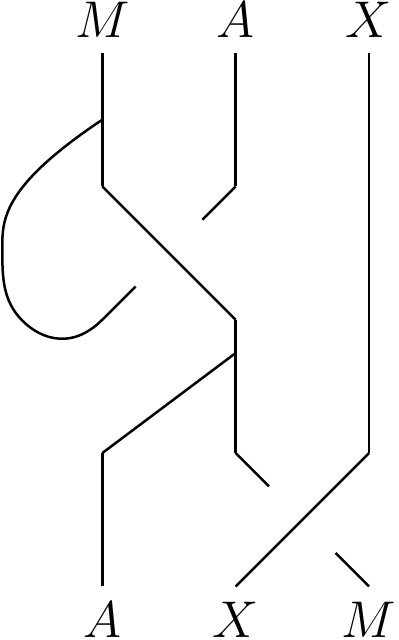}{.3}.\] 
\end{lem}
\begin{proof}
    This is a direct computation via the graphical calculus.
\end{proof}

We now relate the above morphisms $\sigma_{A \otimes X, M}$ to elements in the image of the free module functor.
\begin{lem} \label{lem:BraidingIsBraiding}
    Suppose that $A$ is an \textit{\'etale} algebra object in $\mathcal{C}$. If $X$ and $Y$ are objects in $\mathcal{C}$, then the half-braiding in $\mathcal{C}_A$ satisfies \[ \sigma_{X,\mathcal{F}_A(Y)} = \tau_{X,Y} \circ \mathcal{F}_A(c_{X,Y}) \circ\tau_{X,Y}^{-1} :\mathcal{F}_A(X)\otimes \mathcal{F}_A(Y)  \to \mathcal{F}_A(Y)\otimes \mathcal{F}_A(X).\]
\end{lem}
\begin{proof}
    Expanding the definitions, and using commutativity of the product, both sides of the equation are given by 
    \begin{equation*}
    \att{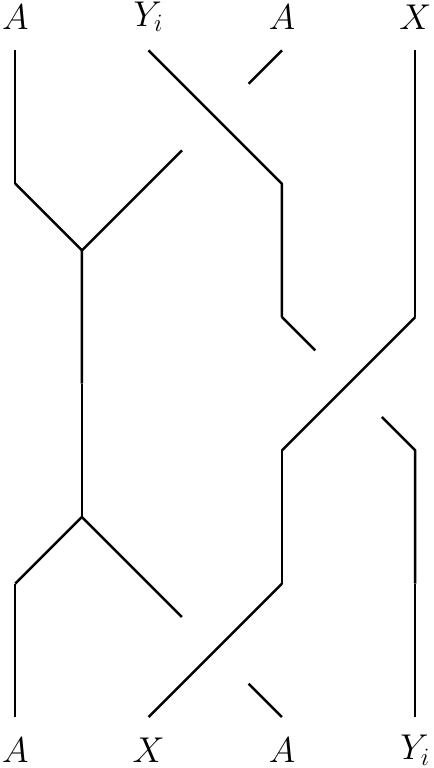}{.25}.
    \end{equation*}
\end{proof}

\begin{cor}\label{cor:OB}
Let $\mathcal{C}$ be a braided tensor category, $A\in \mathcal{C}$ an etale algebra object, and $f\in \operatorname{Hom}_{ \mathcal{C}_{A}}(\mathcal{F}_A(Y_1)\to \mathcal{F}_A(Y_2))$, then 
 \[\att{overproof2}{.3} =\att{overproof1}{.3} .  \] 
\end{cor}
\begin{proof}
    This follows from Lemma \ref{lem:BraidingIsBraiding} combined with naturality of $\sigma$.
\end{proof}

\begin{remark} \label{rem:OvervsUnder2}
    Note that as in Remark \ref{rem:OvervsUnder1}, the above corollary only holds for the \emph{under} strand, not the the \emph{over} strand. This asymmetry is visible in the structure of $\sigma_{X, M}$, because a coupon attached to the last two strands isotopes under $X$, but a coupon attached to the first two strands cannot isotope past the $Y$ strand because $A$ and $X$ strands pass on opposite sides of the $Y$ strand.
\end{remark}

Having constructed a fusion category over $\mathcal{C}$ from an \'etale algebra in $\mathcal{C}$, we now turn to the converse construction. The following two lemmas are standard.

\begin{lem}
If $\mathcal{F}: \mathcal{C} \rightarrow \mathcal{D}$ is a Cauchy-dominant functor of fusion categories, then $\mathcal{F}$ has a right adjoint $\mathcal{F}^R$, which has a canonical structure of a lax monoidal functor.
\end{lem}

\begin{lem} \label{lem:OverToEtale1}
    If $\mathcal{G}: \mathcal{C} \rightarrow \mathcal{D}$ is an lax monoidal functor of fusion categories, then $\mathcal{G}$ induces a functor $\mathcal{G}: \mathrm{Alg}_\mathcal{C} \rightarrow \mathrm{Alg}_\mathcal{D}$.
\end{lem}
    
In particular if $\mathcal{D}$ is a fusion category over $\mathcal{C}$, then $\mathcal{F}^R(\mathbf{1})$ is an algebra object.

\begin{lem} \label{lem:OverToEtale2}
    If $\mathcal{D}$ is a fusion category over $\mathcal{C}$, then $\mathcal{F}^R(\mathbf{1})$ is an \textit{\`etale} algebra object in $\mathcal{C}$.
\end{lem}

Now we can state the main result of \cite[\S 3.1-\S 3.3]{LagrangeUkraine}.

\begin{thm} 
    The assignments $A \mapsto \mathcal{C}_A$ and $\mathcal{D} \mapsto \mathcal{F}^R(\mathbf{1})$ 
    give inverse bijections between the set of isomorphism classes of connected \textit{\'etale} algebras in $A$ and equivalence classes of fusion categories over $\mathcal{A}$.
\end{thm}

Finally we will need one new lemma.

\begin{lem}\label{lem:hack}
    Suppose that $\mathcal{G}: \mathcal{D}_1 \rightarrow \mathcal{D}_2$ is a Cauchy-dominant faithful $\mathcal{C}$-linear functor of fusion categories, then there is an inclusion of algebras $f: A_1 = \mathcal{F}_1^R(\mathbf{1}) \rightarrow \mathcal{F}_2^R(\mathbf{1}) = A_2$. Moreover, if $f$ is an isomorphism on objects, then $\mathcal{D}_1$ and $\mathcal{D}_2$ are equivalent as fusion categories.
\end{lem}
\begin{proof}
    Since $\mathcal{G}$ is a strong monoidal functor (and in particular oplax monoidal), it follows that $\mathcal{G}^R$ is lax monoidal and $G^R(\mathbf{1})$ is an algebra. We have a unit map of algebras $u: \mathbf{1} \rightarrow \mathcal{G}^R(\mathbf{1})$. 
   By Lemma \ref{lem:OverToEtale1} and Lemma \ref{lem:OverToEtale2}, we get a map of algebras $\mathcal{F}_1^R(u): A_1=\mathcal{F}_1^R(\mathbf{1}) \rightarrow \mathcal{F}_1^R(\mathcal{G}^R(\mathbf{1})) \cong  \mathcal{F}_2^R(\mathbf{1})) = A_2$.

    By \cite[Lemma 3.1]{LagrangeUkraine}
    \[\mathrm{FPdim}(\mathcal{D}_1) = \mathrm{FPdim}\edit{(}\mathcal{C})/\mathrm{FPdim}(A_1) = \mathrm{FPdim}\edit{(}\mathcal{C})/\mathrm{FPdim}(A_2) = \mathrm{FPdim}(\mathcal{D}_2).\]
    Hence, the Cauchy-dominant functor $\mathcal{G}$ is an equivalence by \cite[Proposition 6.3.4]{Book}.
\end{proof}

Although Lemma \ref{lem:hack} sufficient for our puproses, it is not very satisfying. The following conjecture would give an improved statement, but requires a careful check that the existing constructions are functorial.

\begin{con} \label{con:Nonhack}
    There's a functor $ \mathrm{Alg}^\text{\'et}_\mathcal{C} \rightarrow \mathrm{Over}_\mathcal{C}$ which on objects sends $A \mapsto \mathcal{C}_A$ and on morphisms sends $f: A \rightarrow B$ to base extension $M \mapsto B \otimes_A M$. There is also a functor $\mathrm{Over}_\mathcal{C} \rightarrow \mathrm{Alg}^\text{\'et}_\mathcal{C}$ which on objects sends $\mathcal{D} \mapsto \mathcal{F}^R(\mathbf{1})$ and sends a morphism $(\mathcal{G},\eta): \mathcal{F}_1 \rightarrow \mathcal{F}_2$ to $\mathcal{F}_1^R(u)$ where $u: \mathbf{1} \rightarrow \mathcal{G}^R(\mathbf{1})$ is the unit map. These functors are inverse equivalences.
\end{con}

In particular, since the right adjoint of faithful functor preserves monomorphisms, the equivalence in Conjecture \ref{con:Nonhack} would induce an order preserving bijection between the poset of fusion categories over $\mathcal{C}$ ordered by faithful Cauchy-dominant $\mathcal{C}$-linear functors and the poset of \textit{\'etale} algebras in $\mathcal{C}$ ordered by inclusion.


\subsection{Unitarity, negligible morphisms and semisimplification} \label{sec:semisimp}
In this section we define the negligble ideal and the semisimplification functor (see \cite{MR1686423,Simp}), and give several of its properties, especially how these behave in the context of unitary fusion categories (see \cite{bicomm}).

\begin{defn}
    Let $\mathcal{C}$ be a spherical category. We define the negligible ideal of $\mathcal{C}$ as
    \[ \operatorname{Neg}(\mathcal{C}) := \left\{ f\in \Hom_\mathcal{C}(X\to Y) : \operatorname{tr}(f\circ g) = 0 \text{ for all } g\in\Hom_\mathcal{C}(Y\to X) \right\}  \]
    where $\operatorname{tr}$ is the categorical trace.
\end{defn}

It is well known that $\operatorname{Neg}(\mathcal{C})$ is a tensor ideal of $\mathcal{C}$ \cite[Lemma 2.3]{Simp}.  Hence we can form the quotient category.
\begin{defn}
    Let $\mathcal{C}$ be a spherical category. We will write
    \[   \overline{\mathcal{C}}:= \mathcal{C} /   \operatorname{Neg}(\mathcal{C}).    \]
    We will call the quotient map $\mathcal{C}\to \overline{\mathcal{C}}$ the semisimplification functor.
\end{defn}
In the case that $\mathcal{C}$ has simple unit, the negligible ideal is the unique maximal ideal of $\mathcal{C}$ \cite[Proposition 3.5]{EHOld}. 


The negligible ideal will be particularly amenable in the context of \textit{rigid $C^*$ tensor categories}.

\begin{defn}
  A rigid $C^*$ tensor category, is a $\mathbb{C}$-linear dagger monoidal category whose underlying category is a $C^*$ category, and which is rigid (see \cite[\S 2.1-\S 2.3]{RigidCStar} for details).
\end{defn}

Note that we do not assume that rigid $C^*$ tensor categories are abelian, or even Cauchy-complete.

\begin{defn}
    A unitary tensor category $\mathcal{C}$ is a rigid $C^*$ tensor category such that all Hom spaces are finite dimensional and $\End(1) = k$.  If in addition if $\mathcal{C}$ is Abelian and has finitely many isomorphism classes of simple objects and we call it a unitary fusion category.
\end{defn}

\begin{lem}
    Any unitary tensor category has a canonical spherical structure. Any Cauchy complete unitary category is Abelian semisimple.
\end{lem}
\begin{proof}
    The endomorphism algebra of any object is a finite-dimensional $C^*$ algebra, and hence is a semisimple algebra. Then \cite[Theorem 3.3]{Tuba} implies that the category is \edit{abelian} semisimple.
    The spherical structure is described in \cite{bicomm} Def. 2.7.
\end{proof}

 In the context of functors into unitary tensor categories, we have the following useful descent lemma.
\begin{prop}\label{prop:descent}
    Let $\mathcal{C}$ be a spherical monoidal category, let $\mathcal{D}$ be a unitary tensor category, and let $\mathcal{F}:\mathcal{C}\to \mathcal{D}$ be a pivotal functor. Then there exists a faithful pivotal functor $\overline{\mathcal{F}}: \overline{\mathcal{C}}\to \mathcal{D}$ such that the following diagram commutes
      \[ \begin{tikzcd}
\mathcal{C} \arrow[d] \arrow[r, "\mathcal{F}"] & \mathcal{D} \\
\overline{\mathcal{C}} \arrow[ur, "\overline{\mathcal{F}}"]& 
\end{tikzcd}\]
where $\mathcal{C}\to \overline{\mathcal{C}}$ is the semisimplification functor. Moreover, $\overline{\mathcal{C}}$ is a unitary category.
\end{prop}

Although this result is well-known to experts, we were not able to find a proof in the literature so we provide one.
\begin{proof}
Suppose that $h \in \Hom_\mathcal{C}(X\to Y)$ is a negligbile morphism. Thus $h$ is in the radical of the trace pairing on the subspace of $\Hom_\mathcal{C}(X\to Y)$ which is the image of $\mathcal{F}$. In a unitary tensor category, in addition to the usual bilinear form $\tr(f \circ g)$, there is also the positive-definite sesquillinear form $\tr(f^\dagger \circ g)$. Since $\dagger$ is an involution, \edit{pairing trivially with all $f^\dagger$ is the same as pairing trivially with all $f$, so} these two forms have the same radical. The restriction of a positive definite form to any subspace is again positive definite, and thus has no radical. Thus $\mathcal{F}(h) = 0$. Moreover, $\overline{\mathcal{C}}$ inherits the unitary structure from $\mathcal{D}$.
\end{proof}


\begin{warning}
    Note that without some version of unitarity the previous proposition is false, since the restriction of a non-degenerate inner product may be degenerate. 
    An explicit counter-example can be found in \cite[Example 1.4]{knots}. 
\end{warning}

\subsection{Combinatorics of \texorpdfstring{$\mathrm{GL}(N)$}{GL(N)} and \texorpdfstring{$\mathrm{SL}(N)$}{SL(N)}}\label{sec:GLNcombo}

We briefly recall some facts about irreducible representations and fusion rules for representations of $\mathrm{SL}(N)$ and $\mathrm{GL}(N)$. 

Note that since $\mathrm{SL}(N)$ is simply connected there's no difference between representations of $\mathrm{SL}(N)$ and representations of $\mathfrak{sl}_n$, but we will often wish to refer to rational or polynomial representations of $\mathrm{GL}(N)$ which are more naturally thought of in terms of the Lie group than the Lie algebra.

The category of representations of a group is graded by the group of central characters $\Hom(Z(G),\mathbb{C})^\times$. In particular, representations of $\mathrm{GL}(N)$ are graded by $\Hom(Z(G),\mathbb{C})^\times = \mathbb{Z}$ and representations of $\mathrm{SL}(N)$ are graded by $\mathbb{Z}/N$. Polynomial representations of $\mathrm{GL}(N)$ are all in non-negative grading.

\begin{lem} \label{lem:GradingSLGL}
    If $X$ and $Y$ are representations of $\mathrm{GL}(N)$ with the same grading, then restriction gives an equality $\Hom_{\mathrm{GL}(N)} (X,Y) = \Hom_{\mathrm{SL}(N)} (X,Y)$
\end{lem}
\begin{proof}
    Since $\mathrm{GL}(N)$ is generated by its center together with $\mathrm{SL}(N)$, a map that respects the action of the center and the action of $\mathrm{SL}(N)$ automatically respects the action of $\mathrm{GL}(N)$.
\end{proof}

The standard Cartan subalgebra $\mathfrak{h}$ of $\mathfrak{gl}(N)$ is given by diagonal matrices. Let $\varepsilon_i\in \mathfrak{h}^*$ be the functional $e_{jj} \mapsto \delta_{i,j}$. A weight $\sum a_i \varepsilon_i$ is dominant if $a_1 \geq a_2 \geq \ldots a_N$ and integral if all the $a_i \in \mathbb{Z}$. A dominant integral weight is non-negative if all the $a_i \in \mathbb{N}$. Irreducible finite-dimensional representations are classified by their highest weight, which must be a dominant integral weight. Polynomial irreducible representations correspond to weights which are in addition non-negative. 

A non-increasing list of natural numbers can be visualized as a Young Diagram. We draw our Young Diagrams in the English notation, that is a Young Diagram is made up of several rows of boxes, which are left-justified and where each row contains no more boxes than the previous ones. 

\begin{defn}
    If $\lambda$ is a Young diagram
    \begin{enumerate}[(a)]
        \item let $\lambda_i$ be the length of the $i$th row,
        \item let $\ell(\lambda)$ be the number of rows,
        \item let $|\lambda|$ be the total number of boxes,
        \item and let $H_\lambda(1,1)= \lambda_1 + \ell(\lambda)-1$ be the hook-length attached to the the top-left box,
        \item $\lambda^T$ is the diagram whose rows are the columns of $\lambda$.
    \end{enumerate}
\end{defn}

We can assign to a Young diagram $\lambda$ with $\ell(\lambda)\leq N$ the non-negative dominant integral weight $\sum \lambda_i \varepsilon_i$, and hence an irreducible polynomial representation $W_\lambda$. The grading of such a representation is $|\lambda|$.

There's an important modification of this construction that allows us to consider rational representations as well. If $(\lambda,\mu)$ is a pair of Young diagrams with $\ell(\lambda) + \ell(\mu) \leq N$ we define the corresponding highest weight to be 
\[   \lambda_1 \varepsilon_1 + \lambda_2   \varepsilon_2 + \cdots \lambda_{\ell(\lambda)}  \varepsilon_{\ell(\lambda)} - \mu_{\ell(\mu)} \varepsilon_{N+1-\ell(\mu)}  -\cdots - \mu_1\varepsilon_N. \] Note that the weight $(\lambda,\mu)$ is always dominant and integral and moreover every dominant integral weight is of this form. Again this gives us a representation $W_{(\lambda,\mu)}$ of $\mathrm{GL}(N)$. The grading is given by $|\lambda| - |\mu|$.

We now turn to the classification of representations of $\mathrm{SL}(N)$. We can restrict a representation $W_\lambda$ or $W_{\lambda,\mu}$ to $\mathrm{SL}(N)$ where we will denote it by $V_\lambda$ or $V_{\lambda,\mu}$.

The Cartan subalgebra for $\mathrm{SL}(N)$ is the subalgebra of traceless diagonal matrices, and hence the space of weights for $\mathrm{SL}(N)$ is a \emph{quotient} of $\{\sum \lambda_i \varepsilon_i\}$ by vectors of the form $\sum c \varepsilon_i$. Again the dominant weights are the ones where the $\lambda_i$ are weakly decreasing, and the integral ones are the ones that have a representative with all $\lambda_i \in \mathbb{Z}$. In particular, if we choose a representative where the last entry is $0$, we see that dominant integral weights for $\mathrm{SL}(N)$ can be classified by Young diagrams $\lambda$ with $\ell(\lambda) < N$.

We now look at the fusion graph, which is an oriented graph with a vertex for every irreducible representation, and an oriented edge from $W$ to $W'$ for each $\dim \Hom(V_{\ydiagram{1}}\otimes W, W')$.

\begin{prop}
    The fusion graph of $\operatorname{Rep}^{\mathrm{poly}}(\mathfrak{gl}(N))$ has vertices indexed by Young diagrams with at most $N$ rows, and an oriented edge from $\lambda$ to $\lambda'$ if you can add a single box to $\lambda$ to turn it into $\lambda'$.
\end{prop}

As $N\rightarrow \infty$ this becomes the Young graph, which has a vertex for every Young diagram and an oriented edge for each way of adding a box.

\begin{prop} \label{prop:GLFusion}
    The fusion graph of $\operatorname{Rep}(\mathfrak{gl}(N))$ has vertices indexed by pairs of Young diagrams $(\lambda,\mu)$ with $\ell(\lambda)+\ell(\mu) \leq N$, and an oriented edge for each way of adding a square to the first diagram or removing a square from the second diagram.
\end{prop}

This graph has a limit as $N \rightarrow \infty$ which we call the double Young graph, which drops the condition on the total length.

\begin{prop}
    The simple objects of $\operatorname{Rep}(\mathfrak{sl}(N))$ are indexed by Young diagrams with at most $N-1$ rows. The fusion graph for tensoring with the defining representation is the oriented graph which has an edge for each way of adding a single box or removing a first column with $N-1$ boxes. 
\end{prop}

We will also require a few lemmas about counting (oriented) paths on these fusion graphs. The first lemma is very well-known.

\begin{defn}
    If $\lambda$ is a Young diagram, a standard tableaux of shape $\lambda$ with entries in $\{1, \ldots, n\}$ is a way of writing entries from $\{1, \ldots, n\}$ into each box in such a way that it is strictly increasing in rows and columns. If we do not specify $n$ then it's assumed that $n = |\lambda|$.
\end{defn}

\begin{lem} \label{lem:easypathstotableaux}
    The number of paths from $(\emptyset)$ to $(\lambda)$ on the Young graph equals the number of standard tableaux of shape $\lambda$.
\end{lem}
\begin{proof}
    The bijection assigns to each tableaux the walk whose $k$th vertex is the diagram consisting of the boxes labeled by $\{1, \ldots, k\}$.
\end{proof}

We now define skew Young diagrams. These combinatorial objects will play a role in Section~\ref{sec:con}. \begin{defn}
    We say that $\mu \subseteq \lambda$ if $\mu_i \leq \lambda_i$ for all $i$. In the special case that $\lambda$ is obtained from $\mu$ by adding one box, we will write $\mu \rightarrow \lambda$.
    If $\mu \subseteq \lambda$, the diagram for $\mu$ sits entirely inside the diagram for $\lambda$ and we define the skew Young diagram of shape $\lambda/\mu$ to consist of the boxes in $\lambda$ which are not in $\mu$. We define a standard tableau for the pair $\mu \subseteq \lambda$ with entries in $\{1, \ldots, n\}$, to be all ways of entering $\{1, \ldots, n\}$ into the boxes of $\lambda$ such that it is strictly increasing in rows and columns inside $\mu$, and also strictly increasing in rows and columns in $\lambda/\mu$. Again, unless otherwise noted $n = |\lambda|$
\end{defn} 
\ytableausetup{boxsize=normal}
\begin{ex}\label{ex:tableau}
    If $\lambda = (2,1)$ and $\mu = (2)$ then there are three standard tableau for the pair $\mu \subseteq \lambda$, namely:
    \[\begin{ytableau}
*(lightgray) 1& *(lightgray) 2  \\
3
\end{ytableau}\;, \qquad 
\begin{ytableau}
*(lightgray) 1& *(lightgray) 3  \\
2
\end{ytableau}\;, \qquad 
\begin{ytableau}
*(lightgray) 2& *(lightgray) 3  \\
1
\end{ytableau}\;.\]

Similarly, if $\lambda = (2,1)$ and $\mu = (1,1)$ then there are three standard tableau for the pair $\mu \subseteq \lambda$, namely:
    \[\begin{ytableau}
*(lightgray) 1&  3  \\
*(lightgray) 2
\end{ytableau}\;, \qquad 
\begin{ytableau}
*(lightgray) 1& 2  \\
*(lightgray) 3
\end{ytableau}\;, \qquad 
\begin{ytableau}
*(lightgray) 2& 1  \\
*(lightgray) 3
\end{ytableau}\;.\]
\end{ex}

The following counting lemmas allow us to relate paths on the Young graph with paths on the double Young graph.
\begin{lem} \label{lem:pathstotableaux}
    Unless $\mu \subseteq \lambda$ there are no paths $n$ from $(\mu, \mu) \to (\lambda,\emptyset)$ (or from $(\mu, \mu^T) \to (\lambda,\emptyset)$) on the double Young graph.
    If $\mu \subseteq \lambda$ then the number of paths from $(\mu, \mu) \to (\lambda,\emptyset)$ (or from $(\mu, \mu^T) \to (\lambda,\emptyset)$) on the double Young graph equals the number of standard tableaux for the pair $\mu \subseteq \lambda$.
\end{lem}
\begin{proof}
    The double Young graph has an automorphism given by taking transpose in the second coordinate, which induces a bijection between paths from  $(\mu, \mu) \to (\lambda,\emptyset)$ and paths from $(\mu, \mu^T) \to (\lambda,\emptyset)$ given by reflecting each diagram that appears in the second coordinate. So we only consider the former case.

    To have a path from $(\mu, \mu) \to (\lambda,\emptyset)$ of length $n$, we must remove $|\mu|$ boxes from $\mu$ to obtain $\emptyset$, and add $n - |\mu|$ boxes to $\mu$ to obtain $\lambda$. We have a bijection between such paths and standard tableaux for the pair $\mu \subseteq \lambda$ which associates to each tableau the path whose $k$th step is the pair whose first coordinate is $\mu$ together with the $\{1, \ldots, k\}$ entries in $\lambda-\mu$ and whose second coordinate is the $\{1, \ldots, k\}$ entries of $\mu$ (in particular, if $k$ lies in $\mu$ we are removing a box from the second coordinate, while if $k$ lies in $\lambda-\mu$ we are adding a box to the first coordinate). 
\end{proof}

\begin{lem} \label{lem:countingskewtableaux}
    Fixing $\lambda$ and a number $k\leq|\lambda|$, the following sets are in bijection:
    \begin{enumerate}
    \item pairs of $\mu \subseteq \lambda$ with $|\mu|=k$ and a standard tableau for $\mu \subseteq \lambda$
    \item pairs of a subset of $\{1, \ldots |\lambda|\}$ of size $k$ and standard tableaux on $\lambda$.
    \end{enumerate}
\end{lem}
\begin{proof}
    We construct a bijection between the two given sets. In order to do this we will need an auxilliary function which reorders the alphabet using the subset $X$.
    
    Let $n = |\lambda|$, define for $X\subseteq\{1,\cdots , n\}$ the unique bijection $f_X : \{1,\cdots , n\}\to \{1,\cdots , n\}$ such that 
    \begin{enumerate}
        \item $f_X$ sends $X$ to $\{1,\cdots , |X|\}$
        \item $f_X$ restricted to $X$ preserves order
        \item $f_X$ restricted to the complement of $X$ preserves order.
    \end{enumerate}

  Let $T$ be a standard tableau for the pair $\mu \subseteq \lambda$. Let $X$ be the subset $\{  T_{i,j} : (i,j) \in \mu\}\subseteq \{1,\cdots, n\}$ consisting of the numbers used to label boxes in $\mu$. We claim that $f_X(T)$, the tableau whose entries are renamed using $f_X$, is a standard $\lambda$ Young tableau. As $f_X$ is order preserving on the entries appearing in both $\mu$ and $\lambda - \mu$, we have that this is still a standard tableau for the pair $\mu \subseteq \lambda$. Further, as every element of $\im(f_X|_X)$ is smaller than every element of $\im(f_X|_{X^c})$ we get that every entry $\mu$ subdiagram is smaller than every entry of the $\lambda/\mu$ subdiagram. Thus $f_X(T)$ is a standard $\lambda$ Young tableau.

  In the reverse direction, let $T$ be a standard $\lambda$ Young tableau, and $X$ a subset $X$ of $\{1,\cdots, n\}$. We define $\mu := \{ (i,j) \in \lambda_1 :  T_{i,j}\in \{1,\cdots , |X|\}\}$ which is a Young diagram as $T$ is standard. From the ordering conditions on $f|_X$ and $f|_{X^c}$ we get that $f_X^{-1}(T)$ is a standard tableau for the pair $\mu \subseteq \lambda$.

  These two maps are clearly inverse to each other, hence we have the above claimed bijection.  
\end{proof}

\begin{ex}
    If $\lambda = (2,1)$ and $k=2$, then the first set in Lemma \ref{lem:countingskewtableaux} consists of the six tableau in Example \ref{ex:tableau}. Under the bijection constructed in the above proof these are sent to
\[\left(\{1,2\},\;\begin{ytableau}
 1&  2  \\
 3
\end{ytableau}\right)\;, \qquad \left(\{1,3\},\;\begin{ytableau}
 1&  2  \\
 3
\end{ytableau}\right)\;, \qquad 
\left(\{2,3\},\;\begin{ytableau}
 1&  2  \\
 3
\end{ytableau}\right)\;,\] 
\[
\left(\{1,2\},\;\begin{ytableau}
 1&  3  \\
 2
\end{ytableau}\right)\;, \qquad 
\left(\{1,3\},\;\begin{ytableau}
 1&  3  \\
 2
\end{ytableau}\right)\;, \qquad
\left(\{2,3\},\;\begin{ytableau}
 1&  3  \\
 2
\end{ytableau}\right)\; .\] 
\end{ex}
\ytableausetup{boxsize = .3em}

The relevance of this subsection to the results of this paper is that the combinatorics of $\mathfrak{sl}_N$ level $L$ are a truncation of the $\mathfrak{sl}_N$ combinatorics.

\begin{lem} \label{lem:stabilityofhom}
    The simple objects of $\overline{\operatorname{Rep}(U_q(\mathfrak{sl}_N))}$ at $q = e^{2\pi i \frac{1}{2(N+L)}}$ are indexed by Young diagrams with $\ell(\lambda) \leq N$ and $\lambda_1 \leq L$. Furthermore, the fusion graph for $\overline{\operatorname{Rep}(U_q(\mathfrak{sl}_N))}$ is the subgraph of the Young graph on these diagrams. 
    
    Moreover, if $\lambda_1 + \mu_1 \leq L$ then $V_\lambda \otimes V_\mu$ has the same decomposition of simples in $\overline{\operatorname{Rep}(U_q(\mathfrak{sl}_N))}$ that it does in $\operatorname{Rep}(\mathfrak{sl}(N))$.
\end{lem}
\begin{proof}
    See \cite{MR1328736} and \cite[\textsection 5]{sawin}.
\end{proof}

\subsection{HOMFLYPT Skein Categories}\label{sec:Hecke}

In this subsection we review the theory of the \textit{HOMFLYPT skein categories} These are a two parameter family of tensor categories that were introduced by Turaev \cite{turaev}. The categories we define in Definition~\ref{def:main} will be extensions of certain HOMFLYPT skein categories, and hence the general theory of the HOMFLYPT skein category will be important for this paper. HOMFLYPT skein categories can also be thought of as a Deligne interpolation category for $\overline{\operatorname{Rep}(U_q(\mathfrak{gl}_N))}$. 

\begin{defn}
    Let $q,r \in \mathbb{C}-\{0\}$. 
    The HOMFLYPT skein category $\mathcal{H}(q,r)$ is 
    the pivotal $\mathbb{C}$-linear monoidal category with objects strings in $\{+,-\}$ and morphisms generated by the two morphisms $\att{braidX}{.2}$ and $\att{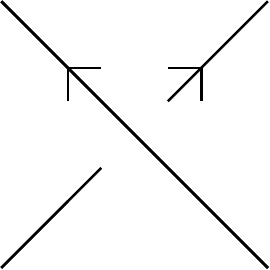}{.2}$ satisfying the relations:
    \begin{alignat*}{3}
\att{loop}{.15} &= \frac{r-r^{-1}}{q-q^{-1}} \qquad&  \att{R1}{.15}&= r\att{R12}{.15} & \att{R21}{.15} &= \att{R22}{.15}
\\
\att{R31}{.10} &= \att{R32}{.10}\qquad &  \att{H11}{.15} - \att{H11inv}{.15} &= (q-q^{-1})\att{H12}{.15}
    \end{alignat*}
\end{defn}
As in Definition~\ref{def:main}, the relations involving unoriented strands are understood to hold for all valid orientations. 

The HOMFLYPT skein category is closely related to quantum versions of $\mathfrak{gl}_N$, but it is easier to make this statement precise by looking at $\overline{\operatorname{Rep}(U_q(\mathfrak{sl}_N))}$. As in \cite{SovietHans} we can define an extension of $\mathcal{H}(q,r)$ when $r$ is of the form $q^N$ for $N\in \mathbb{N}_{\geq 2}$.
\begin{defn}[\cite{SovietHans}]\label{def:SL}
    Let $q \in\mathbb{C}- \{0\}$, and $N\in \mathbb{N}_{\geq 2}$. We define $\mathcal{SH}(q,N)$ as the extension of $\mathcal{H}(q,q^N)$ by the additional generator\edit{s}
\[    \att{kw}{.3}\in \Hom( +^N \to \mathbf{1}),\edit{ \qquad \att{kwdagger}{.3} \in \Hom( \mathbf{1} \to  +^N ) } \]
satisfying the relations
\[   \att{kwrel1}{.25} = q \att{kwrel12}{.25}\qquad \text{ and }\qquad    \att{kwrel2}{.25} =  \att{kwrel22}{.25}, \]
where $p_{1^N}$ is the unique minimal projection in the HOMFLYPT skein subcategory onto the simple $V_{1^N}$. 
\end{defn}
An explicit formula for $p_{1^N}$ in terms of the braid can be found in \cite[Theorem 6]{Paggo}. 

There are two other relations one often sees included \cite{dan,MR2171796,SovietHans}, but they follow automatically.

\begin{lem}
    The following relations hold in $\mathcal{SH}(q,N)$:
    
    \[   \att{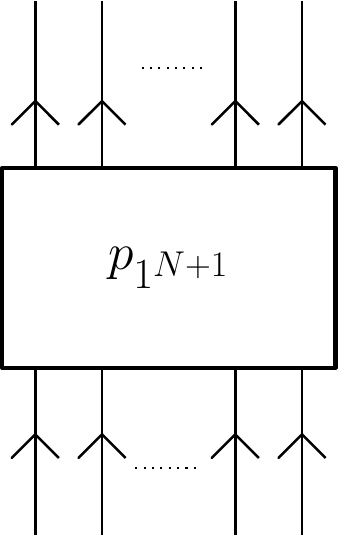}{.25}\quad = \quad 0 \qquad\qquad  \att{kwrel1}{.25}  \quad = \quad q^{-1} \att{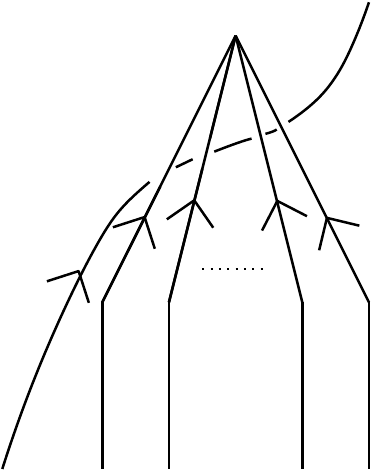}{.25} .  \]

\end{lem}
\begin{proof}
For the first relation we have
   \[  \att{kn0}{.3}\quad = \quad \att{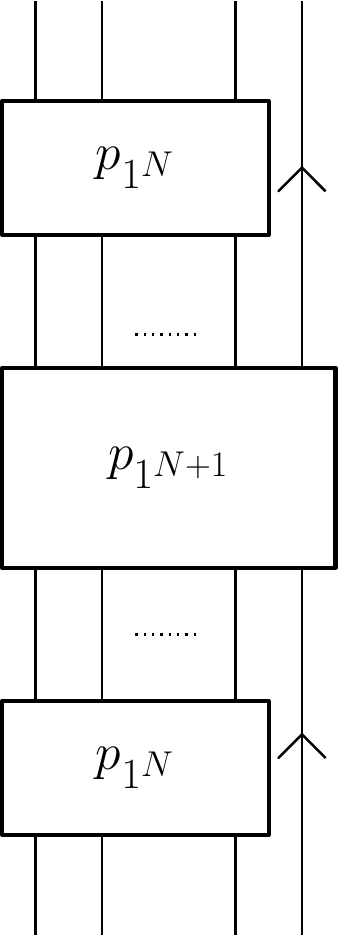}{.3}\quad = \quad \att{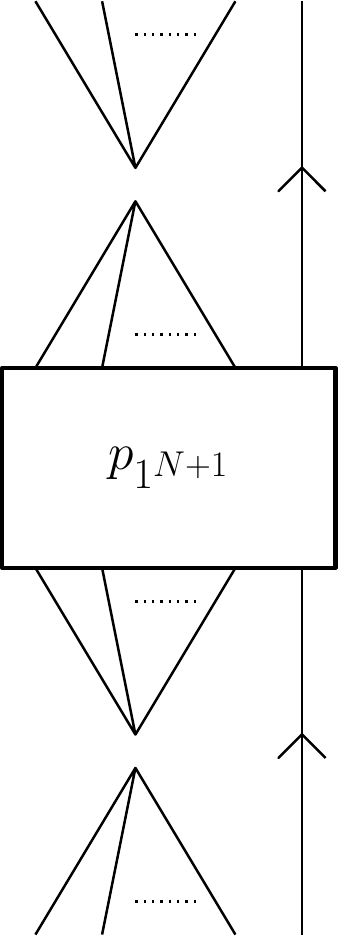}{.3}\quad = \quad \att{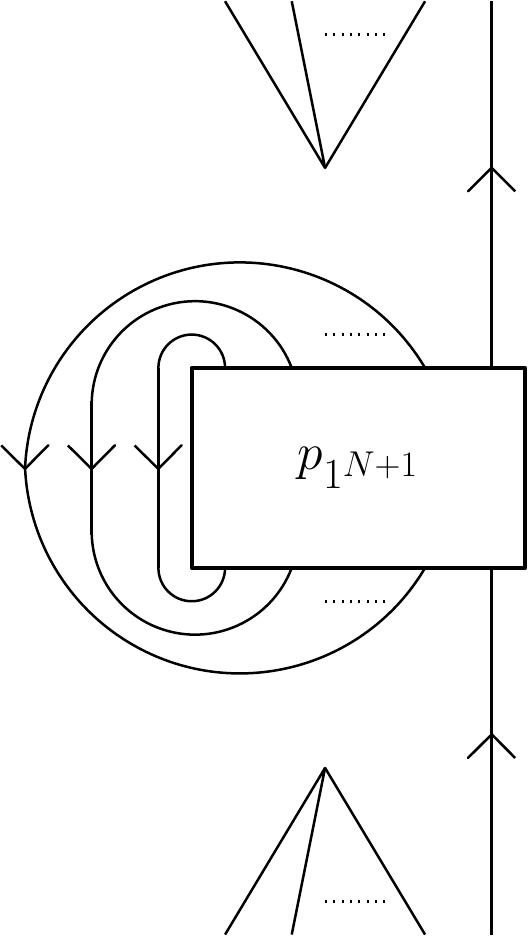}{.3} = 0.  \]
   Here the first equality holds as the projection $p_{1^{N+1}}$ absorbs the projection $p_{1^{N}}$ in all positions \cite[Theorem 5]{Paggo}, and the final equality follows as the quantum trace of $p_{1^{N+1}}$ vanishes \cite[Remark 2.12]{dan}.

   As a consequence of the $p_{1^{N+1}}=0$, it follows from \cite[Corollary 2.1]{Clam} that
   \[ \att{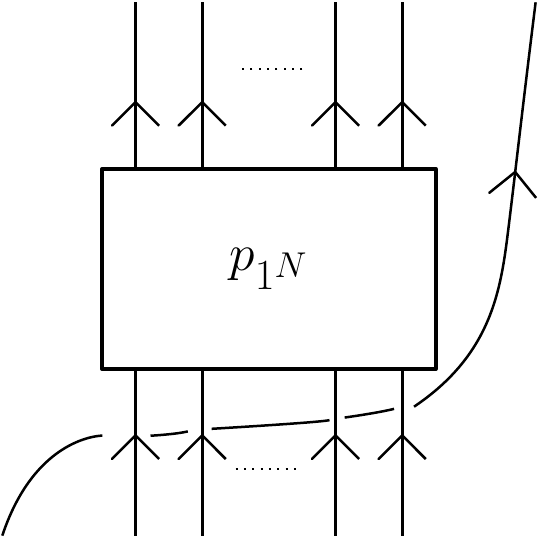}{.25} \quad = \quad q^{-2}\att{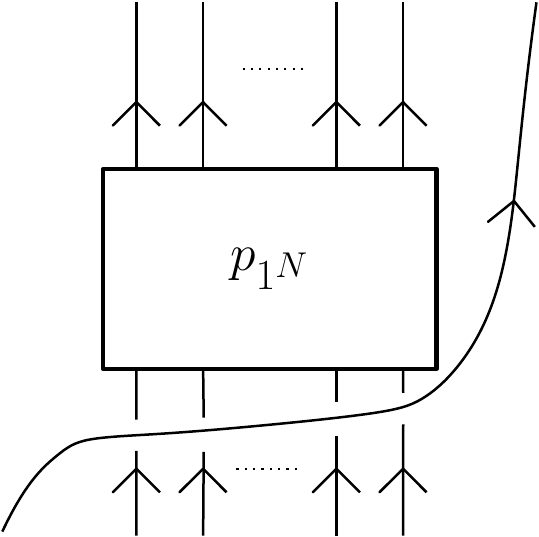}{.25}    \]
  The second relation then follows from the first relation defining relation of $\mathcal{SH}(q, N)$ along with the equality
   \[ \att{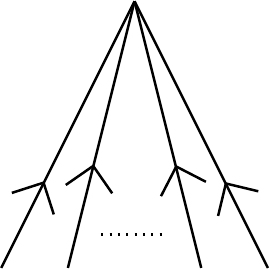}{.25} \quad = \quad \att{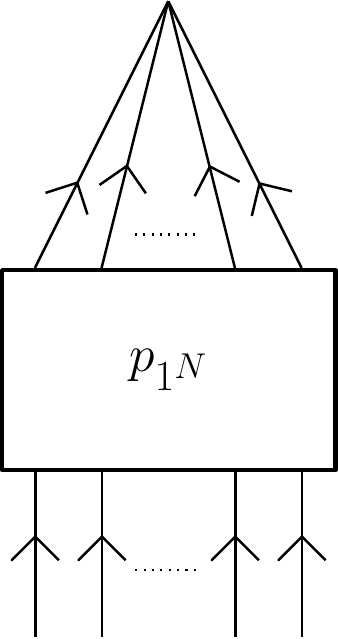}{.25}\]
   which is an immediate consequence of the second defining relation.
\end{proof}

Note that the crossing is not a braiding for $\mathcal{SH}(q,N)$, because it is not natural with respect to the new generator. Instead we have the following.

\begin{lem}\label{lem:fudgefactor} $\mathcal{SH}(q,N)$ is a braided tensor category with braiding $q^{-1/N}$ times the crossing.
\end{lem}

\begin{remark}
    One can incorporate the rescaling in Lemma \ref{lem:fudgefactor} directly into the diagram calculus following \cite{MR2171796}. We choose not to do this, because the benefits are small since the crossings in $\mathcal{E}_R$ are not braidings, and the book-keeping costs of using a modified presentation for the Hecke algebra in this paper (and especially in \cite{HansNew}) are significant.
\end{remark}

The connection between the HOMFLYPT category and the quantum $\mathfrak{sl}_N$ category can now be made precise.

\begin{prop}\cite{SovietHans}\label{prop:sen}
    There is a braided equivalence $\overline{\operatorname{Ab}{\mathcal{SH}(q,q^N)}} \rightarrow \overline{\operatorname{Rep}(U_q(\mathfrak{sl}_{N}))}$ via the unique braided tensor functor sending the strand to the defining representation $V$ and the additional generator to the quantum determinant map.
\end{prop}
\begin{proof}
    By \cite[Prop. 2.2(b)]{SovietHans}, it is enough to check that the functor induces an isomorphism on $\End(+^n) \rightarrow \End(V^{\otimes n})$, but both $\End(+^n)$ and  $\End(V^{\otimes n})$ are the quotient of the Hecke algebra described in \cite[Prop. 3.1]{SovietHans} and the map between them is the identity.
\end{proof}

\subsection{The Hecke-Clifford algebras}\label{Gndef} 

The Hecke-Clifford algebras were defined in \cite{Ol} in the context of centralizer algebras of quantum group versions of isomeric Lie super algebras.

\begin{defn} We define the Hecke-Clifford algebras $G_n$ as the associative $R$-algebras with generators $t_j$, $1\leq j<n$ and $v_j$, $1\leq j\leq n$ and the following relations:

\begin{itemize}
\item[(H)] The generators $t_j$ satisfy the relations of the Hecke algebras $H_n$ of type $A_{n-1}$. This means they satisfy
 the braid relations as well as the quadratic equation $t_j-t_j^{-1}=q-q^{-1}$.

\item[(C)] The elements $v_j$ generate the Clifford algebra $\mathrm{Cliff}(n)$
with relations $v_jv_k+v_kv_j=2\delta_{jk}$. 

\item[(M)] Moreover, we have the additional relations
\[   t_jv_j=v_{j+1}t_j,\qquad  t_jv_{j+1}= v_jt_j-(q-q^{-1})(v_j-v_{j+1}),\qquad t_jv_l=v_lt_j, \hskip 3em l\neq j, j+1.  \]
    
\end{itemize}
\end{defn}

\begin{thm}\label{Gnproperties} The algebra $G_n$ is equal to $H_n\mathrm{Cliff}(n)$ as a vector space and has dimension $2^n n!$. In particular, $G_n$ has a standard basis of the form $h_w c_s$ where $w$ ranges over $S_n$ and $s$ ranges over $\{0,1\}^n$. Similarly, $G_n = \mathrm{Cliff}(n)H_n$ and has a basis of the form $c_s h_w$.
\end{thm}

Observe that the algebra $G_n$ has a $\Z_2$-grading, where the elements $t_j$ have degree 0 and 
the elements $v_j$ have degree 1. It follows from the relations that the map
\begin{equation*}\label{autom}
\al:\quad v_j\mapsto -v_j,\hskip 3em t_j\mapsto t_j,
\end{equation*}
defines an automorphism of order two with eigenspaces $G_n[0]$ and $G_n[1]$.

We will also require generators and relations for $G_n[0]$, which requires a change of variables since the $v_j$ are odd. So we set $e_j=v_jv_{j+1}$, and note that the $e_j \in \Gno$.

\begin{lem} \label{Gn0Presentation}
    The following relations hold in $\Gno$:
    \begin{itemize}
        \item[(E)] The $e$'s satisfy the relations $e_ie_j=e_je_i$ for $|i-j|\neq 1$, $e_j^2=-1$ and $e_je_{j+1}=-e_{j+1}e_j$, $1\leq j<n-1$.
        \item[(M)] We have the mixed relations

        \[  e_jt_j+t_j^{-1}e_j=(q-q^{-1})1, \qquad  t_je_{j+1}=-e_je_{j+1}t_j, \qquad e_jt_{j+1}=-t_{j+1}e_je_{j+1}, \qquad e_jt_{j+1}t_j=t_{j+1}t_je_{j+1}. \]

    \end{itemize}

Moreover, $\Gno$ has a basis of the form $h_w e_s$ where $w$ ranges over $S_n$ and $s$ ranges over elements of $\{0,1\}^{n-1}$ where we again multiply the $e_i$ in lexicographical order.
\end{lem}
\begin{proof}
It is straightforward if somewhat tedious to check that the relations are satisfied. We check the first (M) relation as an example
\begin{align*}
e_jt_j + t_j^{-1}e_j &= v_jv_{j+1}t_j + t_{j}^{-1} v_j v_{j+1}\\
&= v_jt_jv_j + t_jv_jv_{j+1} - (q-q^{-1})v_jv_{j+1}\\
&= t_jv_{j+1}v_j + (q-q^{-1})(v_j - v_{j+1})v_j + t_jv_jv_{j+1} - (q-q^{-1})v_j v_{j+1}\\
&= (q-q^{-1})1.
\end{align*}

Clearly $G_n[0]$ is a product of the Hecke algebra and the even part of the Clifford algebra. We now show that the even part of the Clifford algebra has a basis given by products of $v_jv_{j+1}$ in lexicographic order, by induction on $n$. Consider a basis vector $\prod_{i=1}^k v_{r_i}$ with $k$ even, and the $r_i$ in increasing order. If $r_k \neq n$ then we are done by induction, and otherwise we have \[\prod_{i=1}^k v_{r_i} = (-1)^{r_k-r_{k-1}+1} \prod_{i=1}^{k-2} v_{r_i} \prod_{i=r_{k-1}}^{r_k} v_i v_{i+1},\] and again we're done by induction.
\end{proof}

These relations will appear again in the description of $\End_{\E_{R}}(+^n)$, as given in Lemma \ref{lem:relationspic}.

\section{Conformal Embeddings}\label{sec:con}

A large class of \textit{\'etale} algebra objects in the categories $\overline{\operatorname{Rep}(U_q(\mathfrak{sl}_N))}$ come from conformal embeddings of WZW vertex operator algebras. These \textit{\'etale} algebras first appeared in \cite{Xu} in subfactor language, and appears in \cite{Kril} and \cite{HKL} in tensor category language.

Letting $\mathcal{V}(\mathfrak{g}, k)$ denote the WZW vertex operator algebra for $\mathfrak{g}$ at level $k$ \cite{WZW}, we consider inclusions of the form
\[       \mathcal{V}(\mathfrak{sl}_N, k)\subseteq \mathcal{V}(\mathfrak{g}, 1) ,   \]
where $\mathcal{V}(\mathfrak{g}, 1)$ decomposes as a finite direct sum of irreducible $\mathcal{V}(\mathfrak{sl}_N, k)$-modules. A complete list of such inclusions can be found in \cite{LagrangeUkraine}. These inclusions were first found in \cite{Embeddings1,Embeddings2}. It follows from \cite[Theorem 5.2]{Kril} or \cite[Theorem 3.2]{HKL} that $\mathcal{V}(\mathfrak{g}, 1)$ has the structure of an \textit{\'etale} algebra object in $\operatorname{Rep}(\mathcal{V}(\mathfrak{sl}_N, k))$. The modular tensor category $\operatorname{Rep}(\mathcal{V}(\mathfrak{sl}_N, k))$ is naturally identified with the category of level $k$ integrable representations of $\hat{\mathfrak{sl}_N}$, which is braided equivalent to $\overline{\operatorname{Rep}(U_q(\mathfrak{sl}_N))}$ at $q=e^{2\pi i \frac{1}{2(N+k)}}$ due to work of Kazhdan-Lusztig \cite{KL1, KL3, KL4} and Finkelberg \cite{MR1384612}. This result also follows by combining works of Kazhdan-Wenzl \cite{SovietHans} and Wassermann \cite{Wasser}. We thus get \textit{\'etale} algebra objects in $\overline{\operatorname{Rep}(U_q(\mathfrak{sl}_N))}$ due to the above construction. We direct the reader to \cite[Section 6.2]{LagrangeUkraine} for a more detailed treatment.

For this article, we are interested in the following specific examples of \textit{\'etale} algebra objects coming from conformal embeddings.
\begin{defn}
    For each $N\geq 2$ we define $A$ to be the \textit{\'etale} algebra object in $\overline{\operatorname{Rep}(U_q(\mathfrak{sl}_N))}$ at $q = e^{2\pi i \frac{1}{4N}}$ corresponding to the conformal embedding 
    \[         \mathcal{V}(\mathfrak{sl}_N, N)\subseteq \mathcal{V}(\mathfrak{so}_{N^2-1}, 1)    \]       
from \cite[p. 34]{LagrangeUkraine}
    under the above equivalences.
\end{defn}

One of the main goals of this paper is to study the category $\overline{\operatorname{Rep}(U_q(\mathfrak{sl}_N))}_A$. We will use the unitarity of this category at several points throughout this paper.

\begin{prop}\label{rmk:unitary}
The category $\overline{\operatorname{Rep}(U_q(\mathfrak{sl}_N))}_A$ is a unitary fusion category.
\end{prop}
\begin{proof}
     By \cite{MR4616673} every connected rigid algebra object in a unitary fusion category is a $Q$-system.  As $\overline{\operatorname{Rep}(U_q(\mathfrak{sl}_N))}$ is unitary at $q = e^{2\pi i \frac{1}{4N}}$ \cite{JamsHans} it follows that $\overline{\operatorname{Rep}(U_q(\mathfrak{sl}_N))}_A$ is a unitary fusion category.
\end{proof}

We begin by finding two distinguished objects. The first we obtain from the free module functor.

\begin{defn}
    We define $X:= \mathcal{F}_A(V_\square)\in \overline{\operatorname{Rep}(U_q(\mathfrak{sl}_N))}_A$. 
\end{defn}

Our second distinguished object comes from the subcategory of local modules $\overline{\operatorname{Rep}(U_q(\mathfrak{sl}_N))}_A^0$ (see \cite[Section 3.5]{LagrangeUkraine}). For \textit{\'etale} algebra objects corresponding to embeddings of WZW vertex operator algebras, the subcategory of local modules is fully understood. From \cite[Theorem 5.2]{Kril} we have that 
\begin{equation}\label{eq:loc}
    \overline{\operatorname{Rep}(U_q(\mathfrak{sl}_N))}_A^0 \simeq  \operatorname{Rep}(\mathcal{V}(\mathfrak{so}_{N^2-1}, 1)).
\end{equation} 

We start with some basic facts about the $\operatorname{Rep}(\mathcal{V}(\mathfrak{so}_{N^2-1}, 1))$. The following Lemma is well-known.
\begin{lem} \label{lem:LocalDescription}
    The vector representation $G$ in $\operatorname{Rep}(\mathcal{V}(\mathfrak{so}_{N^2-1}, 1))$ tensor generates a ribbon subcategory equivalent to $\mathrm{SVec}$ where the odd line is given dimension $1$.
\end{lem}
\begin{proof}
    That $G \otimes G \cong 1$ appears in \cite[Eqns. 1.6-1.8]{MR1167302}. We want to check that the $\mathrm{qdim}(G)=1$ and that the twist on $G$ is $-1$, from which it follows that the crossing is negative the identity.
    
    The twist for the irrep $V_\omega$ in $\operatorname{Rep}(\mathcal{V}(\mathfrak{so}_{m}, k))$ is $(q')^{\langle \omega, \omega + 2 \rho\rangle} = (q')^{m-1}$ where $q' = e^{\frac{ \pi \bf{i}}{k+h^\vee}}$, $\rho$ is the Weyl vector, and $h^\vee$ is the dual Coxeter number. So setting $k=1$, $\omega$ the weight of the defining representation, and $h^\vee = m-2$, we get twist $-1$. Similarly, the dimension of the defining representation is $[m-1]_{q'}+1$, and again specializing $q'$ we get dimension $0+1 = 1$.
\end{proof}

\begin{remark}
Note that the category $\operatorname{Rep}(\mathcal{V}(\mathfrak{so}_{N^2-1}, 1))$ is larger than $\mathrm{SVec}$ because it also contains spin representations, again see \cite[Eqns. 1.6-1.8]{MR1167302}.  When $N$ is odd, this category is pointed with underlying group $\mathbb{Z}_2\times\mathbb{Z}_2$. When $N$ is even, then this category is equivalent to an Ising category \cite[Appendix B]{braidedI}. We will not need information regarding these additional spin objects, as the construction of our interpolation category does not use these objects directly. By the main results of this paper, these spin objects will appear as images of idempotents in the category $\mathcal{SE}_N$, but these idempotents lie in a large hom space with the number of strands depending on $N$.
\end{remark}

\begin{defn}
    We define $g$ to be the image of $G$ in $ \overline{\operatorname{Rep}(U_q(\mathfrak{sl}_N))}_A^0$ under the equivalence from Equation~\ref{eq:loc}.
\end{defn}

As $g$ is an object of $\overline{\operatorname{Rep}(U_q(\mathfrak{sl}_N))}_A$, we can apply the forgetful functor to obtain the underlying object in $\overline{\operatorname{Rep}(U_q(\mathfrak{sl}_N))}$. We will denote this underlying object as $B$. 

\begin{remark}
Since $\mathbf{1} \oplus G$ is a commutative super-algebra object, we also have that $A\oplus B$ has the structure of a commutative super-algebra. This algebra is the Free fermion algebra $\mathrm{Fer}(\mathfrak{g})$ \cite[Subsection 3.1]{papi} which makes $A \oplus B$ especially well-behaved.
\end{remark}

There has been a significant body of work dedicated to determining the forgetful functors for various classes of conformal embeddings. This allows us to give an explicit combinatorial description of the objects $A$ and $B$. This result is likely known to experts, however we could not find a proof nor statement in the literature. Here we use the $(\lambda, \mu)$ description of labeling $\mathfrak{sl}_N$ weights introduced in Subsection~\ref{sec:GLNcombo}.

\begin{thm}\label{thm:alg}
    As objects in $\overline{\operatorname{Rep}(U_q(\mathfrak{sl}_N))}$ we have that
    \[     A \cong \bigoplus_{\substack{ H_\mu(1,1) < N \\|\mu| \text{ even }}} V_{(\mu, \mu^T)}\qquad \text{and}\qquad  B \cong \bigoplus_{\substack{H_\mu(1,1) < N \\|\mu| \text{ odd }}} V_{(\mu, \mu^T)}   \]
    where both sums run over Young diagrams.
\end{thm}

\begin{proof}
From \cite[Theorem 3.9]{KacAlg} the simple summands of the object $A$ (resp.  $B$) are parameterised by  conjugacy classes of even (resp. odd) dimensional abelian subalgebras $\mathfrak{a}$ of a Borel subalgebra $\mathfrak{b}$ for $\mathfrak{sl}_N$. The bijection sends the subalgebra $\mathfrak{a}$ to the sum of positive roots appearing in the root space decomposition of $\mathfrak{a}$, which is a dominant weight.

Let $\alpha_i = \varepsilon_i - \varepsilon_{i+1}$. We make our choice of positive roots as
\[\Phi_+ := \{\alpha_1, \alpha_2, \cdots, \alpha_1+\alpha_2, \alpha_2+\alpha_3,\cdots, \alpha_1 +\cdots + \alpha_{N-1}\}.\]
As described in \cite{borel} abelian subalgebras of $\mathfrak{b}$ correspond to subsets $\Psi \subseteq \Phi_+$ such that 
\begin{enumerate}[(a)]
    \item $(\Psi + \Phi_+)\cap \Phi_+ \subseteq \Psi$, and 
    \item $(\Psi + \Psi)\cap \Phi_+ = \emptyset$.
\end{enumerate}

We claim there is a bijection between such subsets, and Young diagrams $\mu$ with $H_\mu(1,1)< N $. This bijection is given by identifying a subset of $\Phi_+$ with a collection of boxes according to the following diagram:
\[ \att{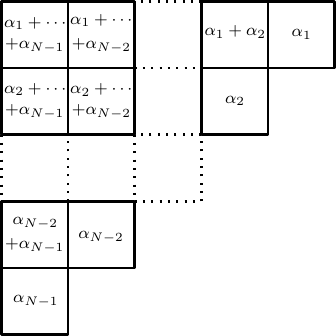}{1.0}\]
Condition (a) is equivalent to the collection being a Young diagram. Condition (b) is equivalent to having the smallest rectangle containing the Young diagram being contained within the full staircase. This is equivalent to the Young diagram having a $(1,1)$ hook length strictly less than $N$. Hence subsets of $\Phi_+$ satisfying conditions (a) and (b) are in bijection with Young diagrams $\mu$ with $H_\mu(1,1) < N$. The dimension of the sub-algebra corresponding to a Young diagram $\mu$ under this bijection is the number of positive roots appearing in $\Psi_{\mu}$, which is clearly seen to be $|\mu|$. Thus summands of the object $A$ correspond to even Young diagrams, and summands of the object $B$ correspond to odd Young diagrams.

Finally we determine the dominant weight corresponding to the Young diagram $\mu$. As described earlier, this dominant weight is given by the sum of the positive roots appearing in $\Psi_{\mu}$. Under our bijection this is given by
\[      \sum_{i=1}^n \mu_i \times \varepsilon_i  - \varepsilon_{n+1-\mu_i}-\cdots  - \varepsilon_{n}.  \]
As $\ell(\mu)   + \ell(\mu^T) =  H_\mu(1,1)+1 \leq N$, the above weight is exactly the dominant weight corresponding to the pair $(\mu, \mu^T)$ under the correspondence described in Section \ref{sec:GLNcombo}.
\end{proof}
This bijection allows a simple determination of the summands of $A$ and $B$ as illustrated in the following example.
\begin{ex}
    In the case of $N=4$ the even Young diagrams with $(1,1)$ hook length less than 4 are
    \[  \emptyset, \quad \ydiagram{2}, \quad \ydiagram{1,1}, \quad \ydiagram{2,2}. \]
    These correspond to the highest weights
    \[  (0,0,0,0), \quad (2,0,-1,-1), \quad (1,1,0,-2), \quad (2,2,-2,-2).       \]
    The odd Young diagrams with $(1,1)$ hook length less than 4 are
    \[  \ydiagram{1}, \quad \ydiagram{2,1}, \quad \ydiagram{3}, \quad \ydiagram{1,1,1}. \]
    These give us the highest weights
    \[   (1,0,0,-1), \quad (2,1,-1,-2), \quad (3,-1,-1,-1), \quad (1,1,1,-3).   \]
\end{ex}

The free module functor $\mathcal{F}_A:\overline{\operatorname{Rep}(U_q(\mathfrak{sl}_N))}\to \overline{\operatorname{Rep}(U_q(\mathfrak{sl}_N))}_A$ is adjoint to the forgetful functor $\operatorname{For}:\overline{\operatorname{Rep}(U_q(\mathfrak{sl}_N))}_A\to \overline{\operatorname{Rep}(U_q(\mathfrak{sl}_N))}$ \cite[Lemma 3.5]{Xu}. This adjunction, along with the explicit description of the objects $A$ and $B$, allows us to determine initial useful information regarding the categories $\overline{\operatorname{Rep}(U_q(\mathfrak{sl}_N))}_A$. The first result shows that $g$ is a summand of $X\otimes X^*$.

\begin{lem}\label{lem:simp}
Let $N\in \mathbb{N}_{\geq 2}$. Then  \[\dim \operatorname{Hom}_{\overline{\operatorname{Rep}(U_q(\mathfrak{sl}_N))}_A}(g \to X\otimes X^*) = 1.\]
\end{lem}
\begin{proof}
    Using the adjunction of $\mathcal{F}_A$ and $\operatorname{For}$, and that $X = \mathcal{F}_A(V_\square)$, we see
    \[  \dim \operatorname{Hom}_{\overline{\operatorname{Rep}(U_q(\mathfrak{sl}_N))}_A}(g \to X \otimes X^*)=  \dim \operatorname{Hom}_{\overline{\operatorname{Rep}(U_q(\mathfrak{sl}_N))}}(B \to V_\square  \otimes V_\square^*).   \]
    We have
    \[   V_\square  \otimes V_\square^*\cong  V_{(\emptyset, \emptyset)} \oplus V_{(\square, \square)}.   \]
    By our description of the summands of $B$, we have that $V_{(\emptyset, \emptyset)}$ has multiplicity zero in $B$, and the object $V_{(\square, \square)}$ has multiplicity one in $B$. Thus
    \[\dim \operatorname{Hom}_{\overline{\operatorname{Rep}(U_q(\mathfrak{sl}_N))}_A}(g \to X\otimes X^*) = 1.\]
 \end{proof}

The same style of argument also gives us the dimensions of the endomorphism algebras of $X^{\otimes n}$ in $\overline{\operatorname{Rep}(U_q(\mathfrak{sl}_N))}_A$. However, the combinatorics are much more involved for this argument. This result will be key later in the paper for producing a basis for these endomorphism algebras.

\begin{remark}
    The bound of $n < \frac{N}{2}$ in the following theorem can be improved to $n<N$. This is proven in \cite{HansNew}. For this paper we only require the weaker bound, which is significantly easier to obtain. 
\end{remark}

\begin{thm}\label{thm:dimEnd}
    Let $N\in \mathbb{N}_{\geq 2}$. Then
    \[ \dim\End_{\overline{\operatorname{Rep}(U_q(\mathfrak{sl}_N))}_A}(X^{\otimes n}) = n!\cdot 2^{n-1}  \]
    for all $n< \frac{N}{2}$.
\end{thm}
\begin{proof}
   As the free module functor $\mathcal{F}_A$ is adjoint to the forgetful functor $\operatorname{For}$, we get that \[    \dim\End_{\overline{\operatorname{Rep}( U_q(\mathfrak{sl}_N)   )}_A}(X^{\otimes n })=\dim\End_{\overline{\operatorname{Rep}( U_q(\mathfrak{sl}_N)   )}_A}(\mathcal{F}_A(V_\square)^{\otimes n }))  = \dim\Hom_{\overline{\operatorname{Rep}( U_q(\mathfrak{sl}_N)   )}}(A\otimes V_\square^{\otimes n }\to V_\square^{\otimes n } ) .  \]

  Supposing that $n < \frac{N}{2}$ by Lemma \ref{lem:stabilityofhom} and Frobenius reciprocity we get
  \[  \dim\Hom_{\overline{\operatorname{Rep}(U_q(\mathfrak{sl}_N)) } }(V_{\lambda}\otimes V_\square ^{\otimes n} \to V_\square ^{\otimes n})=\dim\Hom_{ \operatorname{Rep}(\mathfrak{sl}_N)}  (V_{\lambda}\otimes V_\square ^{\otimes n} \to V_\square ^{\otimes n})  \]
  for all $\lambda$. 
  
  Pick $W_{\lambda_1, \lambda_2}$ to be any object of $\operatorname{Rep}(\mathfrak{gl}_N)$ such that $\operatorname{Res}(W_{\lambda_1, \lambda_2}) = V_{\lambda}$, then by Lemma \ref{lem:GradingSLGL}

  \begin{align*}   \dim\Hom_{ \operatorname{Rep}(\mathfrak{sl}_N)}  (V_{\lambda}\otimes V_\square ^{\otimes n} \to V_\square ^{\otimes n}) 
  & =  \dim\Hom_{ \operatorname{Rep}(\mathrm{GL}_N)}  (W_{ \lambda_1,\lambda_2 }\otimes W_{\square,\emptyset}^{\otimes n} \to W_{\square,\emptyset}^{\otimes n} ).
  \end{align*}

From Theorem~\ref{thm:alg} we have that 
    \[ A = \bigoplus_{\substack{ H_\mu(1,1) <N\\ 2 \text{ divides }|\mu|}} 
   V_{(\mu, \mu^T)}. \]
   We thus have
  \[  \dim\End_{\overline{\operatorname{Rep}( U_q(\mathfrak{sl}_N)   )}_A}(X^{\otimes n}) = \dim\Hom_{\operatorname{Rep}(\mathfrak{gl}_N)}\left(  \left(\bigoplus_{\substack{ H_\mu(1,1) <N\\ 2 \text{ divides }|\mu|}} 
  W_{(\mu, \mu^T)}\right)\otimes W_{(\square, \emptyset)}^{\otimes n}\to W_{(\square, \emptyset)}^{\otimes n} \right). \]

 Hence the dimension of the above hom space is the number of pairs of paths $(\emptyset, \emptyset) \to (\lambda_1, \lambda_2) \leftarrow (\mu, \mu^T)$ of length $n$, running over all even $\mu$ with $H_\mu(1,1) <N$ on the fusion graph of $\mathrm{GL}(N)$ described in Proposition~\ref{prop:GLFusion}. For such a path to exist, we must have $\lambda_2 = \emptyset$, $|\lambda_1| = n$, and $|\mu| \leq n$. As $n < \frac{N}{2}$, and the hook length of a Young diagram is bounded above by its size, we have that $H_\mu(1,1) \leq n <\frac{N}{2} < N$ holds automatically. Hence the only restriction on $\mu$ is that $|\mu| \leq n$. Moreover, since $\ell(\lambda) \leq n < \frac{N}{2}$, we can instead count paths on the double Young graph.

  By Lemma \ref{lem:easypathstotableaux}, the number of paths from $(\emptyset, \emptyset) \to (\lambda_1,\emptyset)$ of length $n$ is the number of standard tableaux on $\lambda_1$. By Lemma \ref{lem:pathstotableaux} the number of paths from $(\mu, \mu^T) \to (\lambda_1,\emptyset)$ of length $n$ is the number of standard tableau for the pair $\mu \subseteq \lambda_1$. Running over all $\mu$ of fixed size $k$, by Lemma \ref{lem:countingskewtableaux} this counts the number of pairs of a subset of size $k$ and a standard tableau on $\lambda$. Thus the total number of paths running over all $\mu$ such that $|\mu| \leq n$ and $2$ divides $|\mu|$ from $(\mu, \mu^T) \to (\lambda_1,\emptyset)$ is $2^{n-1}$ (i.e the number of even sized subsets of $\{1, \ldots, n\}$ times the number of standard Young tableaux on $\lambda_1$. Therefore
  \[\dim\End_{\overline{\operatorname{Rep}(U_q(\mathfrak{sl}_N))}_A}(X^{\otimes n}) = \sum_{\lambda_1 : |\lambda_1| = n} 2^{n-1} \cdot |\{\text{Standard Young tableau on } \lambda_1\}|^2.\]
  As $\sum_{\lambda_1 : |\lambda_1| = n}  \cdot |\{\text{Standard Young tableau on } \lambda_1\}|^2$ is exactly the number of pairs of paths on the Young graph from $\emptyset \to \lambda_1$, running over all $\lambda_1$ with $|\lambda_1| = n$, we have that this quantity is $n!$ by Schur-Weyl duality. Hence $\dim\End_{\overline{\operatorname{Rep}( U_q(\mathfrak{sl}_N)   )}_A}(X^{\otimes n}) = 2^{n-1}\cdot n!$.
  \end{proof}

\section{Existence of the categories \texorpdfstring{$\mathcal{SE}_N$}{SE(N)}}\label{sec:sen}
In Definition~\ref{def:main2} we defined for all $N\in \mathbb{N}_{\geq 2}$ the category $\mathcal{SE}_N$ via generators and relations. \textit{A priori} these categories could be trivial. In this section we show the categories $\mathcal{SE}_N$ are non-zero, by proving that the semi-simplification of $\mathcal{SE}_N$ is a presentation for the category $\overline{\operatorname{Rep}(U_q(\mathfrak{sl}_N))}_A$. As the categories $\overline{\operatorname{Rep}(U_q(\mathfrak{sl}_N))}_A$ are manifestly non-zero, this gives non-triviality for the $\mathcal{SE}_N$ categories. More precisely, we will prove the following theorem.
\begin{thm}\label{thm:special}
    For all $N\in \mathbb{N}_{\geq 2}$ there exists a full and Cauchy-dominant functor 
    \[\Phi: \mathcal{SE}_N \to \overline{\operatorname{Rep}(U_q(\mathfrak{sl}_N))}_A.\]
    This functor descends to a fully faithful Cauchy-dominant functor 
    \[\overline{\Phi}: \overline{\mathcal{SE}_N} \to \overline{\operatorname{Rep}(U_q(\mathfrak{sl}_N))}_A,\]
    and hence a monoidal equivalence
    \[     \operatorname{Ab}(\overline{\mathcal{SE}_N}) \simeq \overline{\operatorname{Rep}(U_q(\mathfrak{sl}_N))}_A.   \]
\end{thm}

This section will be devoted to proving this theorem. We will do this in several parts. First we will show that $\overline{\operatorname{Rep}(U_q(\mathfrak{sl}_N))}_A$ contains morphisms satisfying the defining relations of the generators of $\mathcal{SE}_N$. This gives the existence of the functor $\Phi$, which is Cauchy-dominant by construction. We then show that $\mathcal{SE}_N$ has simple unit, which implies by general theory that $\Phi$ descends to faithful Cauchy-dominant functor $\overline{\Phi}$. Finally, we show $\overline{\Phi}$ is full by showing that it induces an equivalence between $\operatorname{Ab}(\overline{\mathcal{SE}_N})$ and $\overline{\operatorname{Rep}(U_q(\mathfrak{sl}_N))}_A$. As $\Phi$ is equal to $\overline{\Phi}$ composed with the full semisimplification functor $\mathcal{SE}_N \to \overline{\mathcal{SE}_N}$, we have that $\Phi$ is also full.

\begin{remark} \label{rem:Nis2}
    The case $N=2$ is highly degenerate with $A = \mathbf{1}$. Nonetheless, the statement is still true at $N=2$, see Remark \ref{rmk:N2}, and our argument still goes through.
\end{remark}

\subsection{Defining the Functor \texorpdfstring{$\Phi$}{Phi}}\label{sec:functor}

Our goal for this subsection is to find morphisms in $\overline{\operatorname{Rep}(U_q(\mathfrak{sl}_N))}_A$ with $q=e^{2\pi \textbf{i}/4N}$ which satisfy the defining relations of $\mathcal{SE}_N$. This will allow us to define the functor $\Phi: \mathcal{SE}_N\to \overline{\operatorname{Rep}(U_q(\mathfrak{sl}_N))}_A$ appearing in Theorem~\ref{thm:special}.

\begin{remark}
In order to apply diagrammatic techniques, we will work with a strictly pivotal model of $\overline{\operatorname{Rep}(U_q(\mathfrak{sl}_N))}$ for the remainder of the paper. This means that the double dual functor is the identity, and that the associator is trivial \cite[Definition 2.1]{Richard}. We can assume that $\overline{\operatorname{Rep}(U_q(\mathfrak{sl}_N))}$ is strictly pivotal by \cite[Theorem 2.2]{Richard}.  Alternatively, we could use a coherence theorem saying that planar string diagrams can be interpreted in any pivotal monoidal category \cite[Theorem 4.14]{Sel}.

As $\theta_A = \operatorname{id}_A$ \cite[Theorem 5.2]{Kril}, then \cite[Theorem 1.17]{Kril} implies that $\overline{\operatorname{Rep}(U_q(\mathfrak{sl}_N))}_A$ is also strictly pivotal, and that the free module functor $\mathcal{F}_A:  \overline{\operatorname{Rep}(U_q(\mathfrak{sl}_N))}\to \overline{\operatorname{Rep}(U_q(\mathfrak{sl}_N))}_A$ is pivotal. Note that \cite[Theorem 5.2]{Kril} gives that the twist of any \textit{\`etale} algebra coming from a conformal embedding has trivial twist. In particular diagrammatic techniques will work in more general settings (see e.g. \cite{Caleb}).
\end{remark}

The functor $\Phi$ is easy to define on objects. Recall that objects of $\mathcal{SE}_N$ are strings $s$ on the alphabet $\{+, -\}$. The functor $\Phi$ is defined on objects by
\[  \Phi(s) := X^{s_1} \otimes X^{s_2}  \otimes \cdots \otimes  X^{s_n}     \]
where we recall $X = \mathcal{F}_A(V_\square)$. By Lemma~\ref{lem:dom} the functor $\mathcal{F}_A:\overline{\operatorname{Rep}(U_q(\mathfrak{sl}_N))}\to \overline{\operatorname{Rep}(U_q(\mathfrak{sl}_N))}_A$ is Cauchy-dominant. As every object in $\overline{\operatorname{Rep}(U_q(\mathfrak{sl}_N))}$ is a direct summand of a tensor power of $V_\square$, Cauchy-dominance of $\mathcal{F}_A$ guarantees that every object of $\overline{\operatorname{Rep}(U_q(\mathfrak{sl}_N))}_A$ is a direct summand of a tensor power of $X$. As $\Phi(+^n) = X^{\otimes n}$ it follows that $\Phi$ is Cauchy-dominant.

To define the functor on morphisms, we first observe that the free module functor $\mathcal{F}_A: \overline{\operatorname{Rep}(U_q(\mathfrak{sl}_N))} \to \overline{\operatorname{Rep}(U_q(\mathfrak{sl}_N))}_A$ gives us morphisms
\begin{align*}
    \att{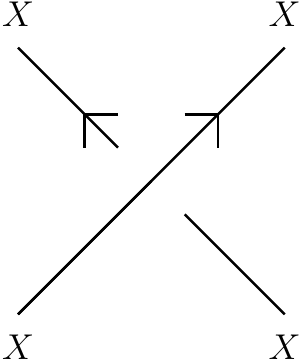}{.25} &:= \tau \circ \mathcal{F}_A\left( \att{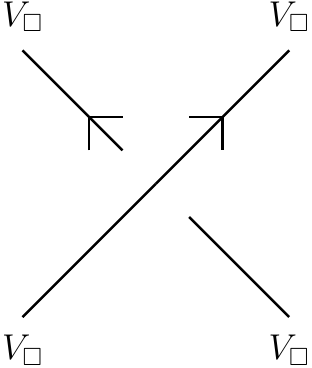}{.25} \right) \circ \tau^{-1} \in \End_{\overline{\operatorname{Rep}(U_q(\mathfrak{sl}_N))}_A}\left(X^{\otimes 2} \right), \\
    \att{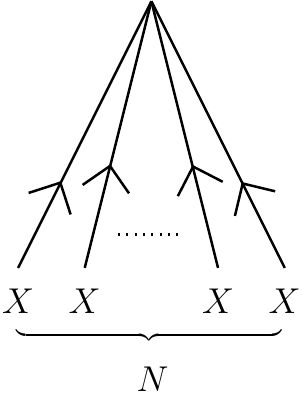}{.25} &:= {\tau'}\circ \mathcal{F}_A\left( \att{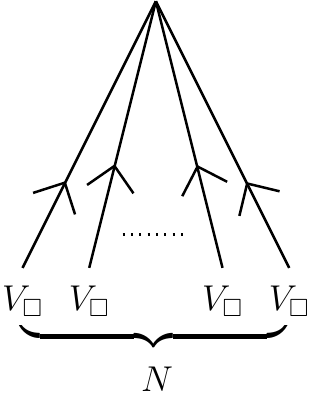}{.25} \right) \circ {\tau'}^{-1} \in \Hom_{\overline{\operatorname{Rep}(U_q(\mathfrak{sl}_N))}_A}\left(X^{\otimes N}\to \mathbf{1}\right)
\end{align*}   
where $\tau$ and $\tau'$ denote the appropriate composition of tensorators
As $\mathcal{F}_A$ is a tensor functor these morphisms satisfy the same relations their sources satisfied in $\overline{\operatorname{Rep}(U_q(\mathfrak{sl}_N))}$. By Proposition \ref{prop:sen}, this gives the following braid relations in $\overline{\operatorname{Rep}(U_q(\mathfrak{sl}_N))}_A$:
\[\att{R21}{.15} = \att{R22}{.15}, \quad \att{R31}{.15} = \att{R32}{.15},  \quad  \att{H1}{.15} = \att{H12}{.15} + (q-q^{-1})\att{H11}{.15}.\] 

Using $q = e^{2 \pi {\bf{i}}/4N}$, we also have
\[\att{loop}{.15} = \frac{q^{2N}-q^{-2N}}{q-q^{-1}} = \frac{2\mathbf{i}}{q-q^{-1}} \quad
    \text{and} \quad \att{R1}{.15} = q^N =\mathbf{i}\att{R12}{.15}.\]

We also get the relations
\[  \att{kwrel1}{.25} = q\att{kwrel12}{.25} \qquad \text{ and } \qquad \att{kwrel2}{.25} = \att{kwrel22}{.25}  \]
where $p_{1_N}$ is the usual formula in terms of braidings for the unique minimal projection in the HOMFLYPT skein category onto the simple $1_N$ in \cite[Theorem 6]{Paggo}.

The remaining generator of $\mathcal{SE}_N$ is the ``new stuff'', which doesn't come from the image of the free module functor. To obtain this ``new stuff'' we observe that from Lemma~\ref{lem:simp}, there exists an invertible object $g\in \overline{\operatorname{Rep}(U_q(\mathfrak{sl}_N))}_A$ such that \[\dim\Hom_{
\overline{\operatorname{Rep}(U_q(\mathfrak{sl}_N))}_A}(X\otimes X^* \to g) = 1.\] We have that $\overline{\operatorname{Rep}(U_q(\mathfrak{sl}_N))}_A$ is semisimple \cite[Proposition 2.7]{LagrangeUkraine}, as $A$ is \textit{\'etale}. Thus, there exists a projection  
\[ \att{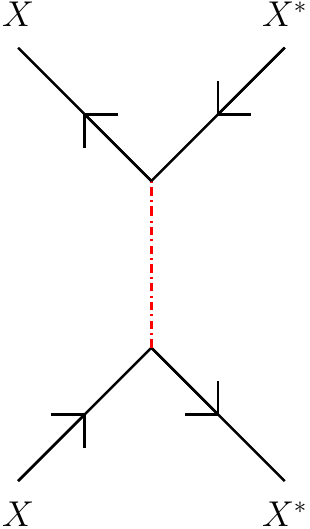}{.2} \in \End_{\overline{\operatorname{Rep}(U_q(\mathfrak{sl}_N))}_A}(X\otimes X^*)\]
onto the object $g$.

Since $g$ appears in $X \otimes X^*$ with multiplicity $1$, it is symmetrically self-dual (write the identity on $X \otimes X^*$ as a sum of projections and rotate by 180-degrees, using that a rotation of a complete set of projections is again a complete set of projections\footnote{See  \cite{MR3973457,MR1909537} for why this need not hold without the multiplicity $1$ assumption, and for some ways that assumption can be weakened}. That is,
\[\left(\att{splittingEq}{.2}\right)^* =   \att{splittingEq}{.2}. \]

\begin{remark}\label{rmk:Short}
    We introduce the following shorthand notation for our diagrams. As the new projection is invariant under 180 degrees rotation, this new notation is unambiguous.
    \[    \att{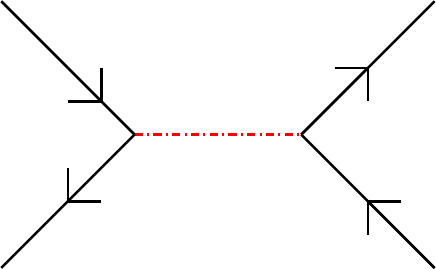}{.2} \quad := \quad  \att{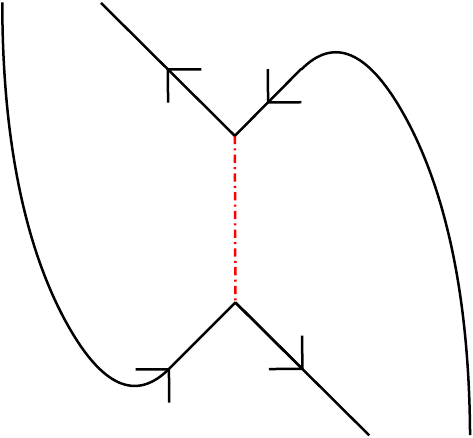}{.2}   \quad = \quad  \att{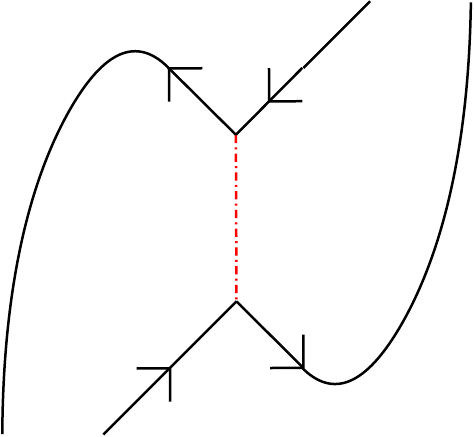}{.2}        \]
    and 
   \[    \att{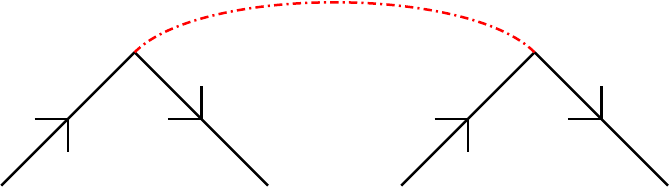}{.2} \quad := \quad  \att{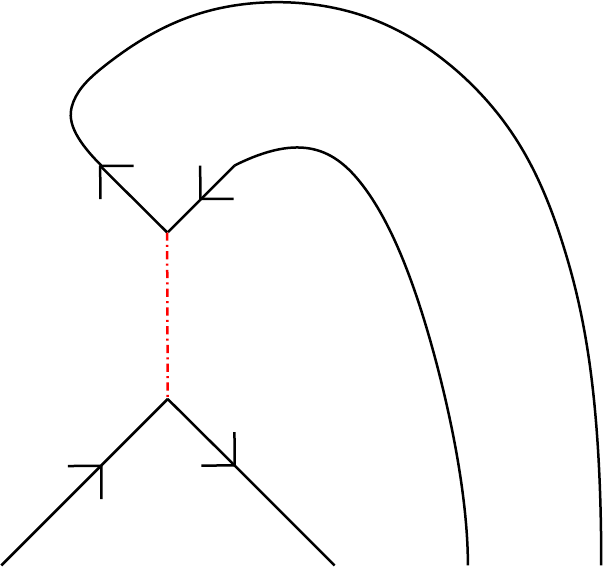}{.2}   \quad = \quad  \att{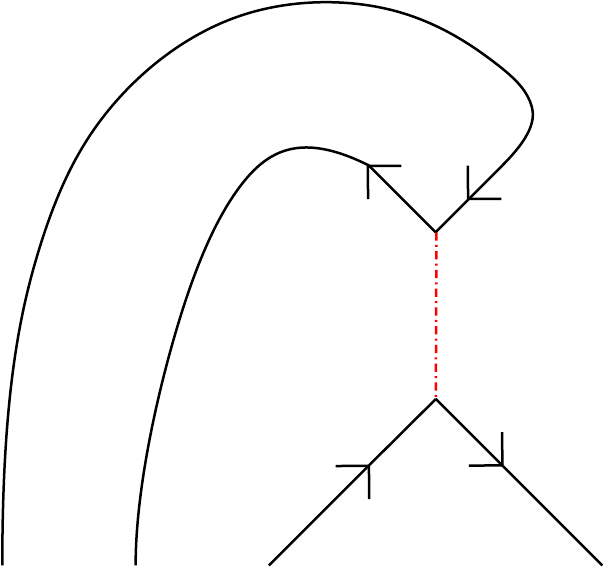}{.2}   .     \]
    
    We also use a similar shorthand for the 180 degree rotation of the above morphism.
\end{remark}

To find additional relations involving this new projection, we recall from Lemma \ref{lem:LocalDescription} that the object $g\in \overline{\operatorname{Rep}(U_q(\mathfrak{sl}_N))}_A$ generates a subcategory which is equivalent as a ribbon category to $\mathrm{SVec}$ with the odd line given dimension $1$.
This gives the relations
\[\att{lolly}{.2} = 0 \quad \text{and} \quad \att{redId}{.2} =  \att{redcap}{.2}\]

Since $X$ is simple (which is a consequence of Theorem~\ref{thm:dimEnd}), we have by Schur's lemma
    \[\att{Tr1}{.2} = c\att{Id1}{.2} \]
for some scalar $c$. The trace of the LHS is the trace of projection onto $g$ and hence $1 = \dim g = c \dim X$. Hence $c = \frac{q-q^{-1}}{2\mathbf{i}}$. 

We obtain the final relation 
\[ \quad \att{half1}{.25} = \att{half2}{.25}\]
from Corollary~\ref{cor:OB} and (R2).

We have now shown that all the defining relations of $\mathcal{SE}_N$ hold in $\overline{\operatorname{Rep}(U_q(\mathfrak{sl}_N))}_A$. Hence we have that $\Phi$ is a tensor functor. We summarise the results of this subsection in the following lemma.
\begin{lem}\label{lem:sur}
    Let $N\in \mathbb{N}_{\geq 2}$ and set $q = e^{2\pi i \frac{1}{4N}}$. Then there exists a Cauchy-dominant pivotal tensor functor
    \[ \Phi:   \mathcal{SE}_N \to \overline{\operatorname{Rep}(U_q(\mathfrak{sl}_N))}_A  \]
    defined on objects by 
    \[  \Phi(s) := X^{s_1} \otimes X^{s_2}  \otimes \cdots \otimes  X^{s_n}         \]
    and on generating morphisms by 
    \[    \Phi\left(\att{braidX}{.25}\right) := \att{braidXX}{.25} \qquad     \Phi\left(\att{kw}{.25}\right) := \att{kwXX}{.25}   \qquad \Phi\left(\att{splittingEq}{.2}\right) := \att{splittingEqXX}{.2} .     \]
\end{lem}

\subsection{Faithfulness}
The functor $\Phi$ from Lemma~\ref{lem:sur} is certainly not faithful. For example, the projection onto $(N+1)$ in the subcategory $\mathcal{H}(q, \mathbf{i})$ is sent to zero. However, we will show that the functor $\Phi$ descends to a faithful functor $\overline{ \mathcal{SE}_{N}} \to \overline{\operatorname{Rep}(U_q(\mathfrak{sl}_N))}_A$ where $\overline{ \mathcal{SE}_{N}}$ is the quotient by the negligible ideal. This will follow by a standard argument, first we show that the  unit of $\mathcal{SE}_N$ is simple, and then it follows that the kernel of $\Phi$ must be equal to the negligible ideal of $\mathcal{SE}_{N}$. 
For use later on, we will prove the stronger result that $\mathcal{E}_R$ has simple unit.

By definition, the hom spaces of $\mathcal{E}_R$ are spanned by morphisms constructed from the generators and duality morphisms, using the operations $\circ$ and $\otimes$. It will be convenient to give these morphisms a name.

\begin{defn}
We will refer to a morphism in $\mathcal{E}_R$ constructed from the generators $\att{braidX}{.15}$, $\att{splittingEq}{.15}$, and duality morphisms, using $\circ$ and $\otimes$ as a \textit{diagram}. We will refer to the hom space of $\mathcal{E}_R$ that $D$ is an element of as the \textit{boundary} of a diagram.
\end{defn}

For example, we have that
\begin{equation}\label{eq:dia}  \att{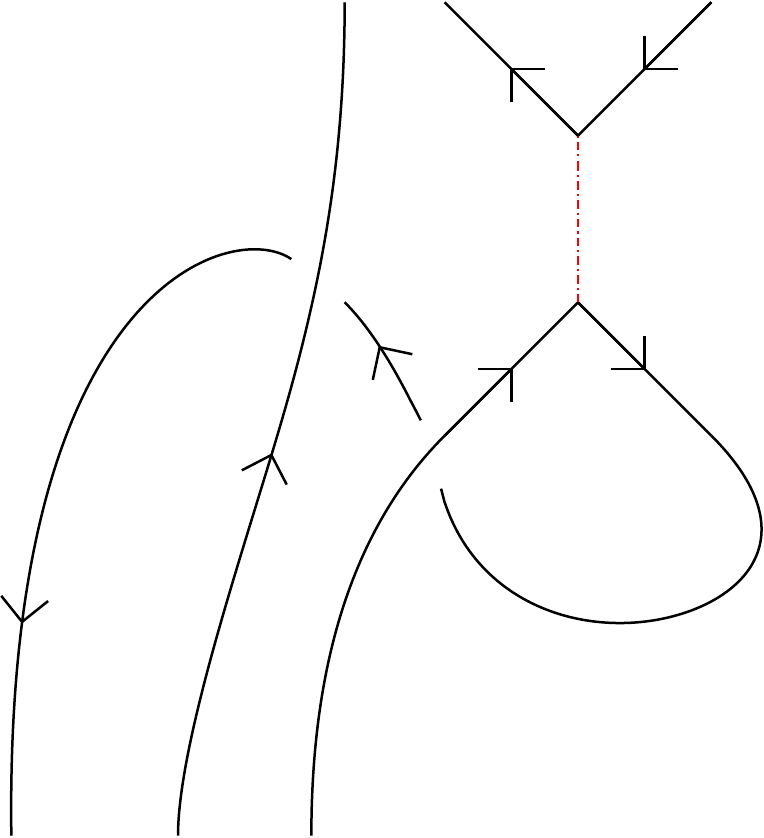}{.2}    \end{equation}
is a diagram with boundary $(-++\to ++-)$, and that
\[  q\att{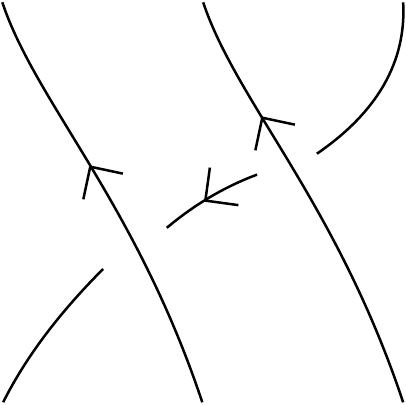}{.2} + (q-q^{-1}) \att{exampleD}{.2}  \]
is an $R$-linear combination of diagrams with boundary $(-++\to ++-)$. It follows from the definition of $\mathcal{E}_R$ as a category given by generators and relations that the hom spaces $\mathcal{E}_R(s_1\to s_2)$ are spanned over $R$ by diagrams with boundary $(s_1\to s_2)$. For a given diagram $D$, we will refer to a sub-diagram of the form $\att{braidX}{.15}$ as a braid term, and a sub-diagram of the form $\att{splittingEq}{.15}$ as a projection term. 

\begin{defn}
     We say two black points in a diagram $D$ are in the same \textit{strand} if they are connected in the diagram obtained from $D$ via removing all projection terms. By construction the two endpoints of a strand in $D$ are either at the boundary, or at a projection term. We say the strand \textit{connects} these two regions.
\end{defn}

The diagram in Equation \ref{eq:dia} has 5 strands.

\begin{defn}
    A self-intersection is a crossing where a strand crosses itself.
\end{defn}

The diagram in Equation \ref{eq:dia} has no self-intersections.

\begin{defn}
    We call a black strand in a diagram a topmost strand (resp. bottommost strand), if it never crosses under (resp. over) strands other than itself.
\end{defn}

In Equation \ref{eq:dia}, the two strands connecting to the top of the generator are both topmost and bottommost. The strand connecting the lower-right of the generator to the lower-left of the diagram is bottommost but not topmost. The remaining two strands are topmost, but not bottommost.

We will often induct on our diagrams with respect to the following partial ordering. 

\begin{defn}\label{def:houseordos}
    We define a partial ordering on diagrams where $D_1 < D_2$ if one of the following holds:
    \begin{itemize}
    \item $D_1$ has fewer projection terms than $D_2$,
    \item they have the same number of projection terms but $D_1$ has fewer self-intersections,
    \item or they have the same number of projection terms and self-intersections, but $D_1$ has fewer braid terms. 
    \end{itemize}
    We say $D_1 \sim D_2$ if they have the same number of projection terms, self-intersections, and braid terms.
\end{defn}

\begin{defn}
    We say that a diagram $D$ is in \textit{standard form} if it has no self-intersections, and every strand connected to a projection term also connects to the boundary.
\end{defn}

Taking the contrapositive, the following lemma follows easily.

\begin{lem}
$D$ is in standard form iff  it has no self-intersections and each connected component contains at most one projection term, and no projection term connects to itself.
\end{lem}

\begin{ex}
   We have that the diagram from Equation~\eqref{eq:dia} contains two connected components. This diagram is in standard form as one of these has one projection term, and the other has no projection terms. On the other hand, the diagram
   \[  \att{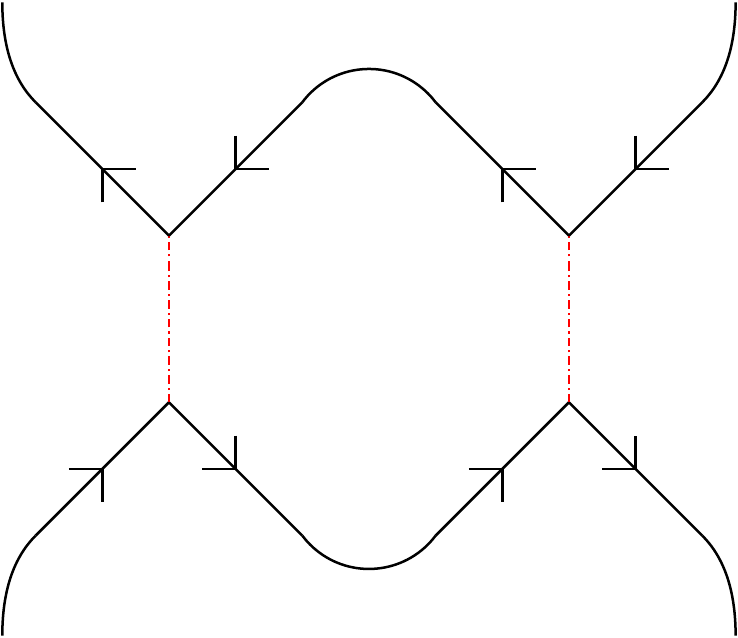}{.2}   \]
   is not in standard form as it has two projection terms connected to each other.
\end{ex}

The following Proposition shows that any diagram in $\mathcal{E}_R$ can be expressed as a linear combination of diagrams in standard form. This technical result will be the key to showing that $\mathcal{E}_q$ has simple unit. In the next section this proposition will also be useful for producing a basis for the hom spaces of $\mathcal{SE}_{N}$.

\begin{prop}\label{prop:standardForm}
    Each Hom space of $\mathcal{E}_R$ is spanned over $R$ by the diagrams in standard form.
\end{prop}

We will need two lemmas for the proof.

\begin{lem} \label{lem:topmost}
    For any diagram $D$ and any choice of strand in $D$, we have that $D = \sum_{i=1}^k a_i D_i$ where $D \sim D_1$, $D_i < D$ for all $1<i\leq k$, and the chosen strand is topmost (resp. bottommost) \edit{in $D_1$}.
\end{lem}
\begin{proof}
    If the chosen strand is not topmost, then we can apply (Hecke) to each of the $k$ times it crosses under another strand, express $D = \sum_{i=1}^{2^k} a_i D_i$ with $D_i < D$ for $i>1$ and where $D_1 \sim D$ but where the strand connecting the two generators is topmost \edit{in $D_1$}. The proof for bottommost is identical.
\end{proof}

\begin{lem}\label{lem:selfintersection}
    If a diagram $D$ has at least one self-intersection, then we have that $D = \sum_{i=1}^k a_i D_i$ where $D_i < D$ for all $1\leq i\leq k$.
\end{lem}
\begin{proof}
    We induct on the number of \edit{self-intersections}. Fix a strand with a self-crossing. By Lemma \ref{lem:topmost}, WLOG we can assume that the strand is topmost.

     We now follow the usual technique for evaluating the HOMFLYPT polynomial. Fix a strand that has self-intersections and fix an end of that strand. Starting at the chosen end of strand move along the strand, and each time you reach a crossing for the first time, if it goes under you do nothing, and if it goes over apply (Hecke) to rewrite it as a linear combination of two other diagrams. After doing this switching $n$ \edit{self-intersections} on that strand, we get $D = \sum_{i=1}^{2^n} a_i D_i$ where $D_i$ has fewer self-intersections than $D$ for $i>1$, and where $D_1$ has the same number of self-intersections as $D$ but the strand is always increasing out of the page. Finally we can apply (R1), (R2), (R3), and (Half-Braid) to $D_1$ to isotope the chosen strand until it has no \edit{self-intersections}.
\end{proof}

\begin{proof}[Proof of Proposition~\ref{prop:standardForm}]
The proof is by induction on the partial ordering $<$. The base cases are the diagrams with no projections or crossings, which are obviously in standard form.

For the induction step we need only show that every diagram $D$ which is not in standard form can be written as a linear combination of diagrams $\sum a_i D_i$ with $D_i< D$. By Lemma \ref{lem:selfintersection}, WLOG we can assume that $D$ has no self-intersections.

We split into two cases: either $D$ contains a projection term which is connected to itself, or two distinct projection terms connected by a black strand.

\begin{bf}Case 1: $D$ has a chosen black strand which connects a chosen projection term to itself.\end{bf}

By Lemma \ref{lem:topmost}, WLOG we can assume the strand is topmost. If there is anything inside the region bounded by the chosen generator and strand, we can use (R2), (R3), and (Half-Braid), to move the strand across each of the generators or strands appearing in this region, to see $D = D'$ where $D' \sim D$ and where the region bounded by the chosen generator and strand is empty.

Altogether, we may assume WLOG that one of the following two configurations appears within our diagram:
\[\att{Tr1}{.2}\quad \text{ or } \quad \att{lolly}{.2}.\]
The relations (Trace) and (Tadpole) show that this diagram is either $0$, or a scalar multiple of a diagram with one less projection term. Thus, we have expressed $D$ as a $R$-linear combination of strictly smaller diagrams.

\begin{bf}Case 2: $D$ has a chosen black strand connecting a first chosen projection to a second distinct chosen projection.\end{bf}

By Lemma \ref{lem:topmost}, we may assume WLOG that the chosen strand is \emph{bottommost}. We repeatedly use \edit{the following local  move} to slide one of the generators across the diagram until they are adjacent. \edit{This local move is a consequence of (Half-Braid) and (R2).}
\[ \edit{\att{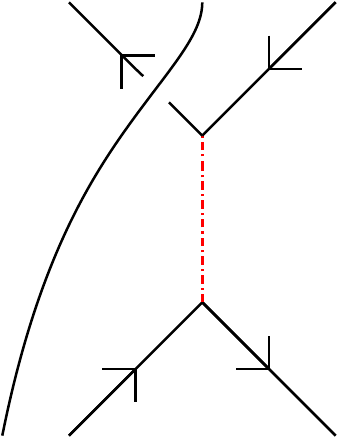}{.2}\quad = \quad \att{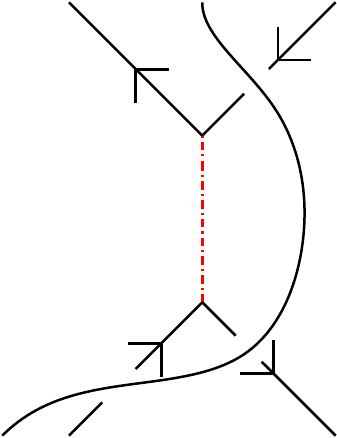}{.2} }\]
(Note that this step will increase the number of crossings, since each time we apply (Half-Braid) we replace a single crossing with three crossings! However, since there were no \edit{self-intersections} \edit{the new diagram} will still have no \edit{self-intersections}.)  We then rotate the generators and apply the (Dual) relation so that the two adjacent generators are in the following configuration

  \[\att{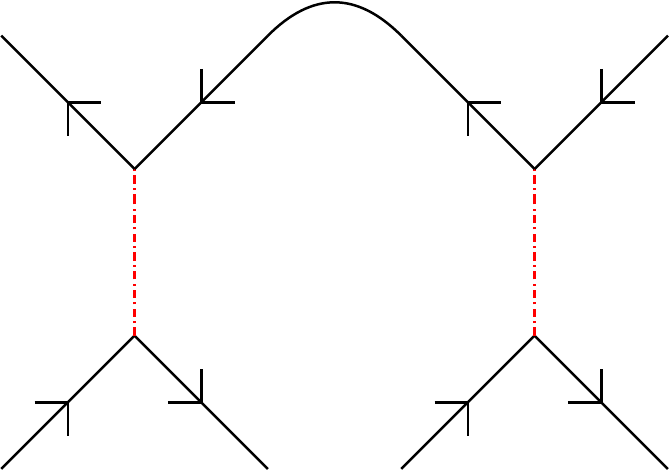}{.2}.\]
This now simplifies
\[\att{splittingEqEv1}{.2} =  \att{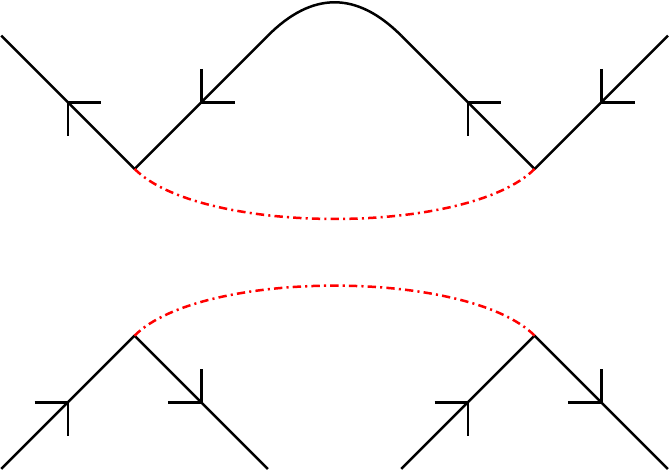}{.2} = \frac{q-q^{-1}}{2\mathbf{i}} \att{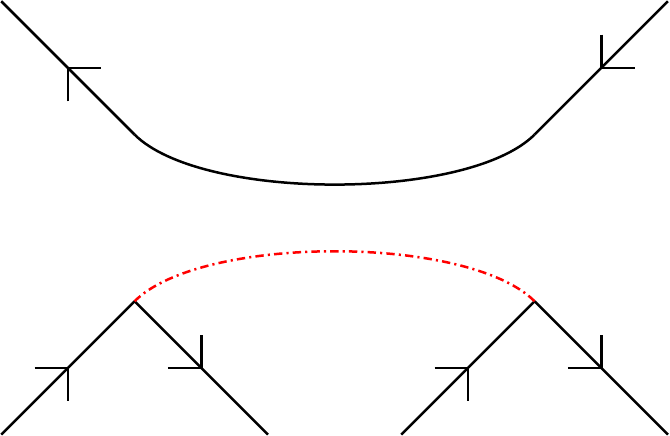}{.2},  \]
using relation ($\mathbb{Z}_2$) in the first step, and relation (Trace) \edit{}in the second step. This expresses the diagram as a scalar multiple of a strictly smaller diagrams.\end{proof}

An immediate corollary shows that if the specialization $\mathcal{E}_q$ is non-trivial, then it has simple unit.

\begin{cor}\label{cor:unit}
We have that $\dim\operatorname{End}_{\mathcal{E}_q}(\mathbf{1} )\leq 1$.
\end{cor}
\begin{proof}
    By Proposition~\ref{prop:standardForm} we have that $\operatorname{End}_{\mathcal{E}_q}(\mathbf{1} )$ is spanned by diagrams in standard form. Further, these diagrams must have no boundary. Clearly, any diagram without a boundary and in standard form can not have any projection terms. Thus $\operatorname{End}_{\mathcal{E}_q}(\mathbf{1} )$ is spanned by closed braid diagrams. It is a well-known result that any closed braid diagram can be evaluated to a scalar multiple of the empty diagram using the HOMFLYPT skein category relations.
\end{proof}

\begin{warning}
    The same argument shows that $\operatorname{End}_{\mathcal{E}_R}(\mathbf{1})$ is a quotient of $R$, but since $R$ is not a field this could a priori be a non-zero proper quotient. We will show later in Theorem~\ref{thm:existence} that in fact $\operatorname{End}_{\mathcal{E}_R}(\mathbf{1}) = R$.
\end{warning}

Recall our aim was to show that the categories $\mathcal{SE}_N$ have simple unit. This will follow now by showing that when $q = e^{2\pi i \frac{1}{4N}}$, the endomorphism algebras of $\mathcal{SE}_N$ are equal to the endomorphism algebras of the subcategory $\mathcal{E}_q$. We prove the following slightly stronger result.

\begin{lem}\label{lem:end}
Let $N\in \mathbb{N}_{\geq 2}$ and $q = e^{2\pi i \frac{1}{4N}}$, and $s_1,s_2$ strings in $\{+,-\}$ such that $\sum s_1 = \sum s_2$. Then
\[      \operatorname{Hom}_{ \mathcal{SE}_N}(s_1 \to s_2) =  \operatorname{End}_{ \mathcal{E}_q}(s_1\to s_2) .   \]
\end{lem}

We will need one definition.

\begin{defn}
Let $f\in\operatorname{Hom}_{ \mathcal{SE}_N}(t_1 \to t_2)$. We define the degree of f as $\sum t_1 - \sum t_2$. 
\end{defn}
The degree of cups, caps, crossings, and projections are all $0$, while $\att{kw}{.2}$ has degree $n$ and \att{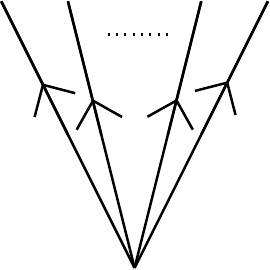}{.2} has degree $-n$.

\begin{proof}[Proof of Lemma~\ref{lem:end}]
    By construction we have that $\mathcal{SE}_N$ is an extension of $\mathcal{E}_q$ by the additional generator \att{kw}{.2}. As usual, we will use the word diagram for homs which are compositions of the basic generators, and diagrams span all morphisms.
    Hence it suffices to show that any diagram in $\operatorname{Hom}_{ \mathcal{SE}_N}(s_1 \to s_2)$ can be expressed in terms of diagrams with no \att{kw}{.2} terms. 

    By assumption we are considering morphisms of degree zero, so the number of $\att{kw}{.2}$ appearing in such a diagram must equal the number of \att{kwdag}{.2}. Each pair of such generators can be brought close to each other (at the cost of scalars) by ($q$-braid), and annihilated by (Pair). This leaves a diagram with no $\att{kw}{.2}$ terms.
\end{proof}

We thus have the desired corollary.
\begin{cor}\label{cor:simpunit}
    Let $N \in \mathbb{N}_{\geq 2}$. Then
    \[      \dim\operatorname{End}_{ \mathcal{SE}_N}(\mathbf{1})   =1.  \]
\end{cor}

This allows us to apply the general theory to produce a faithful functor into $\overline{\operatorname{Rep}(U_q(\mathfrak{sl}_N))}_A$.
\begin{cor}\label{cor:exist}
    For each $N\in \mathbb{N}_{\geq 2}$ there is a faithful functor 
    \[ \overline{\Phi}: \overline{\mathcal{SE}_{ N }} \to \overline{\operatorname{Rep}(U_q(\mathfrak{sl}_N))}_A  \]
     such that the following diagram commutes
     \[\begin{tikzcd}
    \mathcal{SE}_N \arrow{r}{\Phi} \arrow{d}{} & \overline{\operatorname{Rep}(U_q(\mathfrak{sl}_N))}_A    \\
 \overline{\mathcal{SE}_N } \arrow{ur}{\overline{\Phi}}&
\end{tikzcd}\]   
Moreover, $\overline{\mathcal{SE}_{ N }}$ is \edit{a} unitary category.
\end{cor}
\begin{proof}
    We have that $\overline{\operatorname{Rep}(U_q(\mathfrak{sl}_N))}_A$ is a unitary fusion category by Proposition~\ref{rmk:unitary}. The result then follows from Proposition~\ref{prop:descent}. 
\end{proof}

\begin{remark}\label{rmk:N2}
    When $N=2$ we have that 
    \[\att{splittingEq}{.2} - \att{idf2}{.2} + \frac{1}{\sqrt{2}} \att{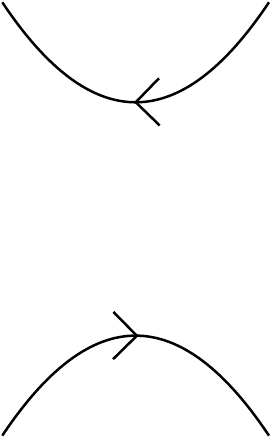}{.2}\]
    is negligible. 
    In this case $\overline{\mathcal{SE}_2}$ is monoidally equivalent to the Temperley-Lieb-Jones category at $\delta = \sqrt{2}$, and the projection onto the red strand is isomorphic to the second Jones-Wenzl projection. When $N >2$ a direct calculation shows that no linear combination \[x \att{splittingEq}{.2} +y \att{idf2}{.2} + z \att{capf2}{.2}\] is negligble, and thus $\overline{\mathcal{SE}_N}$ is strictly larger than $\overline{\operatorname{Rep}(U_q(\mathfrak{sl}_N))}$ for $N \geq 3$.
\end{remark}

\subsection{Fullness}
To complete Theorem~\ref{thm:special}, we have to show that $\overline{\Phi}$ is full. As we know this functor is Cauchy-dominant and faithful, this is equivalent to showing that $ \overline{\Phi}: \overline{\mathcal{E}_{ N }} \to \overline{\operatorname{Rep}(U_q(\mathfrak{sl}_N))}_A$ induces a monoidal equivalence
\[   \operatorname{Ab}(\overline{\mathcal{SE}_N})\simeq \overline{\operatorname{Rep}(U_q(\mathfrak{sl}_N))}_A .  \]
The key idea behind this proof is the Galois-like correspondence between intermediate categories
\[  \mathcal{C} \twoheadrightarrow \mathcal{D} \twoheadrightarrow \mathcal{C}_A   \]
and sub-algebras of $A$, as in Lemma \ref{lem:hack} and Conjecture \ref{con:Nonhack}. For our algebra which corresponds to the conformal embedding $\mathcal{V}(\mathfrak{sl}_N, N)\subset \mathcal{V}(\mathfrak{so}_{N^2-1}, 1)$, the sub-algebras are classified by the results of \cite{papi} and \cite{Kril}. This case is relatively tractable since the free fermion algebra $\operatorname{Fer}(\mathfrak{g})$ is generated as an algebra by $\mathfrak{g}$.

\begin{thm}\label{thm:equiv}

For all $N\in \mathbb{N}_{\geq 2}$ we have that the functor $\overline{\Phi}: \overline{\mathcal{SE}_N}\to \overline{\operatorname{Rep}(U_q(\mathfrak{sl}_N))}_A$ from Corollary~\ref{cor:exist} induces a monoidal equivalence
\[  \operatorname{Ab}(\overline{\mathcal{SE}_N})\simeq \overline{\operatorname{Rep}(U_q(\mathfrak{sl}_N))}_A.    \]
\end{thm}

We will first require some lemmas. We let $\mathcal{I}: \mathcal{SH}(q, q^N) \rightarrow \mathcal{SE}_N$ denote the obvious inclusion map.

\begin{lem} \label{lem:Descent}
    The functors $\mathcal{I}$ and $\Phi$ descend to the semisimplifications. The functors they induce
    \[\begin{tikzcd}
   \overline{\operatorname{Rep}(U_q(\mathfrak{sl}_N))} \cong \operatorname{Ab}(\overline{\mathcal{SH}(q, q^N)})  \arrow{r}{\mathcal{F}_1} & \operatorname{Ab}(\overline{\mathcal{SE}_N}) \arrow{r}{\mathcal{F}_2} &    \overline{\operatorname{Rep}(U_q(\mathfrak{sl}_N))}_A
\end{tikzcd}\] are faithful and Cauchy-dominant.  
\end{lem}
\begin{proof}
    We already have that $\Phi$ descends to the semisimplification and $\overline{\mathcal{SE}_N}$ is unitary. It then follows from Proposition \ref{prop:descent} that $\mathcal{I}$ also descends. Since any tensor ideal is negligible, $\overline{\mathcal{I}}$ and $\overline{\Phi}$ are faithful, hence by \ref{lem:AbFaithful} the functors $\mathcal{F}_1$ and $\mathcal{F}_2$ are faithful. 
    By Lemma \ref{lem:dom}, we also have that $\mathcal{F}_1$ and $\mathcal{F}_2$ are Cauchy dominant.
\end{proof}

In order to apply our general set-up of Subsection~\ref{sec:etale}, we need the following.
\begin{lem} \label{lem:UnitaryTricky}
    $\operatorname{Ab}(\overline{\mathcal{SE}_N})$ is a unitary fusion category.
\end{lem}
\begin{proof}
    We already know $\operatorname{Ab}(\overline{\mathcal{SE}_N})$ is unitary by Corollary \ref{cor:exist} and it is Cauchy-complete by construction, so we need only show that it has finitely many isomorphism classes of simple objects. But $\overline{\operatorname{Rep}(U_q(\mathfrak{sl}_N))}$ is fusion and $\mathcal{F}_1$ is Cauchy-dominant, and so every simple object in $\operatorname{Ab}(\overline{\mathcal{SE}_N})$ is isomorphic to one of finitely many summands of the image of one of the isomorphism classes of simple objects in $\overline{\operatorname{Rep}(U_q(\mathfrak{sl}_N))}$.
\end{proof}

The general set-up now gives the following.
\begin{lem}\label{lem:equiv1}
    There exists an \textit{\'etale} subalgebra object $A' \subseteq A$, such that \[  \operatorname{Ab}(\overline{\mathcal{SE}_N})\simeq \overline{\operatorname{Rep}(U_q(\mathfrak{sl}_N))}_{A'}.    \]
\end{lem}
\begin{proof}
Since it maps the half-braiding to the half-braiding, we have that $\mathcal{F}_2$ is a $\overline{\mathcal{SH}(q, q^N)}$-linear functor. By Lemma \ref{lem:Descent} and Lemma \ref{lem:UnitaryTricky}, it follows from Lemma \ref{lem:OverToEtale2} that $A' = \mathcal{F}_1^\vee(\mathbf{1})$ is an \textit{\'etale} algebra object in $\overline{\operatorname{Rep}(U_q(\mathfrak{sl}_N))}$ and by Lemma \ref{lem:hack} it is a subalgebra of $\mathcal{F}_{\edit{2}}^\vee(\mathbf{1})\cong  A$.
\end{proof}

We now arrive at the key lemma.

\begin{lem} \label{lem:equiv2}
    Suppose that $N\geq 3$. If $A' \subseteq A$ is a subalgebra object, such that \[  \operatorname{Ab}(\overline{\mathcal{SE}_N})\simeq \overline{\operatorname{Rep}(U_q(\mathfrak{sl}_N))}_{A'}, \] 
    then  \[\Hom_{\overline{\operatorname{Rep}(U_q(\mathfrak{sl}_N))}}\left(V_{\ydiagram{2}} \to  A'\otimes V_{\ydiagram{1,1}}\right )\] is non-zero.
\end{lem}
\begin{proof}
In $\mathcal{SH}(q, q^N)$ we have the following projections which map to the projections onto $V_{\ydiagram{1,1}}$ and $V_{\ydiagram{2}}$ in $\overline{\operatorname{Rep}(U_q(\mathfrak{sl}_N))}$

\[  p_{V_{\ydiagram{1,1}}}:= \frac{1}{1+q^{-2}}\left( \att{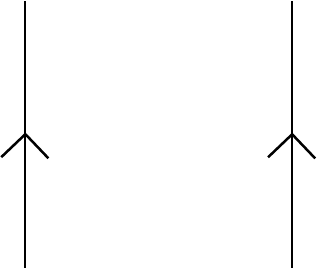}{.2}-q^{-1} \att{braidX}{.2} \right) \qquad \text{ and } \qquad p_{V_{\ydiagram{2}}}:= \frac{1}{1+q^{2}}\left( \att{idid}{.2}+q \att{braidX}{.2} \right)  \]

Furthermore, a direct computation in $\mathcal{SE}_N$ gives that
\[  (1-q^2)\att{idid}{.2} + (q-q^{-1}) \att{braidX}{.2} - \mathbf{i}(1+q^2) \att{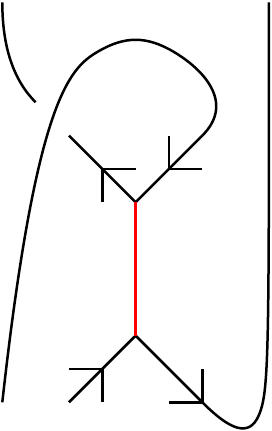}{.2} + \mathbf{i}(q+q^{-1}) \att{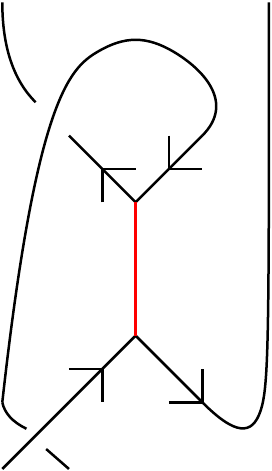}{.2}   \]
is an intertwiner from $ p_{V_{\ydiagram{2}}} \to  p_{V_{\ydiagram{1,1}}}$. Furthermore, when $N\geq 3$, this intertwiner has non-zero trace, and is thus non-zero in $\overline{ \mathcal{SE}_N}$.

We thus get
\begin{align*}
    1 &\leq \dim\Hom_{\operatorname{Ab}(\overline{ \mathcal{SE}_N})} 
 \left(\mathcal{F}_1(V_{\ydiagram{2}} ) \to  \mathcal{F}_1\left(V_{\ydiagram{1,1}} \right) \right) \\
 &= \dim\Hom_{\overline{\operatorname{Rep}(U_q(\mathfrak{sl}_N))}}\left(V_{\ydiagram{2}} \to  \mathcal{F}_1^\vee(\mathcal{F}_1( V_{\ydiagram{1,1}}))\right )\\
 &= \dim\Hom_{\overline{\operatorname{Rep}(U_q(\mathfrak{sl}_N))}}\left(V_{\ydiagram{2}} \to  \mathcal{F}_{A'}^\vee(\mathcal{F}_{A'}( V_{\ydiagram{1,1}}))\right )\\
 &= \dim\Hom_{\overline{\operatorname{Rep}(U_q(\mathfrak{sl}_N))}}\left(V_{\ydiagram{2}} \to  A'\otimes V_{\ydiagram{1,1}}\right ) 
\end{align*}  
\end{proof}

By Lemma~\ref{lem:equiv1} we need to study sub-algebra objects of $A$. To discuss these we introduce the following notation.

\begin{defn}
    An algebra object in a fusion category is called pointed if as an object it is a direct sum of invertible objects.
\end{defn}
The following result gives a nice characterisation of sub-algebras of pointed algebra objects.
\begin{lem} \label{lem:PointedSubalgebra}
    If $A$ is a pointed \textit{\'etale} algebra in a fusion category, and moreover $A = \bigoplus_{i\in G} X_i$ is a multiplicity-free decomposition of $A$, then $G$ is a group under the operation characterized by $X_i \otimes X_j \cong X_{i\cdot j}$. If $H$ is a subgroup of $G$ then $\bigoplus_{i\in H} X_i$ is a subalgebra of $A$, and moreover all subalgebras are of this form.
\end{lem}
\begin{proof}
    Multiplicity-free guarantees that $G$ forms a monoid under tensoring isomorphism classes. Since \textit{\'etale} algebras are self-dual, it forms a group with inverse characterized by $X_{i^{-1}} \cong X_i^*$. 

    If $H$ is a subgroup, then the restriction of the product in $A$ to $X_i \otimes X_j$ must land in $X_{i \cdot j}$, and thus $\bigoplus_{i\in H} X_i$ is a \edit{subalgebra}. If $B \subset A$ is a subalgebra, then by the first statement, the summands which appear again form a group.
\end{proof}

We are now in place to prove our main theorem regarding fullness.
\begin{proof}[Proof of Theorem \ref{thm:equiv}]
In the case of $N=2$ we have that $A = \mathbf{1}$ and the statement of the theorem is trivially true by Remark~\ref{rmk:N2}.

We now assume $N \geq 3$. By Lemma \ref{lem:equiv1} we have that there exists a subalgebra object $A' \subseteq A$, such that \[  \operatorname{Ab}(\overline{\mathcal{SE}_N})\simeq \overline{\operatorname{Rep}(U_q(\mathfrak{sl}_N))}_{A'}.    \]
The sub-algebras of $A$ are classified by \cite[Theorem 2.1]{papi} together with \cite[Theorem 5.2]{Kril}. These results show that either $A' = A$, or that $A'$ is a sub-algebra of the maximal pointed sub-algebra of $A$. Our goal is to show $A=A'$.

If $A'$ is pointed, then Lemma \ref{lem:equiv2} implies that there is a non-trivial invertible object $h$ in $\overline{\operatorname{Rep}(U_q(\mathfrak{sl}_N))}$ such that $h\otimes V_{\ydiagram{1,1}}\cong V_{\ydiagram{2}}$. An explicit description of the invertible objects in $\overline{\operatorname{Rep}(U_q(\mathfrak{sl}_N))}$ and their tensor product formulae can be found in \cite[Section 2.2]{levelterry}. With this information we find that such an $h$ only exists when $N=3$, in which case $h \cong V_{\ydiagram{3}}$. Hence when $N>3$ we have that $A' = A$.

When $N=3$, the maximal pointed algebra is $V_{\emptyset}\oplus V_{\ydiagram{3}} \oplus V_{\ydiagram{3,3}}$, which in this case is also $A$. Since $3$ is prime, the group of simple summands must be $\mathbb{Z}/3$, and thus $A'$ has no non-trivial subalgebras by Lemma \ref{lem:PointedSubalgebra}. We have that $A' \neq 1$, because $\overline{\mathcal{SE}_N}$ is strictly larger than $\overline{\operatorname{Rep}(U_q(\mathfrak{sl}_N))}$ (see Remark \ref{rem:Nis2}). Thus $A'=A$, and by Lemma \ref{lem:hack}, $\mathcal{F}_2$ is an equivalence. In particular, we have
 \[ \operatorname{Ab}(\overline{\mathcal{SE}_N}) \simeq \overline{\operatorname{Rep}(U_q(\mathfrak{sl}_N))}_{A'}  \simeq \overline{\operatorname{Rep}(U_q(\mathfrak{sl}_N))}_A    \]
as desired.

\end{proof}

As the endomorphism algebras of $\mathcal{SE}_N$ and $\mathcal{E}_q$ are the same at $q=e^{2\pi i \frac{1}{4N}}$ we obtain the following corollary. This result will be useful for constructing bases for the Hom spaces of $\mathcal{E}_q$ in the next section.
\begin{cor}\label{cor:full}
   Let $N\in \mathbb{N}_{\geq 2}$ and set $q = e^{2 \pi i \frac{1}{4N}}$. Then we have
   \[ \dim \End_{\overline{\mathcal{E}_q}}(+^{n}) =\dim \End_{\overline{\operatorname{Rep}(U_q(\mathfrak{sl}_N))}_A}(X^{\otimes n}).     \]
\end{cor}
\begin{proof}
It follows from Theorem~\ref{thm:equiv} that the functor $\overline{\Phi}: \overline{\mathcal{SE}_N} \to \overline{\operatorname{Rep}(U_q(\mathfrak{sl}_N))}_A$ is fully faithful. Thus
\[       \dim \End_{\overline{\mathcal{SE}_N}}(+^{n}) =\dim \End_{\overline{\operatorname{Rep}(U_q(\mathfrak{sl}_N))}_A}(X^{\otimes n}).     \]
We have from Lemma~\ref{lem:end} that $\End_{\mathcal{SE}_N}(+^{n}) = \End_{\mathcal{E}_q}(+^{n})$, and by definition they have the same trace, thus the negligible ideal of these two endomorphism algebras are equal. Hence we have
\[\End_{\overline{\mathcal{SE}_N}}(+^{n}) = \End_{\overline{\mathcal{E}_q}}(+^{n})\]
which completes the proof.
\end{proof}

\section{Existence of the interpolating category \texorpdfstring{$\mathcal{E}_R$}{E}}\label{sec:eq}

In the previous section we showed that the categories $\mathcal{SE}_N$ were non-trivial for all $N\in \mathbb{N}_{\geq 2}$. In particular, this implies that the subcategories $\mathcal{E}_q$ are non-trivial for $q = e^{2 \pi i \frac{1}{4N}}$. In this section we show that $\mathcal{E}_q$ are non-trivial for all $q\in \mathbb{C}-\{-1,0,1\}$. In fact, we show a stronger result that the Hom spaces for $\mathcal{E}_R$ are free $R$-modules of specific ranks. 

That our putative $R$-basis spans is a diagrammatic argument which is the hardest part of the construction. Independence is easier and follows from general \textit{Deligne interpolation} techniques similar to that used in \cite{Deligne}. Namely, given an $R$-linear dependence, we can specialize each coefficient to $q= e^{2 \pi i \frac{1}{4N}}$ and use our results from the last section to see that these rational \edit{functions vanish} for infinitely many values and hence must be identically $0$. As a consequence, we obtain that the endomorphism algebras of $+^n$ in $\mathcal{E}_q$ are Hecke-Clifford algebras. This will be important in \cite{HansNew}, which studies the representation theory of these endomorphism algebras.

For the following definition we remind our reader of the shorthand notation of Remark~\ref{rmk:Short}.
\begin{defn}
For ease of notation, we define the following morphism in $\operatorname{End}_{\mathcal{E}_R}(++)$:
\[   \att{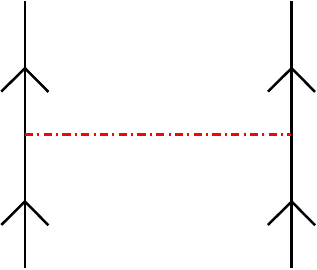}{.25} := - \frac{2 }{q-q^{-1}} \att{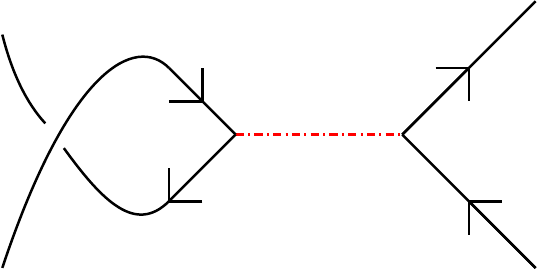}{.25} .  \]
We will refer to this morphism as a \textit{ladder}.
    \end{defn}

Note that either crossing can be used in the ladder morphism, due to the relations (Hecke) and (Tadpole). Direct computation gives the following relations between braids and ladders.

\begin{lem}\label{lem:relationspic}
The following relations hold in $\mathcal{E}_R$
\begin{align*}
    \text{(Over-Braid)} \att{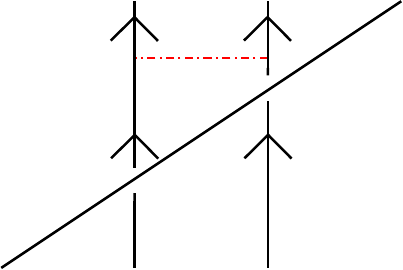}{.3} &= \att{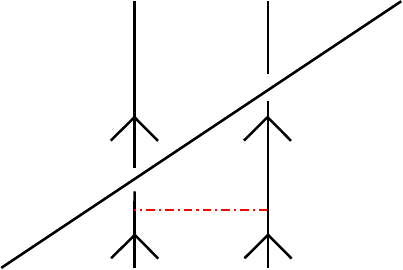}{.3} \qquad &&
   \text{(Exchange) } \att{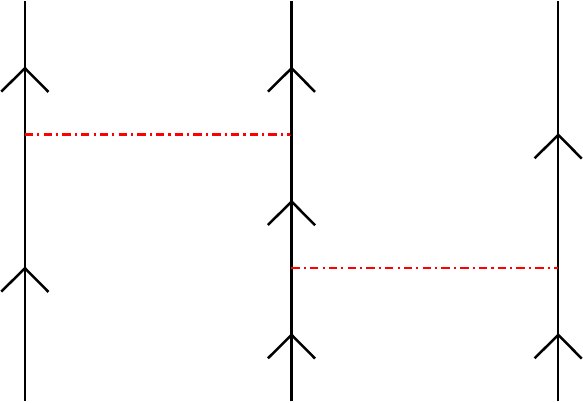}{.2} = -\att{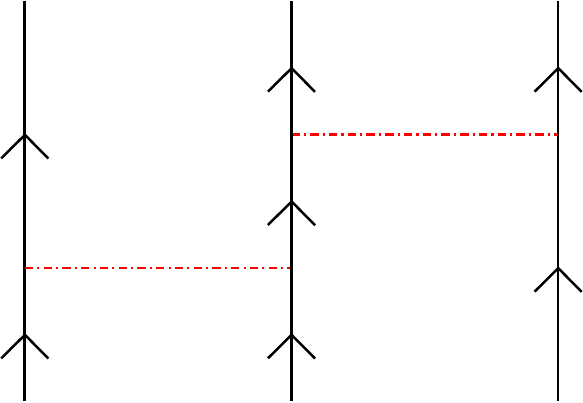}{.2}\\
    \text{(Stack) } \att{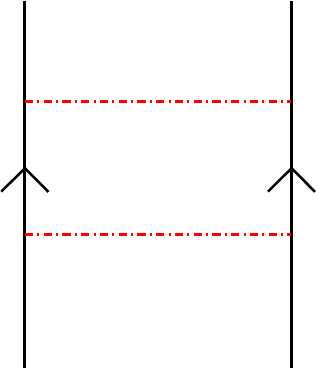}{.2} &= - \att{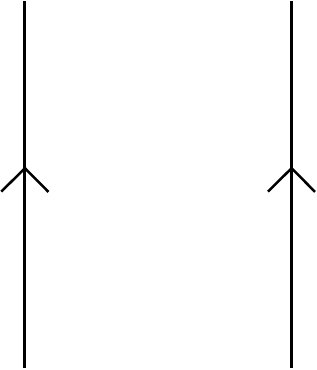}{.2}\qquad 
    && 
     \text{(Slide) } \att{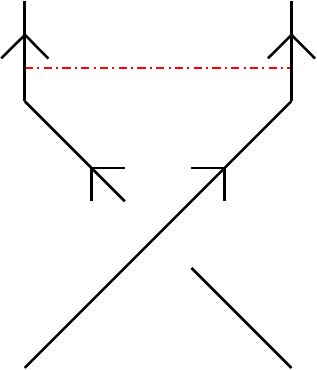}{.2} = -\att{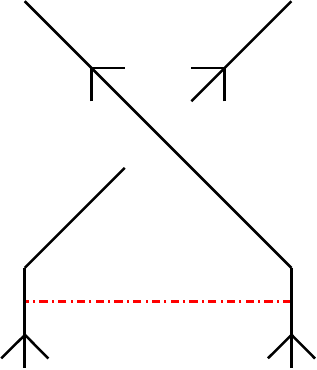}{.2} + (q-q^{-1})\att{splittingEndR5}{.25}\\
  \text{(Commute) } \att{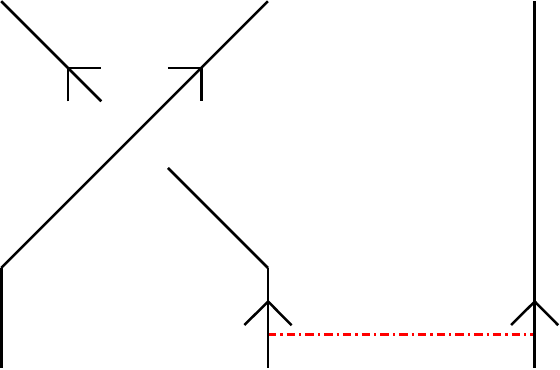}{.2}&= \att{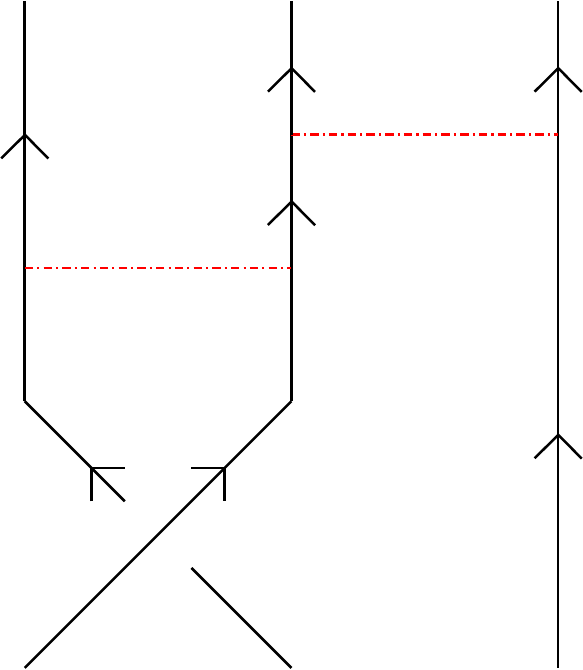}{.2}  .
\end{align*}    
\end{lem}
\begin{proof}
This is by direct computation. We include the derivation of (Commute), and leave the remainder to the reader. Using the defining relations of $\mathcal{E}_R$ we have
    \[ \att{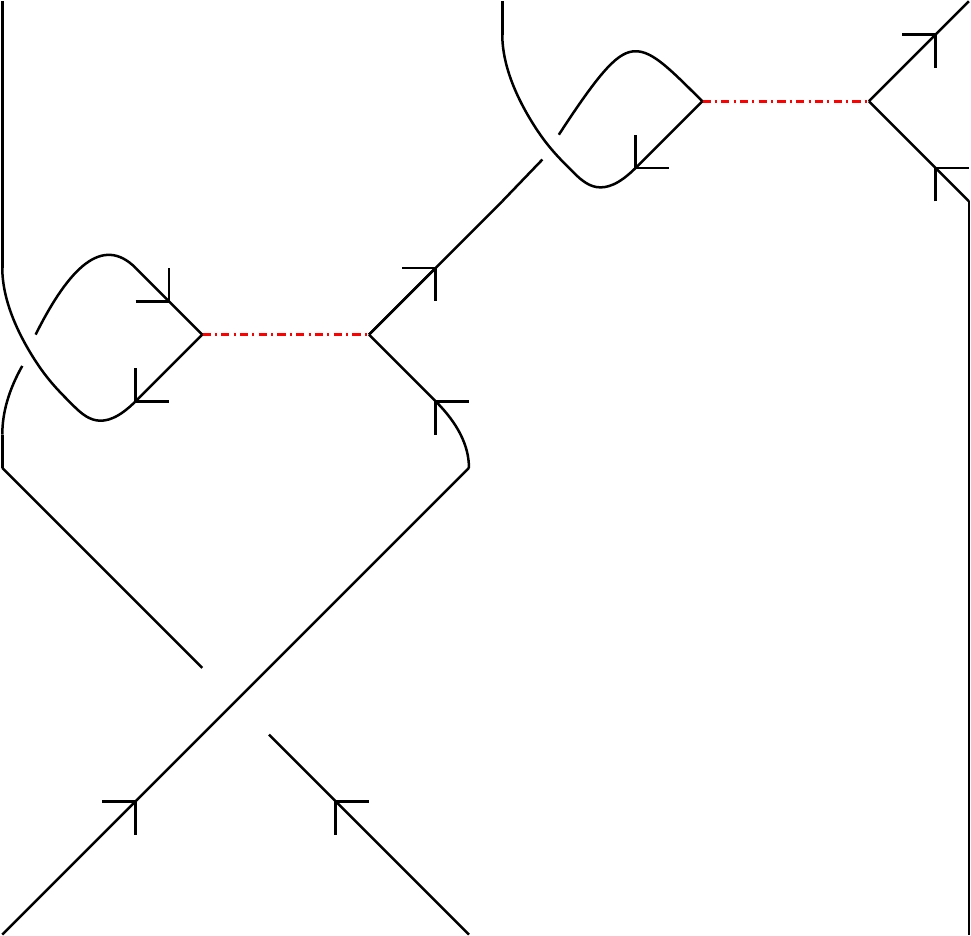}{.13}  =\att{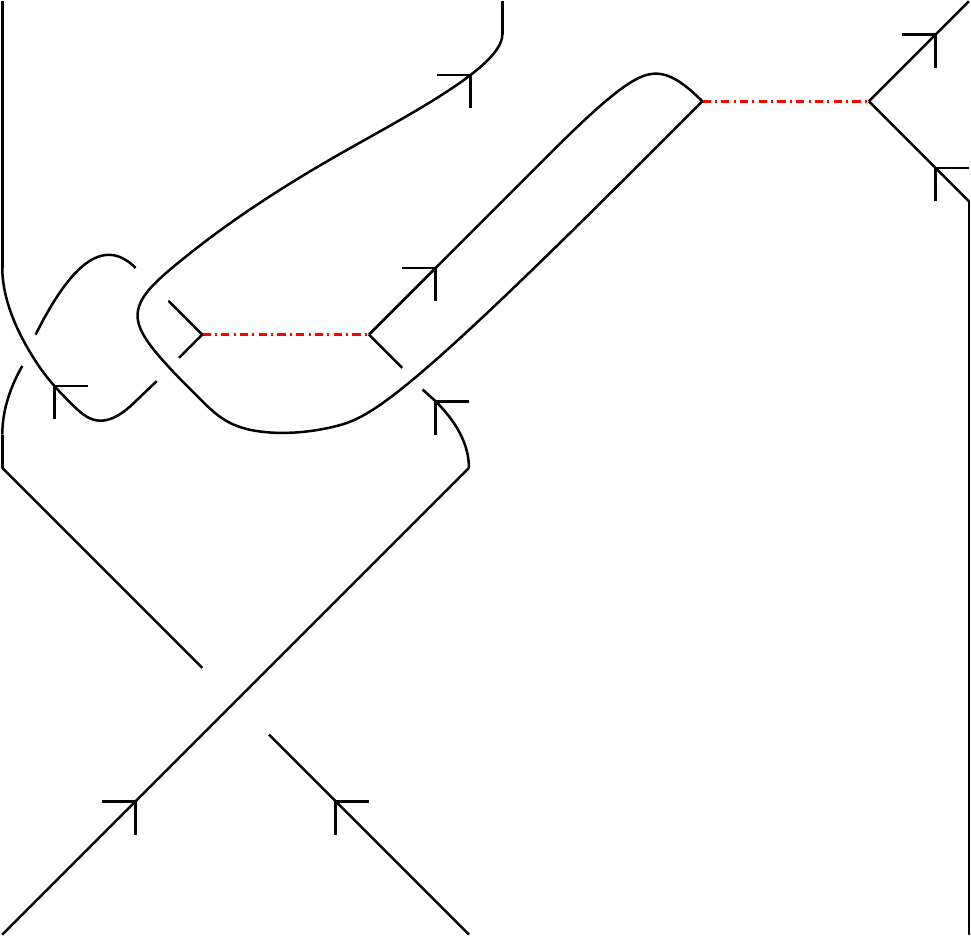}{.13} =\att{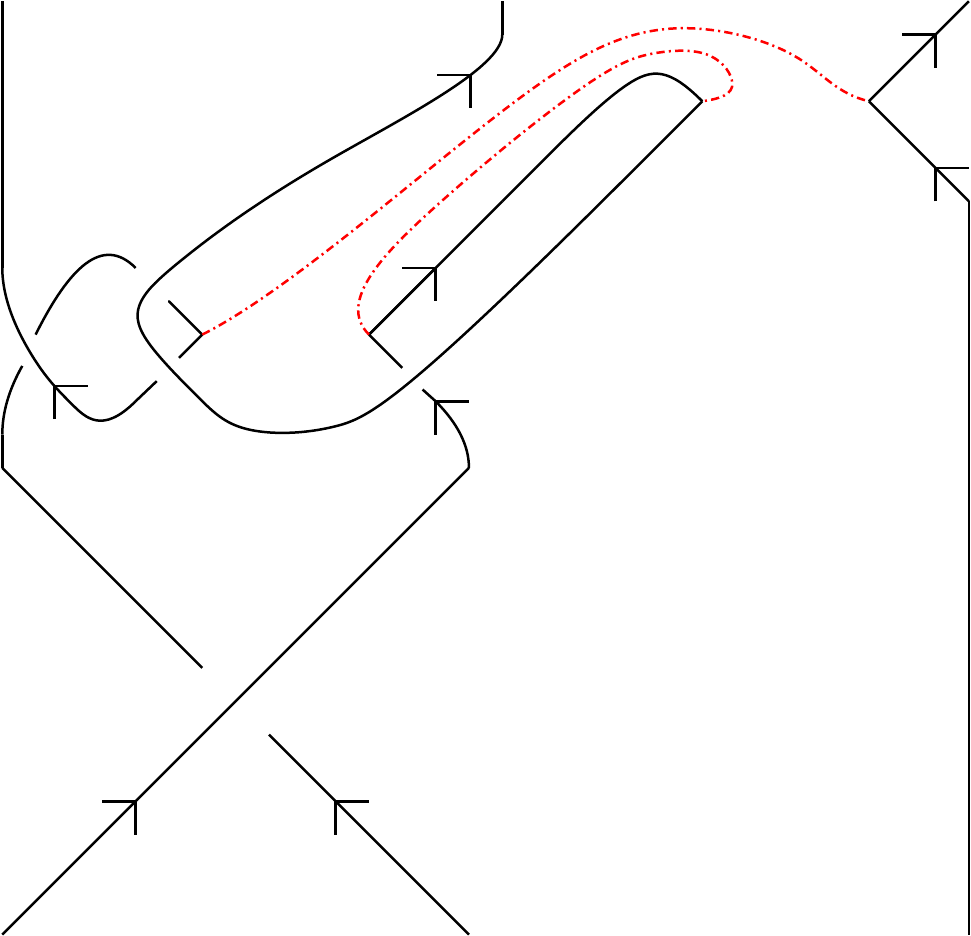}{.13} = \frac{q-q^{-1}}{2\mathbf{i}} \att{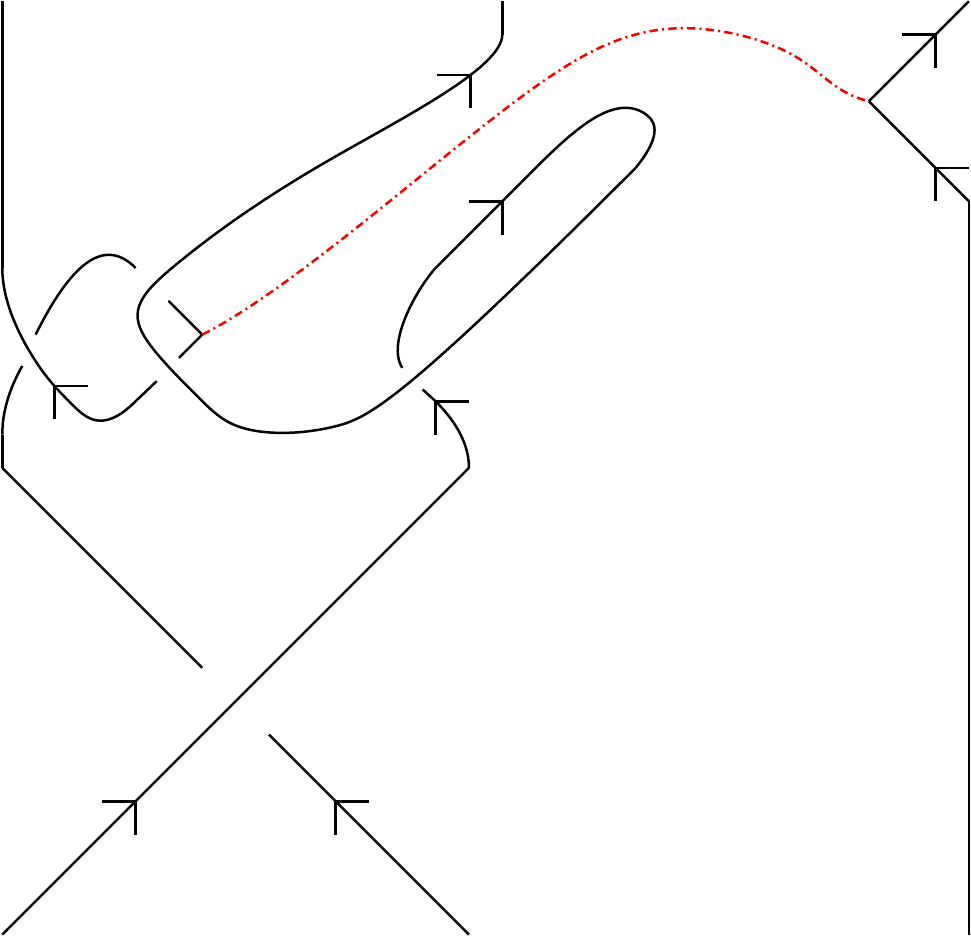}{.13} =-\frac{q-q^{-1}}{2} \att{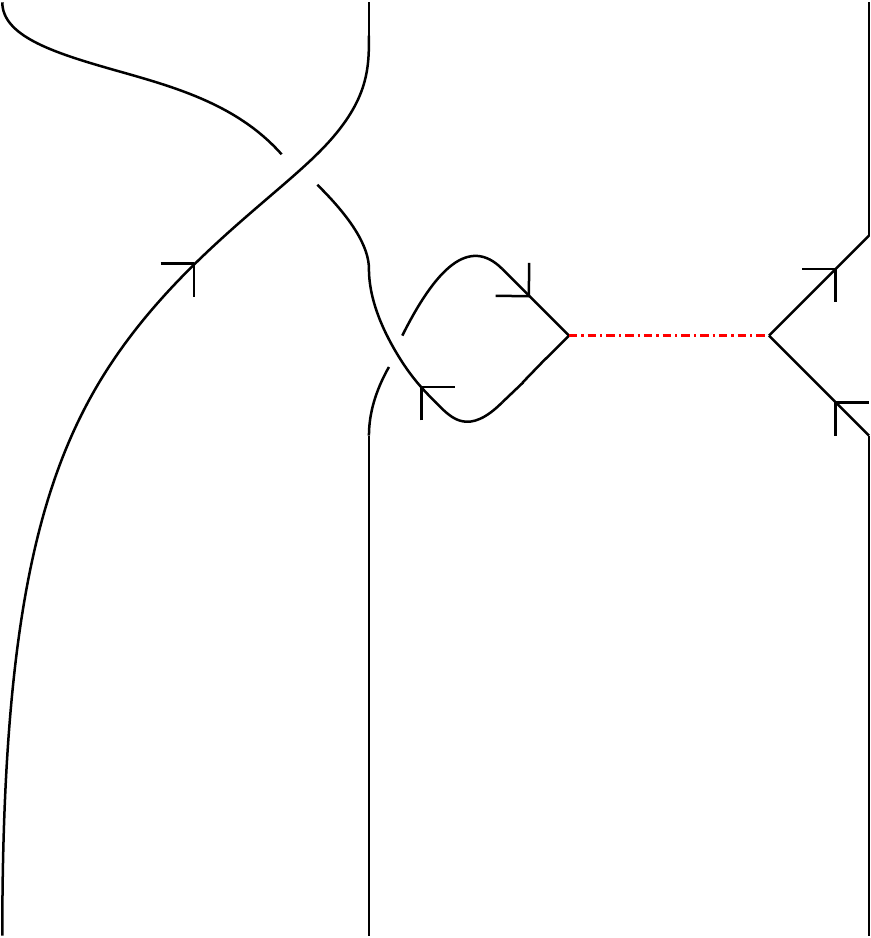}{.13}.  \]
    
Multiplying both sides by $\left(-\frac{2}{q-q^{-1}}\right)^2$ gives the desired relation.
\end{proof}

\begin{remark}\label{traceremark}
    It is worth pointing out that we could also define the category $\mathcal{E}_R$ as the rigid monoidal $R$-linear category generated by $\att{braidX}{.15}$ and $\att{splittingEnd2}{.15}$ along with the HOMFLYPT skein category relations, the five relations above, and the trace relations
    \[\att{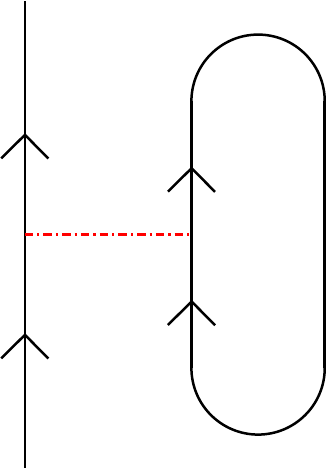}{.2} \quad=\quad 0 \quad=\quad \att{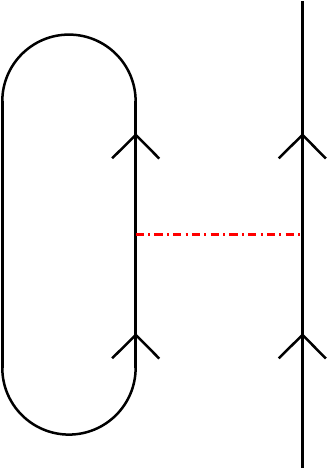}{.2} \qquad \text{and}\qquad \att{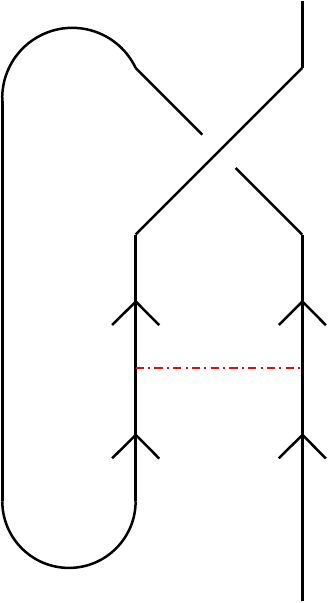}{.2} \quad=\quad \mathbf{i} \; \att{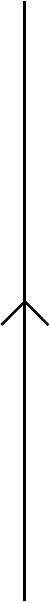}{.2}\quad  = \quad \att{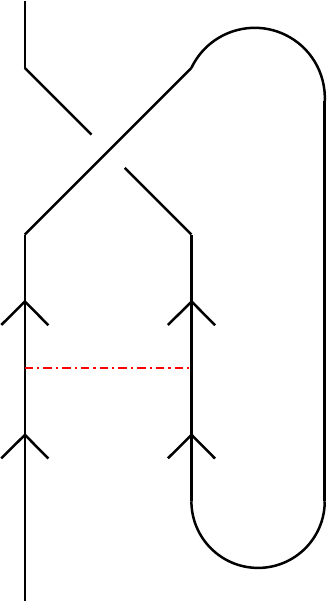}{.2}\]
\end{remark} 

\begin{cor}\label{cor:hom}
There's a ring homomorphism $\Psi: \Gno \rightarrow \operatorname{End}_{\mathcal{E}_R}(+^n)$ 
defined on generators by 
    \[        t_i\mapsto   \att{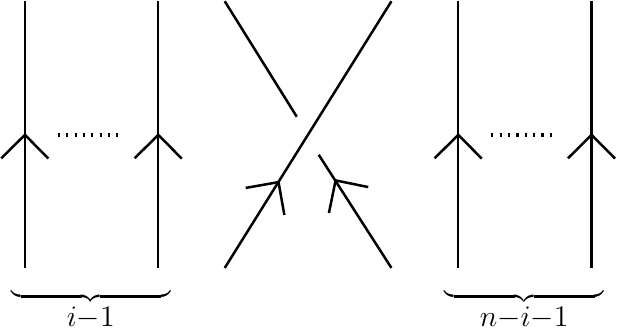}{.3}\qquad \text{ and }\qquad  e_i \mapsto  \att{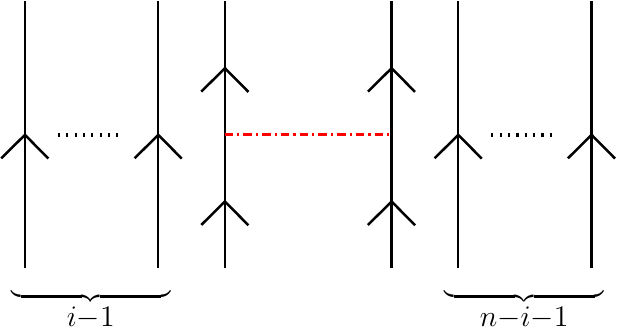}{.3}.        \]
\end{cor}

Recall that $\Gno$ has a basis consisting of elements of the form $h_w\cdot e_s$ where $h_w$ is the basis element of the Hecke algebra corresponding to $w\in S_n$ and $e_s$ is the basis of $\operatorname{Cliff}(n)$ indexed by $s\in \{0,1\}^{n-1}$. In a slight abuse of notation we will use the same labeling for the corresponding elements in $\mathcal{E}_R$ under the map $\Psi$.

Our goal is to construct a basis for the hom spaces of $\mathcal{E}_R$ in terms of \edit{ladders, braids, cups, and caps}. 
We begin with an intermediate result which shows that the algebras $\operatorname{End}_{\mathcal{E}_R}(+^n)$ are spanned by elements with do not contain any cups or caps. This result is analogous to the classical result that any oriented $(n,n)$ tangle can be written as a braid.

\begin{lem}\label{lem:straight}
    Let $D$ be a diagram with boundary $(+^n, +^n)$. Then $D$ can be expressed as an $R$-linear combination of diagrams that do not contain cups or caps.
\end{lem}
\begin{proof}
    We induct on the partial ordering of diagrams from Definition~\ref{def:houseordos}. The base cases consist of diagrams with no projection terms. This is exactly the classical result that any oriented $(n,n)$ tangle can be written as a braid.

    Let $D$ be a diagram with boundary $(+^n, +^n)$. By Proposition~\ref{prop:standardForm} we may write $D$ as a $R$-linear combination of diagrams in standard form, each of which is less than or equal to $D$ in the partial ordering. Hence we may assume that $D$ is in standard form.
    
    We select a \edit{cup or cap} in $D$. As $D$ is in standard form, this rigidity map must be on a black strand which either connects the bottom boundary to the top boundary, the bottom boundary to a projection term, or a projection term to the top boundary. In any case, we repeatedly use (Hecke) to express $D$ as a sum of the diagram $D$ with the black strand pulled to the top of the diagram, along with a $R$-linear combination of diagrams strictly smaller than $D$ in the partial ordering. By induction, all of these smaller diagrams can be expressed as a $R$-linear combination of diagrams which contain no cups or caps. For the diagram with the black strand on top, we can use the braid relations, along with (Over-Braid) and the zig-zag rigidity relation to remove all cups and caps on this strand. 

    Repeating this process for all black strands in $D$ gives an expression for $D$ as in the statement of the Lemma.
\end{proof}

Using this intermediate result, we can show that we
obtain the following spanning set in terms of the diagrams $e_s$ and $h_w$ for the endomorphism algebras in $\mathcal{E}_R$. 
It should be noted that these spanning diagrams are no longer in standard form.

\begin{lem}\label{lem:span}
    We have that
   \[ \operatorname{End}_{\mathcal{E}_R}(+^n) = \operatorname{span}_{R}\{h_w\cdot e_s :  w\in S_n,s\in \{0,1\}^{n-1} \}.\]
\end{lem}
\begin{proof}
Let $D$ be a diagram with boundary $(+^n, +^n)$. By Lemma~\ref{lem:straight} we may assume that $D$ contains no rigidity morphisms. Hence $D$ can be written as a word in the morphisms $h_w$ and $e_s$. The relations (Slide), (Over-Braid), and (Commute) allow us to write this as a word in $h_w$'s, composed with a word in $e_s$. It is well known that any word in $h_w$ can be reduced to an element of $\operatorname{span}_{\mathbb{C}}\{h_w :  w\in S_n \}$ using the HOMFLYPT skein category relations. It is immediate from relations (Exchange) and (Stack) that a word in the $e_s$ can be reduced to an element of $\operatorname{span}_{\mathbb{C}}\{ e_s: s\in \{0,1\}^{n-1} \}$.
\end{proof}

In general, showing linear independence of a set of morphisms in a category given by generators and relations is a difficult problem. In our setting we are fortunate to know from Corollary~\ref{cor:full} that when $q$ is specialised to $q = e^{2\pi i \frac{1}{4N}}$, the endomorphism algebras $\operatorname{End}_{\overline{\mathcal{E}_q}}(+^n)$ have the same dimension as the endomorphism algebras $\operatorname{End}_{\overline{\operatorname{Rep}(U_q(\mathfrak{sl}_N))}_A}(X^{\otimes n})$. The dimension of the later endomorphism algebras we have already computed in Theorem~\ref{thm:dimEnd}. This gives that our spanning set is in fact a basis when $q$ is specialised to certain roots of unity.
\begin{cor}\label{cor:basis}
    Let $N\in \mathbb{N}_{\geq 2}$ and set $q = e^{2\pi i \frac{1}{4N}}$. Then
    \[\{h_w\cdot e_s :  w\in S_n,s\in \{0,1\}^{n-1} \}\]
    is a basis for $\End_{\mathcal{E}_{q}}(+^n)$ for all $N >\edit{2n}$.
\end{cor}
\begin{proof}
    By Lemma~\ref{lem:span}, the set $\{h_w\cdot e_s :  w\in S_n,s\in \{0,1\}^{n-1} \}$ spans $\End_{\mathcal{E}_R(+^n)}$, and thus also spans $\End_{\mathcal{E}_q(+^n)}$. To see that this set forms \edit{forms} a basis we need only check that the number of elements in this set (namely, $2^{n-1} n!$) is no larger than the dimension of $\End_{\mathcal{E}_q(+^n)}$. Specifically,
    \[ \dim\End_{\mathcal{E}_{q}}(+^n)\geq \dim\End_{\overline{\mathcal{E}_{q}}}(+^n)=\dim\End_{\overline{\operatorname{Rep}(U_q(\mathfrak{sl}_N))}_A}(X^{\otimes n}) = 2^{n-1}\cdot n! \]
    Here the first equality follows from Corollary~\ref{cor:full}, and the second equality is Theorem~\ref{thm:dimEnd}.  
   
\end{proof}

Using the cups, caps and the half-braiding we can build an isomorphism between the endomorphism algebras of $\mathcal{E}_R$, and the other hom spaces. This isomorphism will allow us to obtain spanning sets of all hom spaces of $\mathcal{E}_R$, and \edit{bases} of small hom spaces in $\mathcal{E}_q$. 

\begin{defn}
    Suppose that $s_1$ and $s_2$ are strings of $\pm$ with $\sum s_1 = \sum s_2$. Let $n=\frac{|s_1|+|s_2|}{2}$ be the average of the lengths of the strings. We define \edit{a} vector space isomorphism $\beta_{s_1,s_2}: \operatorname{Hom}_{\mathcal{E}_q}(n\to n) \rightarrow \operatorname{Hom}_{\mathcal{E}_q}(s_1\to s_2)$ by using the half-braiding as in the following picture 
    \[\beta_{++-+-,-++}: \att{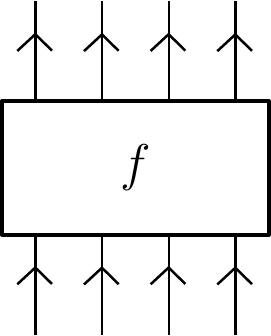}{.4}\qquad \mapsto \qquad  \att{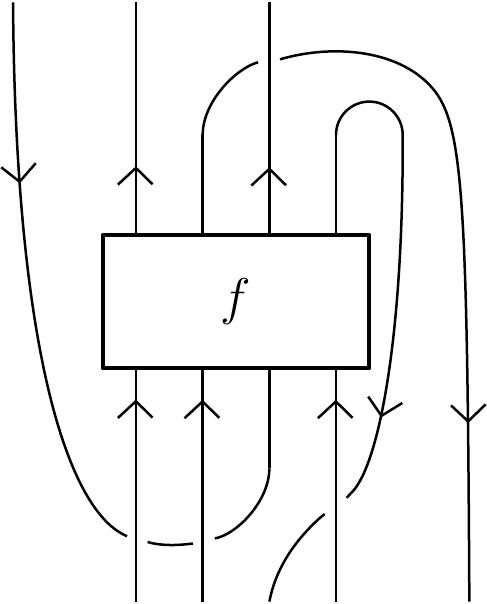}{.4}.\]
\end{defn}

\begin{defn}\label{def:iso}
    If $\sum s_1 \neq \sum s_2$ let $B_{s_1,s_2}$ denote the empty set. If $\sum s_1 = \sum s_2$ let $B_{s_1,s_2}$ denote the set $\{\beta_{s_1,s_2}(h_w\cdot e_s) :  w\in S_{\frac{|s_1|+|s_2|}{2}},s\in \{0,1\}^{\frac{|s_1|+|s_2|}{2}-1} \}$.
\end{defn}

Applying the isomorphisms $\beta_{s_1,s_2}$, we obtain the following two corollaries. 

\begin{cor}\label{cor:localbasis}
    For each $s_1\in \{+,-\}^n, s_2\in \{+,-\}^m$ we have that $B_{s_1,s_2}$ is a basis of $\operatorname{Hom}_{\mathcal{E}_q}(s_1\to s_2)$ for all $q = e^{2 \pi i \frac{1}{4N}}$ with $N > \frac{n+m}{4}$ 
\end{cor}
\begin{proof}
If $\sum s_1 \neq \sum s_2$ then a parity check shows that there are no diagrams with boundary $(s_1, s_2)$. Hence $\operatorname{Hom}_{\mathcal{E}_q}(s_1\to s_2)=0 $ in this case. If $\sum s_1 = \sum s_2$ then the result follows from applying $\beta_{s_1,s_2}$ to Corollary~\ref{cor:basis}
\end{proof}

\begin{cor}
    We have that $\operatorname{Hom}_{\mathcal{E}_R}(s_1\to s_2)$ is spanned over $R$ by $B_{s_1,s_2}$.
\end{cor}

A simple algebraic geometry style argument now shows that the sets $B_{s_1,s_2}$ form a basis for the hom spaces of $\mathcal{E}_R$.

\begin{thm}\label{thm:existence}
    For all objects $s_1, s_2$ we have that $B_{s_1,s_2}$ is an $R$-basis for $\operatorname{Hom}_{\mathcal{E}_R}(s_1\to s_2)$.
\end{thm}
\begin{proof}
    We already know that $B_{s_1,s_2}$ spans, we need only show that it's linearly independent. Suppose we have a linear dependence $\sum c_i b_i = 0$ for $c_i \in R$ and $b_i \in B_{s_1,s_2}$. By Corollary \ref{cor:localbasis} we have that the rational function $c_i$ is zero after specializing to $q = e^{2 \pi i \frac{1}{4N}}$ with $N > \frac{|s_1|+|s_2|}{4}$. But a rational function with infinitely many zeros is identically zero, so $c_i = 0$, and thus $B_{s_1,s_2}$ is linearly independent.
\end{proof}

In particular, this result shows that the endomorphism algebras in $\mathcal{E}_q$ are Hecke-Clifford algebras.

\begin{cor}\label{cor:gniso}
    For all $q\in \mathbb{C}-\{-1,0,1\}$, the map 
    \[\Psi:\Gno \to \operatorname{End}_{\mathcal{E}_q}(+^n)\]
    from Corollary~\ref{cor:hom} is an isomorphism.
\end{cor}

\bibliography{Cells}
\bibliographystyle{alpha}
\end{document}